\documentclass[11pt,letterpaper,reqno]{amsart}
\usepackage{amsmath,amsthm,amsfonts,amssymb,mathtools,algpseudocode}
\usepackage{comment,bm}
\usepackage{multicol}

\usepackage{bookmark,hyperref}
\hypersetup{pdfstartview={FitH}}

\newtheorem{Theorem}{{\bf Theorem}}[section]
\newtheorem{Example}[Theorem]{{\bf Example}}
\newtheorem{Corollary}[Theorem]{{\bf Corollary}}
\newtheorem{Algorithm}[Theorem]{{\bf Algorithm}}
\newtheorem{Proposition}[Theorem]{{\bf Proposition}}
\newtheorem{Definition}[Theorem]{{\bf Definition}}
\newtheorem{Lemma}[Theorem]{{\bf Lemma}}
\newtheorem{Remark}[Theorem]{{\bf Remark}}

\numberwithin{equation}{section}

\newcommand{\Ol}{\mbox{\Large $\mathcal O$}}
\newcommand{\Os}{\mbox{\footnotesize $\mathcal O$}}
\newcommand{\Il}{\mbox{$\mathcal R$}}
\newcommand{\Ilb}{\mbox{\Large $\mathcal R$}}
\newcommand{\Jl}{\mbox{\Large $\mathcal J$}}

\newcommand{\ubce}{\text{UBCE}}

\newcommand{\ve}{\text{vec}}

\def\restrict#1{\raise-.5ex\hbox{\ensuremath|}_{#1}}

\begin{document}

\title[Convergence of the Cyclic and Quasi-cyclic Block Jacobi Methods]{Convergence of the Cyclic and Quasi-cyclic Block Jacobi Methods}

\author{Vjeran Hari}\thanks{Vjeran Hari, Department of Mathematics, Faculty of Science, University of Zagreb, Bijeni\v{c}ka 30, 10000 Zagreb, Croatia}
\author{Erna Begovi\'{c}~Kova\v{c}}\thanks{Erna Begovi\'{c}~Kova\v{c}, Faculty of Chemical Engineering and Technology, University of Zagreb, Maruli\'{c}ev trg 19, 10000 Zagreb, Croatia}
\thanks{This work has been fully supported by Croatian Science Foundation under the project 3670.}
\date{19 February 2017}


\begin{abstract}
The paper studies the global convergence of the block Jacobi me\-thod for symmetric matrices.
Given a symmetric matrix $A$ of order $n$, the method generates a sequence of matrices by the rule $A^{(k+1)}=U_k^TA^{(k)}U_k$, $k\geq0$, where $U_k$ are orthogonal elementary block matrices. A class of generalized serial pivot strategies is introduced, significantly enlarging the known class of weak wavefront strategies, and appropriate global convergence proofs are obtained. The results are phrased in the stronger form: $S(A')\leq c S(A)$, where $A'$ is the matrix obtained from $A$ after one full cycle, $c<1$ is a constant and $S(A)$ is the off-norm of $A$.  Hence, using the theory of block Jacobi operators, one can apply the obtained results to prove convergence of block Jacobi methods for other eigenvalue problems, such as the  generalized eigenvalue problem. As an example, the results are applied to the block $J$-Jacobi method. Finally, all results are extended to the corresponding quasi-cyclic strategies.
\end{abstract}

\maketitle

\section{Introduction}

The main incentive for writing this paper was a need to expand the class of ``convergent strategies'' for the block Jacobi method for symmetric matrices. With a large choice of classes at our disposal we can prove global convergence of other block-wise or element-wise Jacobi-type methods and even apply it to related problems, e.g., to the generalized eigenvalue or singular value problem (see~\cite{nov+sin-15}). The techniques we are about to employ use the theory of block Jacobi operators, which was described in~\cite{har-15}.

Over the last two decades the Jacobi method has emerged as a method of choice for the eigenvalue computation for dense symmetric matrices. This is mostly due to its inherent parallelism and high relative accuracy on well-behaved matrices. Although the original method is very old~\cite{jac-1846} and it had been one of the first methods to be implemented on computers, it was forgotten in the 1970s after appearance of the QR and Divide and conquer method. Already in 1971 Sameh~\cite{sam-71} showed how to adapt the serial Jacobi method to parallel processing. Later, in 1992, Demmel and Veseli\'{c}~\cite{dem-ves-92} proved high relative accuracy of the method on well-behaved symmetric positive definite matrices. Following their breakthrough, the method came back to the focus of the current researcher. Drma\v{c} and Veseli\'{c}~\cite{drm+ves-04a,drm+ves-04b} showed that, even on standard one-processor computers, the method can be modified to become faster than the QR method while still retaining its distinguished property: high relative accuracy. Nowadays, the Jacobi method is well understood. On the one hand, its asymptotic convergence was considered in~\cite{wil-62,har-91,rhe+har-93} and its global convergence was studied in~\cite{for+hen-60,han-63,hen+zim-68,naz-75,har+ves-87,fer-89,luk+par-89,shr+sch-89,mas-90,har-07}. On the other hand, its high relative accuracy was considered in~\cite{ves+har-89,dem-ves-92,math-95,mat-08}, while its efficiency was investigated in~\cite{drm+ves-04a,drm+ves-04b}. The method has also been implemented as a standard LAPACK routine.

With the development of CPU and GPU parallel computing platforms, it has been found that a sensible way of increasing numerical efficiency of the method involves using a one-sided algorithm, together with BLAS~3 subroutines, which can nicely exploit cache memory hierarchy. The matrix description of such a method is called the block Jacobi method. This block method is always implemented as a one-sided block (Jacobi or J-Jacobi) algorithm because high efficiency and high relative accuracy are warranted then. However, in the global and asymptotic convergence considerations the results are  cast in terms of a two-sided block Jacobi method. The first global convergence results for the block Jacobi methods were given
in~\cite{drm-07,buj+drm-12,har+sin-10,har+sin-11,har-15}. These papers considered the most common serial pivot strategies and the strategies equivalent to them.

The aim of this paper is to further develop the global convergence theory for the block Jacobi method and to provide a large class of usable pivot strategies for which the convergence can be established rigorously. In general, our class consists of more than $4\cdot 2!\cdot 3!\cdots m!$ cyclic strategies, where $m$ is the number of block-columns in the block-matrix partition of a symmetric matrix of order $n$. These strategies include the weak-wavefront ones from~\cite{shr+sch-89} and many others. As a byproduct of this research, we can now prove that every cyclic (element-wise or block) Jacobi method for symmetric matrices of order $4$ is globally convergent (see~\cite{B,beg+har15a}). In addition, we consider a similar class of quasi-cyclic strategies and derive the corresponding convergence results. The block analogue of the strategy that is used in the LAPACK implementation of the Jacobi method lies in that class.

The convergence results are given in the ``stronger form'',
$$S(A')\leq c S(A), \quad 0\leq c <1.$$
Here, $A$ is the initial symmetric matrix of order $n$, $A'$ is obtained from $A$ after applying one sweep of some cyclic or quasi-cyclic block Jacobi method, $S(\ )$ is departure from the diagonal form, and $c$ is a constant depending on $n$ and the block-matrix partition, but not on $A$. Such a result allows for the use of the theory of block Jacobi operators. Hence, it can be utilized to prove the global convergence of other Jacobi-type methods, designed for different eigenvalue problems. As an application, we will apply it to the block J-Jacobi method from~\cite{har+sin-11}.

Some of the results presented here can be found in the unpublished thesis~\cite{B}.

The paper is divided into six sections of the main text and an appendix. In Section~2 we present the basic concepts linked with a block Jacobi method for symmetric matrices. Special attention is paid to cyclic and quasi-cyclic pivot strategies, and to the ways of enlarging significantly the number of ``convergent strategies''. The concepts of equivalent, weak-equivalent and permutation equivalent strategies are used. Another useful tool is introduced, the so-called block Jacobi annihilators and operators for symmetric matrices, and some basic results related to them are proved. In Section~3 we introduce a class of generalized serial strategies and prove the corresponding global convergence results. In Section~4 we briefly introduce a similar class of quasi-cyclic pivot strategies and prove the appropriate convergence results. As an application, in Section~5 we prove the global convergence of the block J-Jacobi method under the strategies from the newly introduced classes. Section~6 announces the future work.
Finally, to make the paper easier to read, we move all lengthy and technical proofs to the Appendix.

\section{Basic concepts and notation}

We introduce the basic definitions linked with the block Jacobi method for symmetric matrices. Special attention is paid to the cyclic and quasi-cyclic pivot strategies. Later we deal with more advanced concepts like the block Jacobi annihilators and operators.

\subsection{Block Jacobi method}

Let $A$ be a square matrix of order $n$ and let $\pi$ be an integer partition of $n$,
\begin{equation}\label{pi}
\pi=(n_1,n_2,\ldots,n_m), \quad n_i\geq1, \ 1\leq i\leq m; \quad n_1+n_2+\cdots+n_m=n.
\end{equation}
Then $\pi$ determines the block-matrix partition of $A$,
\begin{equation}\label{blokmatrica}
A=\left[\begin{array}{cccc}
      A_{11} & A_{12} & \ldots & A_{1m} \\
      A_{21} & A_{22} &  & A_{2m} \\
      \vdots &  & \ddots & \vdots \\
      A_{m1} & A_{m2} & \ldots & A_{mm} \\
    \end{array}
  \right] \begin{array}{c}
             n_1 \\
             n_2 \\
             \vdots \\
             n_m \\
           \end{array},
\end{equation}
where the diagonal blocks $A_{11},\ldots,A_{mm}$ are square matrices of order
$n_1,\ldots,n_m$, respectively. Relation~\eqref{blokmatrica} will be schematically denoted by $A=(A_{rs})$.

Since we consider the global convergence of the block Jacobi method for symmetric matrices, we assume that $A$ is symmetric. A block Jacobi method is determined by the partition $\pi$, some pivot strategy and the algorithm. The partition is chosen in accordance with the capacity of the hierarchical cache memory of the computer.
Typically, the code presumes $n_1=n_2=\cdots =n_{m-1}$, $n_m=n-(m-1)n_1\leq n_1$.
In our analysis we consider it arbitrary but unchanged over the iterations. Actually, there are situations when it makes sense to change $\pi$ during the process, but these are linked with the asymptotic convergence of the method.

The block Jacobi method uses orthogonal \emph{elementary block matrices} as transformation matrices. An orthogonal elementary block matrix $\mathbf{U}_{ij}$ has the form (see~\cite{har-15})
\begin{equation}\label{Uij}
\mathbf{U}_{ij} = \left[\begin{array}{ccccc}
             I &  &  &  &  \\
              & U_{ii} &  & U_{ij} &  \\
              &  & I &  &  \\
              & U_{ji} &  & U_{jj} &  \\
              &  &  &  & I
           \end{array}\right]
\begin{array}{c}  \\ n_i \\  \\ n_j \\ \\ \end{array}
\text{if} \ i<j, \ \text{or} \
\mathbf{U}_{ij} =
\left[\begin{array}{ccc} I &  & \\ & U_{ii} &  \\ &  & I \end{array}\right]
\begin{array}{c} \\ n_i \\ \\ \end{array}
\ \text{if} \ i=j,
\end{equation}
where it is presumed that the block-matrix partition is determined by $\pi$ from the relation~\eqref{pi}. Since $i$ and $j$ address the blocks, they can be called \emph{block pivot indices}, but for brevity we simply call them \emph{pivot indices}. Similarly, $(i,j)$ is the \emph{pivot pair} and
\begin{equation}\label{piv_mat_E}
\widehat{U}_{ij}=\left[\begin{array}{cc} U_{ii} & U_{ij}\\ U_{ji} & U_{jj}\end{array}\right] \ \text{if} \ i<j, \qquad \text{or} \ \widehat{U}_{ij}=U_{ii} \ \text{if} \  i=j,
\end{equation}
is the \emph{pivot submatrix} of $\mathbf{U}_{ij}$. When $(i,j)$ is understood we will also write
$\widehat{U}$ instead of $\widehat{U}_{ij}$.
We can build an orthogonal elementary block matrix using the function $\mathcal{E}$ which imbeds any
orthogonal matrix $\widetilde{U}$ of order $n_i+n_j$ (or $n_i$ if $i=j$)  into the identity matrix $I_n$, so that $\mathbf{U}_{ij}=\mathcal{E}(i,j,\widetilde{U})$ implies $\widehat{U}_{ij}=\widetilde{U}$. The mapping $\mathcal{E}$ depends on the partition $\pi$.

Each \textit{block Jacobi method} is an iterative processes of the form
\begin{equation}\label{jacobiblokagm}
A^{(k+1)}=U_k^TA^{(k)}U_k, \quad k\geq0; \quad A^{(0)}=A,
\end{equation}
where $U_k$, $k\geq0$, are orthogonal elementary block matrices.
Let $A^{(k)}=(A_{rs}^{(k)})$. We say that $A^{(k+1)}$ is obtained or generated from $A^{(k)}$ at \emph{step} $k$ via the recursion~\eqref{jacobiblokagm}. Let  $U_k=\mathcal{E}(i(k),j(k),\widehat{U}_k)$. Then $i(k)$, $j(k)$ are the pivot indices and $(i(k),j(k))$ is the pivot pair at step $k$. For brevity, we will often omit $k$ and denote the pivot indices simply by $i$, $j$ and the pivot pair by $(i,j)$. The way of selecting the pivot pair at each step is called a \emph{pivot strategy}.

At step $k$ the block Jacobi method diagonalizes the pivot submatrix of
$A^{(k)}$. Thus, if $i<j$, the pivot blocks $A_{ij}^{(k)}$ and $A_{ji}^{(k)}$ are annihilated and the affected diagonal blocks $A_{ii}^{(k)}$ and $A_{jj}^{(k)}$ are diagonalized. If $\widehat{A}^{(k)}$ denotes the pivot submatrix of order $n_i+n_j$ at step $k$, it is transformed as follows:
\begin{equation}\label{blokpivotmatr}
\left[\begin{array}{cc}
    \Lambda_{ii}^{(k+1)} & 0 \\
    0 & \Lambda_{jj}^{(k+1)} \\
\end{array}
\right]=\left[
\begin{array}{cc}
    U_{ii}^{(k)} & U_{ij}^{(k)} \\
    U_{ji}^{(k)} & U_{jj}^{(k)} \\
\end{array}
\right]^T\left[
\begin{array}{cc}
    A_{ii}^{(k)} & A_{ij}^{(k)} \\
    (A_{ij}^{(k)})^T & A_{jj}^{(k)} \\
\end{array}
\right]\left[
\begin{array}{cc}
    U_{ii}^{(k)} & U_{ij}^{(k)} \\
    U_{ji}^{(k)} & U_{jj}^{(k)} \\
\end{array}\right],
\end{equation}
where $\Lambda_{ii}^{(k+1)}$ and $\Lambda_{jj}^{(k+1)}$ are diagonal. If $i=j$, then just
$A_{ii}^{(k)}$ is diagonalized. For the diagonalization of the pivot submatrix, one can choose any method for solving the symmetric eigenvalue problem. Typically, one applies a standard (element-wise) Jacobi method for its high relative accuracy~\cite{mat-08} and efficiency on nearly diagonal matrices.

As has been explained in~\cite{har-15}, it is preferable to preprocess the initial matrix by $m$ block Jacobi steps with pivot pairs $(1,1),\ldots,(m,m)$, so that in the starting matrix the diagonal blocks are actually diagonal submatrices. This preprocessing is depicted below for the case $\pi=(3,2,1,2)$:
$$A =
\mbox{\scriptsize $\displaystyle \left[\begin{array}{ccc|cc|c|cc}
        x & x & x & x & x & x & x & x   \\
        x & x & x & x & x & x & x & x   \\
        x & x & x & x & x & x & x & x   \\ \hline
        x & x & x & x & x & x & x & x   \\
        x & x & x & x & x & x & x & x   \\ \hline
        x & x & x & x & x & x & x & x   \\ \hline
        x & x & x & x & x & x & x & x   \\
        x & x & x & x & x & x & x & x
      \end{array} \right]$} \quad\longmapsto\quad
A^{(0)}  =
\mbox{\scriptsize $\displaystyle \left[\begin{array}{ccc|cc|c|ccc}
        x & 0 & 0 & x & x & x & x & x  \\
        0 & x & 0 & x & x & x & x & x  \\
        0 & 0 & x & x & x & x & x & x  \\ \hline
        x & x & x & x & 0 & x & x & x  \\
        x & x & x & 0 & x & x & x & x  \\ \hline
        x & x & x & x & x & x & x & x  \\ \hline
        x & x & x & x & x & x & x & 0  \\
        x & x & x & x & x & x & 0 & x
      \end{array}
    \right]$}. \label{A0}
$$
Once the diagonal blocks are diagonalized, all later steps will preserve that property. It means that
at each step the pivot indices will satisfy $i<j$, which unifies and simplifies the algorithm. In this regard the blocks $A_{ii}$ and $A_{jj}$ on the right side of the relation~\eqref{blokpivotmatr} can be replaced by $\Lambda_{ii}$ and $\Lambda_{jj}$, respectively.
Therefore, in the sequel it is presumed that the diagonal blocks of each $A^{(k)}$  are diagonal and for the pivot indices $i<j$ holds.

\subsection{Pivot strategies}\label{sec:2.2}

Each pivot strategy can be identified with a function
$I:\mathbb{N}_0\rightarrow \mathcal{P}_m,$
where $\mathbb{N}_0=\{0,1,2,3,\ldots\}$ and $\mathcal{P}_m=\{(r, s) | 1\leq r<s\leq m\}$.
If $I$ is a periodic function with period $T$, we say that $I$ is a \textit{periodic pivot strategy}. In this paper we consider two types of periodic strategies: cyclic and quasi-cyclic ones.

If $T=M\equiv\frac{m(m-1)}{2}$ and $\{(i(0),j(0)),(i(1),j(1)),\ldots,(i(T-1),j(T-1))\}=\mathcal{P}_m$, then we say that
the pivot strategy is \emph{cyclic}. It immediately follows that, during any $M$ successive steps of the method, all off-diagonal blocks are annihilated exactly once. Such block Jacobi method is also said to be cyclic and the transition from $A^{((r-1)M)}$ to $A^{(rM)}$ is called the $r$th \emph{cycle} or \emph{sweep} of the method.

If $T\geq M$ and $\{(i(0),j(0)),(i(1),j(1)),\ldots,(i(T-1),j(T-1))\}=\mathcal{P}_m$, then the strategy is called \emph{quasi-cyclic}. Thus, during any $T$ successive steps of the method, each off-diagonal block is annihilated at least once. The corresponding block Jacobi method is called \emph{quasi-cyclic} and the transition from $A^{((r-1)T)}$ to $A^{(rT)}$ is called the $r$th \emph{quasi-cycle} or \emph{sweep} of the method.

Let us examine cyclic and quasi-cyclic strategies more closely.

For $\mathcal{S}\subseteq \mathcal{P}_m$ we denote by $\mathcal{\Ol}(\mathcal{S})$ the set of all finite sequences containing the elements of $\mathcal{S}$, assuming that each pair from $\mathcal{S}$ appears at least once in each sequence.
If $I$ is a cyclic or quasi-cyclic strategy with period $T$, then $\mathcal{O}_I$ stands for the sequence $I(0),I(1),\ldots,I(T-1)\in\mathcal{\Ol}(\mathcal{P}_m)$, generated by the first $T$ steps (i.e., by the first \emph{sweep}) of the method.
Conversely, if $\mathcal{O}\in\mathcal{\Ol}(\mathcal{P}_m)$, $\mathcal{O}=(i_0,j_0),(i_1,j_1),\ldots,(i_{T-1},j_{T-1})$, then the periodic stra\-te\-gy $I_{\mathcal{O}}$ is defined by $I_{\mathcal{O}}(k)=(i_{\tau(k)},j_{\tau(k)})$, where $\tau(k)$ is the unique integer satisfying $0\leq \tau(k)\leq T-1$ and $k\equiv$ \linebreak $\tau(k)(\!\!\! \mod \ T)$, $k\geq0$.

These two functions, $\mathcal{O}\mapsto I_{\mathcal{O}}$ and $I\mapsto \mathcal{O}_I$, enable us to investigate the cyclic and quasi-cyclic strategies by studying the sequences from $\mathcal{\Ol}(\mathcal{P}_m)$. Note that, if $I$ is cyclic, then $\mathcal{O}_I$ is simply an \emph{ordering} of $\mathcal{P}_m$. We will also use the term \emph{pivot ordering} in this case, while if $I$ is quasi-cyclic, we will use the term \emph{pivot sequence}.

An \emph{admissible transposition} on $\mathcal{O}\in\mathcal{\Ol}(\mathcal{S})$, $\mathcal{S}\subseteq\mathcal{P}_m$, is any transposition of two adjacent terms in $\mathcal{O}$,
\[
(i_r,j_r),(i_{r+1},j_{r+1})\rightarrow(i_{r+1},j_{r+1}),(i_r,j_r),
\]
provided that the sets $\{i_r,j_r\}$ and $\{i_{r+1},j_{r+1}\}$ are disjoint. We also say that such pairs $(i_r,j_r)$ and $(i_{r+1},j_{r+1})$ \emph{commute}. The number of pairs in $\mathcal{O}$ is denoted by $|\mathcal{O}|$ and it is called the \emph{length} of $\mathcal{O}$.

\begin{Definition}
Two sequences $\mathcal{O},\mathcal{O}'\in\mathcal{\Ol}(\mathcal{S})$, $\mathcal{S}\subseteq \mathcal{P}_m$, are said to be
\begin{itemize}
\item[(i)] \emph{equivalent} (we write $\mathcal{O}\sim\mathcal{O}'$) if one can be obtained from the other by a finite set of admissible transpositions;
\item[(ii)] \emph{shift-equivalent} ($\mathcal{O}\stackrel{\mathsf{s}}{\sim}\mathcal{O}'$) if $\mathcal{O}=[\mathcal{O}_1,\mathcal{O}_2]$ and $\mathcal{O}'=[\mathcal{O}_2,\mathcal{O}_1]$, where $[ \ , \  ]$ stands for concatenation and the length of $\mathcal{O}_1$ is called shift length;
\item[(iii)] \emph{weak equivalent} ($\mathcal{O}\stackrel{\mathsf{w}}{\sim}\mathcal{O}'$) if there exist $\mathcal{O}_i\in\mathcal{\Ol}(\mathcal{S})$, $0\leq i\leq r$, such that every two adjacent terms in the sequence $\mathcal{O}=\mathcal{O}_0,\mathcal{O}_1,\ldots,\mathcal{O}_r=\mathcal{O}'$ are equivalent or shift-equivalent.
\end{itemize}
\end{Definition}

One can check that $\sim$, $\stackrel{\mathsf{s}}{\sim}$ and $\stackrel{\mathsf{w}}{\sim}$ are equivalence relations on $\mathcal{\Ol}(\mathcal{S})$. If three or more sequences are connected by $\sim$ or $\stackrel{\mathsf{s}}{\sim}$ one can omit the mid terms because of the transitivity property of equivalence relation. Hence, if $\mathcal{O}\stackrel{\mathsf{w}}{\sim}\mathcal{O}'$, then there is a sequence $\mathcal{O}=\mathcal{O}_0,\mathcal{O}_1,\ldots,\mathcal{O}_r=\mathcal{O}'$
such that
\begin{equation}\label{weak_eq}
\text{either} \ \mathcal{O}_0\sim
\mathcal{O}_1\stackrel{\mathsf{s}}{\sim}\mathcal{O}_2\sim\mathcal{O}_3\stackrel{\mathsf{s}}{\sim}\mathcal{O}_4\ldots \qquad \text{or} \ \mathcal{O}_0\stackrel{\mathsf{s}}{\sim}
\mathcal{O}_1\sim\mathcal{O}_2\stackrel{\mathsf{s}}{\sim}\mathcal{O}_3\sim\mathcal{O}_4\ldots.
\end{equation}
Two pivot strategies $I_{\mathcal{O}}$ and $I_{\mathcal{O}'}$ are equivalent (resp.\@ shift-equivalent, weak equivalent) if the corresponding sequences $\mathcal{O}$ and $\mathcal{O}'$ are equivalent (resp.\@ shift-equivalent, weak equivalent).

The most common cyclic pivot strategies are the row-cyclic one, $I_{\text{row}}=I_{\mathcal{O}_{\text{row}}}$, and the column-cyclic one, $I_{\text{col}}=I_{\mathcal{O}_{\text{col}}}$, which are defined by the ``row-wise'' and ``column-wise''  orderings of $\mathcal{P}_m$:
\begin{align*}
\mathcal{O}_{\text{row}} & = (1,2),(1,3),\ldots,(1,m), (2,3),\ldots,(2,m),\ldots,(m-1,m) \qquad \text{and} \\
\mathcal{O}_{\text{col}} & = (1,2),(1,3),(2,3),\ldots,(1,m), (2,m),\ldots,\ldots,(m-1,m).
\end{align*}
The common name for them is \emph{serial strategies}. The cyclic pivot strategies which are equivalent (resp.\@ weak-equivalent) to the serial ones are also called wavefront (resp.\@ weak-wavefront) strategies (see~\cite{shr+sch-89}).

\begin{Definition}\label{tm: def_reverse}
Let $\mathcal{O}\in\mathcal{\Ol}(\mathcal{S})$, $\mathcal{S}\subseteq \mathcal{P}_m$, $\mathcal{O}=(i_0,j_0),(i_1,j_1),\ldots,(i_{r},j_{r})$. Then
\[
\mathcal{O}^{\leftarrow}=(i_{r},j_{r}),\ldots,(i_1,j_1),(i_0,j_0) \ \in\mathcal{\Ol}(S)
\]
is the reverse (or inverse) sequence to $\mathcal{O}$. If $\mathcal{S}=\mathcal{P}_m$, we say that the pivot strategy $I_{\mathcal{O}}^{\leftarrow}=I_{\mathcal{O}^{\leftarrow}}$ is reverse (inverse) to $I_{\mathcal{O}}$.
\end{Definition}

Obviously, we have $\mathcal{O}^{\leftarrow\leftarrow}=\mathcal{O}$ and hence $I_{\mathcal{O}^{\leftarrow\leftarrow}}=I_{\mathcal{O}}$ for $\mathcal{O}\in\mathcal{\Ol}(\mathcal{P}_m)$.

\begin{Lemma}\label{tm: lema_2.3}
Let $\mathcal{O},\mathcal{O}'\in\mathcal{\Ol}(\mathcal{S})$, $\mathcal{S}\subseteq\mathcal{P}_m$. Then $\mathcal{O}'\stackrel{\mathsf{w}}{\sim}\mathcal{O}$ if and only if
${\mathcal{O}'}^{\leftarrow} \stackrel{\mathsf{w}}{\sim} \mathcal{O}^{\leftarrow}$.
\end{Lemma}

\begin{proof}
From formula~\eqref{weak_eq}, we see that it is sufficient to prove the assertion for relations $\sim$ and $\stackrel{\mathsf{s}}{\sim}$. Let $\mathcal{O}$ and $\mathcal{O}^{\leftarrow}$ be as in Definition~\ref{tm: def_reverse}.

Suppose $\mathcal{O}'$ is obtained from $\mathcal{O}$ by applying one admissible transposition.
Then for some $0\leq t< r$ we have $\{i_t,j_t\}\cap \{i_{t+1},j_{t+1}\}=\emptyset$. If $t\geq 1$, then
\begin{align*}
\mathcal{O}' & = (i_0,j_0),\ldots,(i_{t-1},j_{t-1}),(i_{t+1},j_{t+1}),(i_t,j_t),\ldots,(i_{r},j_{r}), \\
{\mathcal{O}'}^{\leftarrow} & = (i_{r},j_{r}),\ldots,(i_t,j_t),(i_{t+1},j_{t+1}),(i_{t-1},j_{t-1}),\ldots,(i_0,j_0),
\end{align*}
and obviously ${\mathcal{O}'}^{\leftarrow}\sim\mathcal{O}^{\leftarrow}$. If $t=0$,
$\mathcal{O}'= (i_1,j_1), (i_0,j_0),\ldots,(i_{r},j_{r})$ and ${\mathcal{O}'}^{\leftarrow} =
(i_{r},j_{r}),\ldots,$ \linebreak $(i_0,j_0),(i_1,j_1)$, so we also have ${\mathcal{O}'}^{\leftarrow}\sim \mathcal{O}^{\leftarrow}$. If $\mathcal{O}'$ is the result of applying more than one admissible transposition to $\mathcal{O}$, then the proof proceeds by applying the above argument several times.

Suppose that ${\mathcal{O}'}^{\leftarrow}\sim \mathcal{O}^{\leftarrow}$ holds. By the implication proved in the preceding paragraph we have
${\mathcal{O}'}^{\leftarrow\leftarrow}\sim \mathcal{O}^{\leftarrow\leftarrow}$, and this is the same as $\mathcal{O}'\sim \mathcal{O}$.

Now, let $\mathcal{O}'\stackrel{\mathsf{s}}{\sim}\mathcal{O}$. Suppose
$\mathcal{O}'=(i_{t+1},j_{t+1}),\ldots,(i_{r},j_{r}),(i_{0},j_{0}),\ldots,(i_{t},j_{t})$. Then
${\mathcal{O}'}^{\leftarrow}=(i_{t},j_{t}),\ldots,(i_0,j_0),(i_{r},j_{r}),\ldots,(i_{t+1},j_{t+1})$, which is shift-equivalent to $\mathcal{O}^{\leftarrow}=(i_{r},j_{r}),\ldots,(i_1,j_1),$ \linebreak $(i_0,j_0)$ with shift equal to $r-t$. Now we know that
${\mathcal{O}'}^{\leftarrow}\stackrel{\mathsf{s}}{\sim}\mathcal{O}^{\leftarrow}$ implies
${\mathcal{O}'}^{\leftarrow\leftarrow}\stackrel{\mathsf{s}}{\sim} \mathcal{O}^{\leftarrow\leftarrow}$ and this is the same as
$\mathcal{O}'\stackrel{\mathsf{s}}{\sim} \mathcal{O}$.
$\qquad$ \end{proof}

To visually depict an ordering $\mathcal{O}$ of $\mathcal{P}_m$ we make use of the symmetric matrix $\mathsf{M}_{\mathcal{O}}=(\mathsf{m}_{rt})$ of order $m$, defined by the rule
\[
\mathsf{m}_{i(k)j(k)}=\mathsf{m}_{j(k)i(k)}=k, \quad k=0,1,\ldots,M-1; \quad M=m(m-1)/2.
\]
We set $\mathsf{m}_{rr}=-1$, $1\leq r\leq m$, but since the pairs $(r,r)$ do not appear in $\mathcal{O}$, we will rather use $*$ to represent $-1$.

\begin{Example}
As an illustration, we depict the matrices
$\mathsf{M}_{\mathcal{O}^{\leftarrow}_{\text{row}}}$ and $\mathsf{M}_{\mathcal{O}^{\leftarrow}_{\text{col}}}$ for $m=5$.
\[
\mathsf{M}_{\mathcal{O}^{\leftarrow}_r}=\left[
              \begin{array}{ccccc}
                * & 9 & 8 & 7 & 6  \\
                9 & * & 5 & 4 & 3  \\
                8 & 5 & * & 2 & 1  \\
                7 & 4 & 2 & * & 0  \\
                6 & 3 & 1 & 0 & *
              \end{array}\right], \qquad
\mathsf{M}_{\mathcal{O}^{\leftarrow}_c}= \left[
              \begin{array}{ccccc}
                * &  9 & 8 & 6 & 3  \\
                9 &  * & 7 & 5 & 2  \\
                8 &  7 & * & 4 & 1  \\
                6 &  5 & 4 & * & 0  \\
                3 &  2 & 1 & 0 & *
              \end{array}\right].
\]
These two matrices give us information on the order in which the off-diagonal blocks
in the block-matrix from~\eqref{blokmatrica} are  annihilated during each cycle.
\end{Example}

\subsubsection{Permutation equivalent strategies}\label{sec:2.2.1}

Let us introduce yet another equivalence relation on $\mathcal{\Ol}(\mathcal{P}_m)$ and on the set on cyclic and quasi-cyclic pivot strategies.

Two pivot orderings $\mathcal{O},\mathcal{O}'\in\mathcal{\Ol}(\mathcal{P}_m)$ are
permutation equivalent if $\mathsf{M}_{\mathcal{O}'}=\mathsf{P} \mathsf{M}_{\mathcal{O}}\mathsf{P}^T$ holds for some permutation matrix $\mathsf{P}$. In that case we write   $\mathcal{O}'\stackrel{\mathsf{p}}{\sim}\mathcal{O}$ and
$I_{\mathcal{O}'}\stackrel{\mathsf{p}}{\sim}I_{\mathcal{O}}$.
Recall that each permutation matrix $\mathsf{P}$ of order $m$ is defined by some permutation $\mathsf{p}$ of the set $\mathcal{S}_m=\{1,2,\ldots,m\}$ via the relation
\begin{equation}\label{perm}
\mathsf{P}e_r = e_{\mathsf{p}(r)}, \quad 1\leq r\leq m.
\end{equation}
Here $I_m=[e_1,\ldots,e_m]$. The mapping $\mathsf{p}\mapsto \mathsf{P}$ is an isomorphism between the symmetric group on the set $\mathcal{S}_m$ and the group of permutation matrices of order $m$. If $X=(x_{rt})$ is any square matrix of order $m$, then  $\mathsf{P}X\mathsf{P}^T=(x_{\mathsf{p}^{-1}(r),\mathsf{p}^{-1}(t)})$. Hence, if  $\widetilde{\mathcal{O}}$ is permutation equivalent to $\mathcal{O}$ and $\mathsf{M}_{\mathcal{O}}=(\mathsf{m}_{rt})$,
$\mathsf{M}_{\widetilde{\mathcal{O}}}=(\widetilde{\mathsf{m}}_{rt})$, then $\widetilde{\mathsf{m}}_{\mathsf{p}(r)\mathsf{p}(t)}=\mathsf{m}_{rt}$ holds for all $r$, $t$.
This relation shows that the $(r,t)$-element of $\mathsf{M}_{\mathcal{O}}$ becomes the $(\mathsf{p}(r),\mathsf{p}(t))$-element of $\mathsf{M}_{\widetilde{\mathcal{O}}}$. Therefore, if
$\widetilde{\mathcal{O}}\stackrel{\mathsf{p}}{\sim}\mathcal{O}$ and
$\displaystyle \mathcal{O}=(i_0,j_0),\ldots,(i_{M-1},j_{M-1})$, then
$\displaystyle \widetilde{\mathcal{O}}=(\mathsf{p}(i_0),\mathsf{p}(j_0)),\ldots,$ $(\mathsf{p}(i_{M-1}),\mathsf{p}(j_{M-1}))$.
Here it is presumed that in the case $\mathsf{p}(i)>\mathsf{p}(j)$, the pair $(\mathsf{p}(i),\mathsf{p}(j))$ in the ordering $\widetilde{\mathcal{O}}$ is replaced by $(\mathsf{p}(j),\mathsf{p}(i))$.

\begin{Example}
Let $m=4$, $\mathcal{S}_m =\{1,2,3,4\}$, $\mathcal{O} = (1,2),(2,3),(2,4),(3,4),(1,3)$, $(1,4)$ and $\mathsf{P} = [e_2,e_4,e_3,e_1]$, so that
$\mathsf{p}=\displaystyle \left(\begin{array}{cccc} 1&2&3&4 \\ 2&4&3&1\end{array} \right)$,
$\mathsf{p}^{-1}=\left(\begin{array}{cccc} 1&2&3&4 \\ 4&1&3&2\end{array} \right)$,
$\mathsf{P}^T = [e_4,e_1,e_3,e_2]$.
Let $\widetilde{\mathcal{O}}$ be such that
$\mathsf{M}_{\widetilde{\mathcal{O}}}=\mathsf{P}\mathsf{M}_{\mathcal{O}}\mathsf{P}^T$. Since
\[
\mathsf{M}_{\widetilde{\mathcal{O}}}=
\left[\begin{array}{c}e_4^T \\ e_1^T \\ e_3^T \\ e_2^T \end{array}\right]
\left[\begin{array}{cccc}*&0&4&5\\ 0&*&1&2\\ 4&1&*&3\\ 5&2&3&*\end{array}\right]
[e_4,e_1,e_3,e_2]=
\left[\begin{array}{cccc}*&5&3&2\\ 5&*&4&0\\ 3&4&*&1\\ 2&0&1&*\end{array}\right],
\]
we have $\widetilde{\mathcal{O}}=(2,4),(3,4),(1,4),(1,3),(2,3),(1,2)$.
On the other hand, we have
\begin{align*}
\mathcal{O}(\mathsf{p}) & \equiv (\mathsf{p}(1),\mathsf{p}(2)),(\mathsf{p}(2),\mathsf{p}(3)),(\mathsf{p}(2),\mathsf{p}(4)),\!(\mathsf{p}(3),\mathsf{p}(4)),(\mathsf{p}(1),\mathsf{p}(3)),\!(\mathsf{p}(1),\mathsf{p}(4)) \\
& = (2,4),(3,4),(1,4),(1,3),(2,3),(1,2) = \widetilde{\mathcal{O}}.
\end{align*}
\end{Example}

Now it is easy to extend the notion of permutation equivalence from pivot orderings to pivot sequences from
$\mathcal{\Ol}(\mathcal{P}_m)$.

\begin{Definition}
Let $\mathcal{O}=(i_0,j_0), (i_1,j_1),\ldots,(i_{T-1},j_{T-1})\in\mathcal{\Ol}(\mathcal{P}_m)$, $T\geq M$. The sequence
$\mathcal{O}'\in\mathcal{\Ol}(\mathcal{P}_m)$ is permutation equivalent to $\mathcal{O}$, and we write $\mathcal{O}'\stackrel{\mathsf{p}}{\sim}\mathcal{O}$, if there is a permutation $\mathsf{q}$ of the set $\mathcal{S}_m$ such that
$\mathcal{O}'=(\mathsf{q}(i_0),\mathsf{q}(j_0)),(\mathsf{q}(i_1),\mathsf{q}(j_1)),\ldots,$ $(\mathsf{q}(i_{T-1}),\mathsf{q}(j_{T-1}))$. Then $\mathcal{O}'$ is denoted by $\mathcal{O}(\mathsf{q})$.
\end{Definition}

Since $\mathcal{O} = \mathcal{O}(\mathsf{e})$, where $\mathsf{e}$ is the identity permutation, we have $\mathcal{O}\stackrel{\mathsf{p}}{\sim}\mathcal{O}$. If $\mathcal{O}'=\mathcal{O}(\mathsf{q})$ then $\mathcal{O}=\mathcal{O}'(\mathsf{q}^{-1})$.
Also if $\mathcal{O}'=\mathcal{O}(\mathsf{q})$ and $\mathcal{O}''=\mathcal{O}'(\mathsf{q'})$ then $\mathcal{O}''=\mathcal{O}(\mathsf{q'}\circ \mathsf{q})$ where $\circ$ denotes the binary operation in the permutation group, which is simply the composition of functions. We conclude that $\stackrel{\mathsf{p}}{\sim}$ is an equivalence relation on the set $\mathcal{\Ol}(\mathcal{P}_m)$.
Note that $\mathcal{O}_r = \mathcal{O}_{r-1}(\mathsf{q}_r)$, $1\leq r\leq t$, implies $\mathcal{O}_t = \mathcal{O}_{0}(\mathsf{q})$
with $\mathsf{q}=\mathsf{q}_t\circ \mathsf{q}_{t-1}\circ\cdots\circ\mathsf{q}_1$.

\begin{Lemma}
Let $\mathcal{O},\mathcal{O}_1,\mathcal{O}_2,\mathcal{O}_3,\mathcal{O}_4 \in\mathcal{\Ol}(\mathcal{P}_m)$.
\begin{itemize}
\item[(i)] If $\mathcal{O}\stackrel{\mathsf{w}}{\sim}\mathcal{O}_1\stackrel{\mathsf{p}}{\sim}\mathcal{O}_2$,
then there is $\mathcal{O}'\in\mathcal{\Ol}(\mathcal{P}_m)$ such that
$\mathcal{O}\stackrel{\mathsf{p}}{\sim}\mathcal{O}'\stackrel{\mathsf{w}}{\sim} \mathcal{O}_2$.
\item[(ii)] If $\mathcal{O}\stackrel{\mathsf{p}}{\sim} \mathcal{O}_3\stackrel{\mathsf{w}}{\sim}\mathcal{O}_4$,
then there is $\widetilde{\mathcal{O}}\in\mathcal{\Ol}(\mathcal{P}_m)$ such that
$\mathcal{O}\stackrel{\mathsf{w}}{\sim}\widetilde{\mathcal{O}}\stackrel{\mathsf{p}}{\sim} \mathcal{O}_4$.
\end{itemize}
\end{Lemma}

\begin{proof}
$(i)$ It is sufficient to prove this assertion for the two cases:
(a) $\stackrel{\mathsf{w}}{\sim}$ is reduced to $\sim$
and (b) $\stackrel{\mathsf{w}}{\sim}$ is reduced to $\stackrel{\mathsf{s}}{\sim}$.

(a) It is sufficient to assume that $\mathcal{O}_1$ is obtained from $\mathcal{O}$ by applying one admissible transposition.
Let
$\mathcal{O}=(i_0,j_0),\ldots,(i_r,j_r),(i_{r+1},j_{r+1}),\ldots,(i_{T-1},j_{T-1})$ with $\{i_r,j_r\}$ $\cap\{i_{r+1},j_{r+1}\}=\emptyset$, so that
$\mathcal{O}_1=(i_0,j_0),\ldots,(i_{r+1},j_{r+1}),(i_r,j_r),\ldots,(i_{T-1},j_{T-1})$. If $\mathcal{O}_2=\mathcal{O}_1(\mathsf{q})$, then $\mathcal{O}_2=(\mathsf{q}(i_0),\mathsf{q}(j_0)),\ldots,(\mathsf{q}(i_{r+1}),\mathsf{q}(j_{r+1})),(\mathsf{q}(i_r),\mathsf{q}(j_r)),\ldots,$ $(\mathsf{q}(i_{T-1}),$ $\mathsf{q}(j_{T-1}))$.
Since $\mathsf{q}$ is a bijection from $\mathcal{S}_m$ onto itself, we have $\{\mathsf{q}(i_r),\mathsf{q}(j_r)\}\cap\{\mathsf{q}(i_{r+1}),$ $\mathsf{q}(j_{r+1})\}=\emptyset$. Therefore, we can set
$\mathcal{O}'=\mathcal{O}(\mathsf{q})$, i.e.,
\[
\mathcal{O}'=(\mathsf{q}(i_0),\mathsf{q}(j_0)),\ldots,(\mathsf{q}(i_r),\mathsf{q}(j_r)),(\mathsf{q}(i_{r+1}),\mathsf{q}(j_{r+1})),\ldots,(\mathsf{q}(i_{T-1}),\mathsf{q}(j_{T-1})).
\]

(b) Let $\mathcal{O}$ be as in case (a) and assume that the cut has been made behind the term $(i_r,j_r)$, so that
$\mathcal{O}_1=(i_{r+1},j_{r+1}),\ldots,(i_{T-1},j_{T-1}),(i_0,j_0),\ldots,(i_r,j_r)$.
Let $\mathcal{O}_2=\mathcal{O}_1(\mathsf{q})=(\mathsf{q}(i_{r+1}),\mathsf{q}(j_{r+1})),$ \linebreak $\ldots,(\mathsf{q}(i_{T-1}),\mathsf{q}(j_{T-1})),(\mathsf{q}(i_0),\mathsf{q}(j_0)),\ldots,(\mathsf{q}(i_r),\mathsf{q}(j_r))$. Obviously, we can define $\mathcal{O}'$ to be the same as in case (a). To obtain $\mathcal{O}_2$ from $\mathcal{O}'$, one has to make the cut just behind the term $(\mathsf{q}(i_r),\mathsf{q}(j_r))$, i.e., to use the shift $r+1$.

(ii) The proof is quite similar to the proof of $(i)$. First, consider the case when $\stackrel{\mathsf{w}}{\sim}$ is reduced to $\sim$, and $\sim$ is given by one admissible transposition that interchanges the terms at positions $r$ and $r+1$, as in case (a) above. Let  $\mathcal{O}$ be as in case $(i)$ and denote $\mathcal{O}_3=\mathcal{O}(\widetilde{\mathsf{q}})$. Then
$\mathcal{O}_4$ is the same as $\mathcal{O}_2$, provided that $\mathsf{q}$ is replaced by $\widetilde{\mathsf{q}}$ and
$\widetilde{\mathcal{O}}$ is like $\mathcal{O}_1$.
Next, consider the case when $\stackrel{\mathsf{w}}{\sim}$ is reduced to $\stackrel{\mathsf{s}}{\sim}$ and  $\mathcal{O}_3=\mathcal{O}(\widetilde{\mathsf{q}})$. Assume that the cut has been made as in case (b), behind the term $(\widetilde{\mathsf{q}}(i_r),\widetilde{\mathsf{q}}(j_r))$.
Then $\mathcal{O}_4$ is the same as $\mathcal{O}_2$ (from case (b)), provided that $\mathsf{q}$ is replaced by $\widetilde{\mathsf{q}}$ and $\widetilde{\mathcal{O}}$ is like $\mathcal{O}_1$.
$\qquad$ \end{proof}

\begin{Proposition}\label{tm:peq1}
Let $\mathcal{O},\mathcal{O}',\mathcal{O}_1,\mathcal{O}_2,\mathcal{O}_3,\ldots,\mathcal{O}_{2t} \in\mathcal{\Ol}(\mathcal{P}_m)$, $t\geq 1$. If
\begin{equation}\label{chain}
\mathcal{O} \stackrel{\mathsf{w}}{\sim} \mathcal{O}_1\stackrel{\mathsf{p}}{\sim}\mathcal{O}_2 \stackrel{\mathsf{w}}{\sim}\mathcal{O}_3\stackrel{\mathsf{p}}{\sim}\mathcal{O}_4
\stackrel{\mathsf{w}}{\sim}\cdots\stackrel{\mathsf{w}}{\sim} \mathcal{O}_{2t-1}\stackrel{\mathsf{p}}{\sim} \mathcal{O}_{2t}\stackrel{\mathsf{w}}{\sim}\mathcal{O}',
\end{equation}
then there exist $\mathcal{O}_0',\widetilde{\mathcal{O}}_0\in\mathcal{\Ol}(\mathcal{P}_m)$ such that
$\mathcal{O}\stackrel{\mathsf{p}}{\sim} \mathcal{O}_0' \stackrel{\mathsf{w}}{\sim}\mathcal{O}'$ and
$\mathcal{O}\stackrel{\mathsf{w}}{\sim} \widetilde{\mathcal{O}}_0\stackrel{\mathsf{p}}{\sim}\mathcal{O}'$.

If
$\mathcal{O}_{2r}=\mathcal{O}_{2r-1}(\mathsf{q}_r)$, $1\leq r\leq t$, then
$\mathcal{O}_0' = \mathcal{O}(\mathsf{q})$ and $\mathcal{O}' = \widetilde{\mathcal{O}}_0(\mathsf{q})$ with $\mathsf{q}=\mathsf{q}_t\circ\mathsf{q}_{t-1}\circ\cdots\circ\mathsf{q}_1$.
\end{Proposition}

\begin{proof} Using assertion $(i)$ (resp.\@ assertion $(ii)$) of the previous lemma, one can gradually move all appearances of $\stackrel{\mathsf{p}}{\sim}$ to the left (resp.\@ right) end of the chain (\ref{chain}). First, the leftmost (resp.\@ rightmost) $\stackrel{\mathsf{p}}{\sim}$ is moved. The leftmost (resp.\@ rightmost) part of the chain takes the form
\[
\mathcal{O}\stackrel{\mathsf{p}}{\sim}\mathcal{O}_1'\stackrel{\mathsf{w}}{\sim}\mathcal{O}_2 \stackrel{\mathsf{w}}{\sim}\mathcal{O}_3\stackrel{\mathsf{p}}{\sim}\mathcal{O}_4\cdots \qquad (\text{resp.} \ \cdots\mathcal{O}_{2t-3}\stackrel{\mathsf{p}}{\sim}\mathcal{O}_{2t-2} \stackrel{\mathsf{w}}{\sim}\mathcal{O}_{2t-1}\stackrel{\mathsf{w}}{\sim}\widetilde{\mathcal{O}}_{2t}
 \stackrel{\mathsf{p}}{\sim}\mathcal{O}').
\]
Note that two consecutive $\stackrel{\mathsf{w}}{\sim}$ can be replaced by one, so we can remove $\mathcal{O}_2$ (resp.\@ $\mathcal{O}_{2t-1}$) from the obtained chain. Next, $\stackrel{\mathsf{p}}{\sim}$ that links $\mathcal{O}_3$ and $\mathcal{O}_4$
(resp.\@ $\mathcal{O}_{2t-3}$ and $\mathcal{O}_{2t-2}$) is moved. The leftmost (resp.\@ rightmost) part of the chain takes the form
\[
\mathcal{O}\stackrel{\mathsf{p}}{\sim}\mathcal{O}_1'\stackrel{\mathsf{p}}{\sim}\mathcal{O}_3' \stackrel{\mathsf{w}}{\sim}\mathcal{O}_4\stackrel{\mathsf{w}}{\sim}\mathcal{O}_5\cdots \qquad (\text{resp.} \ \cdots\mathcal{O}_{2t-4} \stackrel{\mathsf{w}}{\sim}\mathcal{O}_{2t-3} \stackrel{\mathsf{w}}{\sim}\widetilde{\mathcal{O}}_{2t-2}\stackrel{\mathsf{p}}{\sim}\widetilde{\mathcal{O}}_{2t} \stackrel{\mathsf{p}}{\sim}\mathcal{O}').
\]
Continuing this way one ultimately obtains
\begin{align}\label{chain_fin}
& \mathcal{O}\stackrel{\mathsf{p}}{\sim}\mathcal{O}_1'\stackrel{\mathsf{p}}{\sim}\mathcal{O}_3'
\stackrel{\mathsf{p}}{\sim}\cdots\stackrel{\mathsf{p}}{\sim} \mathcal{O}_{2t-1}'\stackrel{\mathsf{w}}{\sim}\mathcal{O}_{2t}\stackrel{\mathsf{w}}{\sim}\mathcal{O}' \\
& (\text{resp.} \ \mathcal{O} \stackrel{\mathsf{w}}{\sim}\mathcal{O}_{1}\stackrel{\mathsf{w}}{\sim}\widetilde{\mathcal{O}}_{2}\stackrel{\mathsf{p}}{\sim}\cdots
\stackrel{\mathsf{p}}{\sim}\widetilde{\mathcal{O}}_{2t-2}\stackrel{\mathsf{p}}{\sim}\widetilde{\mathcal{O}}_{2t}\stackrel{\mathsf{p}}{\sim}\mathcal{O}'). \nonumber
\end{align}
Here $\mathcal{O}_{2t}$ (resp.\@ $\mathcal{O}_{1}$) can be removed. Note that $\stackrel{\mathsf{p}}{\sim}$ is an equivalence relation. Hence by the transitivity property, the leftmost part of the chain
$\mathcal{O}\stackrel{\mathsf{p}}{\sim}\mathcal{O}_1'\stackrel{\mathsf{p}}{\sim}\mathcal{O}_3'
\stackrel{\mathsf{p}}{\sim}\cdots\stackrel{\mathsf{p}}{\sim}\mathcal{O}_{2t-1}'$ can be replaced by $\mathcal{O}\stackrel{\mathsf{p}}{\sim}\mathcal{O}_{2t-1}'$ (and similarly for the rightmost part of the chain). To complete the proof of the first assertion one has to rename $\mathcal{O}_{2t-1}'$ as $\mathcal{O}_{0}'$
(resp.\@ $\widetilde{\mathcal{O}}_2$ as $\widetilde{\mathcal{O}}_0$)

The proof of the second assertion is the same, but one can use more information. Now, from the proof of the preceding lemma we know that in the final chain (\ref{chain_fin}) we have $\mathcal{O}_{1}'=\mathcal{O}(\mathsf{q}_1)$ and $\mathcal{O}_{2r-1}'=\mathcal{O}_{2r-3}'(\mathsf{q}_r)$, $2\leq r\leq t$
(resp.\@ $\mathcal{O}'=\widetilde{\mathcal{O}}(\mathsf{q}_t)$ and $\widetilde{\mathcal{O}}_{2r}=\widetilde{\mathcal{O}}_{2r-2}(\mathsf{q}_{r-1})$, $2\leq r\leq t$).
Hence $\mathcal{O}_{0}'=\mathcal{O}_{2t-1}'=\mathcal{O}(\mathsf{q})$ (resp.\@ $\mathcal{O}'=\widetilde{\mathcal{O}}_{2}(\mathsf{q})=\widetilde{\mathcal{O}}_{0}(\mathsf{q})$) with $\mathsf{q}=\mathsf{q}_t\circ \mathsf{q}_{t-1}\circ\cdots\circ\mathsf{q}_1$.
$\qquad$ \end{proof}

Obviously, the first or the last (or both) appearance of $\stackrel{\mathsf{w}}{\sim}$ in the chain (\ref{chain}) can be omitted provided that
$\mathcal{O}=\mathcal{O}_1$ and $\mathcal{O}' = \mathcal{O}_{2t}$  (resp.\@ $\mathcal{O}=\mathcal{O}_1$ and $\mathcal{O}' = \mathcal{O}_{2t}$).

Two sequences $\mathcal{O},\mathcal{O}'\in\mathcal{\Ol}(\mathcal{P}_m)$ can be linked via a long
chain like (\ref{chain}), which may include all equivalence relations $\sim$, $\stackrel{\mathsf{s}}{\sim}$, $\stackrel{\mathsf{w}}{\sim}$ and $\stackrel{\mathsf{p}}{\sim}$ introduced so far.
Proposition~\ref{tm:peq1} shows that each such chain can be reduced to a short chain that uses just one
$\stackrel{\mathsf{p}}{\sim}$ and one $\stackrel{\mathsf{w}}{\sim}$. Furthermore, weak equivalence can be written in the most compact form (\ref{weak_eq}).

\begin{Definition}\label{tm:can_chain}
Two sequences $\mathcal{O}, \mathcal{O}'\in\mathcal{\Ol}(\mathcal{P}_m)$ are connected by a \emph{chain of equivalence relations} if there exist
$\mathcal{O}_1,\ldots ,\mathcal{O}_r\in \mathcal{\Ol}(\mathcal{P}_m)$, $r\geq 0$, such that in the sequence
$\mathcal{O}$, $\mathcal{O}_1,\ldots,\mathcal{O}_r,\mathcal{O}'$ each two neighboring terms are linked with
$\sim$, $\stackrel{\mathsf{s}}{\sim}$, $\stackrel{\mathsf{w}}{\sim}$ or $\stackrel{\mathsf{p}}{\sim}$.
The chain is in \emph{canonical form} if it looks like (\ref{weak_eq}) with one $\stackrel{\mathsf{p}}{\sim}$ placed in front of it or behind it.
\end{Definition}

We conclude that every chain of equivalence relations can be reduced to the canonical form.

Let us return to the block matrix $A$ from relation (\ref{blokmatrica}). To a partition $\pi=(n_1,\ldots,n_m)$ we associate the $n$-tuple $(s_1,\ldots,s_m)$,
\[
s_i = s_{i-1}+n_{i}, \quad 2\leq i\leq m; \quad  s_1=n_1.
\]
Note that the sequence of the first $n$ natural numbers, $1,2,\ldots,n$, can be written as $1,\ldots,s_1,$ $s_1+1,\ldots ,s_2,\ldots,s_{m-1}+1,\ldots ,s_m$, which is the same as
\[
s_1-n_1+1,\ldots,s_1, s_2-n_2+1,\ldots,s_2,\ldots,s_{m}-n_m+1,\ldots,s_m.
\]
Here, $s_r-n_r+1,\ldots,s_r$ are indices of the columns (resp.\@ rows) which make the $r$th block-column (resp.\@ block-row) of $A$.
Let  $\mathcal{O}\in\mathcal{P}_m$ and $\widetilde{\mathcal{O}}=\mathcal{O}(\mathsf{p})$, so that $\widetilde{\mathcal{O}}=$ $(\mathsf{p}(i_0),\mathsf{p}(j_0))$, $(\mathsf{p}(i_1),\mathsf{p}(j_1))$, \ldots, $(\mathsf{p}(i_{T-1}),\mathsf{p}(j_{T-1}))$.
The permutation $p$ of the set $\mathcal{S}_n=\{1,2,\ldots,n\}$ associated with $\mathsf{p}$ is defined by:
{\small
\[
\mbox{\normalsize $p$}=\left(\!\!\begin{array}{cccccccccc}
s_1-n_1+1, & \kern-1em \ldots & s_1, & \! s_2-n_2+1, & \kern-1em \ldots & s_{m}-n_m+1, & \kern-1em \ldots & \! s_{m} \\
s_{\mathsf{p}(1)}-n_{\mathsf{p}(1)}+1, & \kern-1em \ldots & s_{\mathsf{p}(1)}, & \! s_{\mathsf{p}(2)}-n_{\mathsf{p}(2)}+1, &
\kern-1em \ldots & s_{\mathsf{p}(m)}-n_{\mathsf{p}(m)}+1, & \kern-1em \ldots & \! s_{\mathsf{p}(m)}
\end{array}\!\! \right).
\]
}
Using the same rule (\ref{perm}), we obtain the permutation matrix $P$ of order $n$,
associated with $p$. It satisfies $Pe_t=e_{p(t)}$, $1\leq t\leq n$. The matrix $P$ has the form $P=[E_{\mathsf{p}(1)} \cdots E_{\mathsf{p}(m)}]$, where each $E_{\mathsf{p}(k)}$ is an $n\times n_{\mathsf{p}(k)}$ matrix (i.e., a single block-column) that differs from the zero matrix only in its $\mathsf{p}(k)$th block-row,
\begin{eqnarray*}
E_{\mathsf{p}(k)}\ &=&\ \left[\begin{array}{ccccccc}
\ 0\ &\ \cdots\ &\ 0\ &\ I_{\mathsf{p}(k)}\ &\ 0\ &\ \cdots\ &\ 0\
\end{array}\right]^T,\quad  1\leq k\leq m.\\
&& \hspace{3ex}n_1\hspace{13ex}n_{\mathsf{p}(k)}\hspace{13ex}n_m
\end{eqnarray*}
Let $A^{(T)}$ be the matrix obtained from $A$ by applying one sweep of the quasi-cyclic block Jacobi method defined by $I_{\mathcal{O}}$. The iterative process has the form (\ref{jacobiblokagm}), hence using $P$ and $P^T$ we can write
\begin{align}
PA^{(T)}P^T & = P(U_{T-1}^TU_{T-2}^T\cdots U_0^TP^T)\, PAP^T\, (PU_0\cdots U_{T-2}U_{T-1})P^T \label{p-process}\\
& = (PU_{T-1}^TP^T)\cdots(PU_0^TP^T)\,P A P^T(PU_0P^T)\cdots(PU_{T-1}P^T) \nonumber \\
& = \widetilde{U}_{T-1}^T\cdots \widetilde{U}_0^T (PAP^T)\, \widetilde{U}_0\cdots \widetilde{U}_{T-1},\nonumber
\end{align}
where
\begin{equation}\label{UTAU}
\widetilde{U}_k=PU_kP^T, \quad 0\leq k\leq T-1.
\end{equation}
Each $\widetilde{U}_k$ is an orthogonal elementary block matrix whose pivot pair is $(\mathsf{p}(i_k),\mathsf{p}(j_k))$. We can interpret the process (\ref{p-process}) as a quasi-cyclic block Jacobi method defined by $I_{\widetilde{\mathcal{O}}}$. When it is applied to $\widetilde{A}=PAP^T$, after one sweep it results in $\widetilde{A}^{(T)}=PA^{(T)}P^T$. Indeed, at step $k$ of that process we have
\begin{equation}\label{A'=UTAU}
\widetilde{A}^{(k+1)}= \widetilde{U}_{k}^T \widetilde{A}^{(k)}\widetilde{U}_{k},\quad 0\leq k\leq T-1,
\end{equation}
where $\widetilde{A}^{(k)}= P A^{(k)} P^T$ for any $k$. For the matrix $\widetilde{A}^{(k+1)}$ we know that its $(\mathsf{p}(i_k),\mathsf{p}(j_k))$-pivot submatrix of order $n_{\mathsf{p}(i_k)}+n_{\mathsf{p}(j_k)}$ is diagonal. Hence, it is a quasi-cyclic block Jacobi method. Moreover, if $U_k$ is in the class $\ubce_{\pi}(\varrho)$ from Section~2.3, the same will be true for  $\widetilde{U}_k$.

The process (\ref{A'=UTAU}) is a block Jacobi method for the matrix $\widetilde{A}$ which
carries block-matrix partition defined by $\pi_{\mathsf{p}} = (n_{\mathsf{p}(1)},\ldots,n_{\mathsf{p}(m)})$. Thus, if $A$ is replaced by $\widetilde{A}$ and $\pi$ by $\pi_{\mathsf{p}}$,
then the  block method (\ref{A'=UTAU}) is defined by the pivot sequence $\widetilde{\mathcal{O}}=\mathcal{O}(\mathsf{p})\in\mathcal{\Ol}(\mathcal{P}_m)$. We can formally write $(\pi,\mathcal{O},A)$ $\stackrel{\mathsf{p}}{\mapsto}$ $(\pi_{\mathsf{p}},\mathcal{O}(\mathsf{p}), PAP^T)$. This is equivalent to
$(\pi_{\mathsf{p}^{-1}},$ $\mathcal{O}(\mathsf{p^{-1}}), P^TAP)$ $\stackrel{\mathsf{p}}{\mapsto}$ $(\pi, \mathcal{O},A)$.

We end this subsection with two remarks.

First, the reverse ordering $\mathcal{O}^{\leftarrow}$ is not the same as $\mathcal{O}(\widetilde{\mathsf{e}})$, where
\begin{equation}\label{etilde}
\widetilde{\mathsf{e}}=\left(\begin{array}{cccc} 1 & 2 & \ldots & m \\ m & m-1 & \ldots & 1\end{array}\right).
\end{equation}
Examples which confirm this claim are those from Section~\ref{sec:3.2}.

Second, if in Lemma~\ref{tm: lema_2.3} the equivalence relation $\stackrel{\mathsf{w}}{\sim}$
is replaced by $\stackrel{\mathsf{p}}{\sim}$, the assertion will remain to hold. The proof is trivial.
Indeed, if in the chain (\ref{weak_eq}), which may now include $\stackrel{\mathsf{p}}{\sim}$,
only one sequence $\mathcal{O}_t$ is replaced by $\mathcal{O}_t^{\leftarrow}$, then all the sequences have to be replaced by their inverses. Otherwise, the chain can be broken into two chains  that are not mutually connected.

\subsection{Global convergence}

A block Jacobi method is \emph{convergent} on $A$ if the obtained sequence of matrices $(A^{(k)})$ converges to some diagonal matrix $\Lambda$. The method is \emph{globally convergent} if it is convergent on every symmetric matrix $A$. This definition assumes that the partition $\pi$ is arbitrary. In particular one can take
$\pi =(1,1,\ldots,1)$, which means that it is the proper generalization of the standard notion of global convergence. The words ``global'' and ``globally'' are often omitted. For example, if one says that the block method converges for some pivot strategy, this means that the method converges for every initial symmetric matrix. For the global convergence considerations, it is irrelevant whether the diagonal blocks of the initial matrix are diagonal submatrices. Namely, after some iteration (within the first sweep) this property will be fulfilled and it will remain to hold until convergence.
To measure how much the method has converged, we use the quantity
$$S(A)= \frac{\sqrt{2}}{2}\|A -\text{diag}(A)\|_F = \left[\sum_{s=1}^{n-1} \sum_{t=s+1}^n |a_{st}|^2 \right]^{\frac{1}{2}},$$
where $\|X\|_F=\sqrt{\text{trace}(X^T X)}$ stands for the Frobenius norm of $X$. In the definition of $S(A)$ we could have used blocks instead of elements, but since the diagonal blocks are diagonal submatrices, this reduces to the same quantity.
Obviously, the convergence of a block Jacobi method on $A$ implies that $S(A^{(k)})\rightarrow 0$ as $k\rightarrow\infty$.
The converse is true provided that the diagonal elements of
$\text{diag} (\Lambda_{ii}^{(k+1)},\Lambda_{jj}^{(k+1)})$ from (\ref{blokpivotmatr}) are always ordered in some prescribed order, typically nonincreasingly.

\begin{Theorem} \label{tm:gl_conv-opp}
Let $A$ be a symmetric matrix and $A^{(k)}$, $k\geq 0$, be the sequence obtained by applying the block Jacobi method to $A$. Let the pivot strategy be cyclic or quasi-cyclic and assume that $\lim_{k\rightarrow\infty}S(A^{(k)})=0$.
\begin{itemize}
\item[(i)] If the algorithm that diagonalizes the pivot submatrix always delivers
$\text{diag} (\Lambda_{ii}^{(k+1)},$ $\Lambda_{jj}^{(k+1)})$ with nonincreasingly (resp.\@ nondecreasingly) ordered diagonal elements,
then $\Lambda = $ \linebreak $\lim_{k\rightarrow\infty}A^{(k)}$  and the diagonal elements of  $\Lambda$ are nonincreasingly (resp.\@ nondecreasingly) ordered.
\item[(ii)] If the algorithm that diagonalizes the pivot submatrix is any standard (i.e., element-wise) globally convergent Jacobi method,
then $\Lambda = \lim_{k\rightarrow\infty}A^{(k)}$.
\end{itemize}
\end{Theorem}

\begin{proof}
The proof has been moved to the Appendix.
$\qquad$ \end{proof}

Theorem~\ref{tm:gl_conv-opp} implies that the global convergence problem  of the block Jacobi method reduces to the convergence of the sequence $S(A^{(k)})$, $k\geq 0$, to zero.

By inspecting the proofs of the results related to the global convergence of the standard  cyclic Jacobi method~\cite{shr+sch-89,han-63} one finds out that they hold for block methods, too. We summarize those results as follows.

\begin{Theorem} \label{tm:SSekv}
If a block Jacobi method converges for some cyclic strategy, then it converges for all strategies that are weak equivalent to it. The block methods defined by equivalent cyclic strategies generate the same matrices after each full cycle and within the same cycle they produce the same sets of orthogonal elementary matrices.
\end{Theorem}

Indeed, the proof for the standard Jacobi method essentially uses the fact that com\-mu\-ting pivot pairs results in commuting the Jacobi rotations. Similarly, the proof for the block method uses the fact that commuting pivot pairs $(i,j)$ and $(p,q)$ imply commuting orthogonal elementary matrices $\mathbf{U}_{ij}$ and $\mathbf{U}_{pq}$. For the convergence of the diagonal elements one should presume conditions like those in Theorem~\ref{tm:gl_conv-opp} for the kernel algorithms. The second part of the theorem holds because it presumes that the block Jacobi method
uses the same kernel algorithm.

Theorem~\ref{tm:SSekv} also holds for the quasi-cyclic methods provided that the care is taken for the blocks that are annihilated more than once within a sweep.

A sufficient condition for the global convergence of the serial standard Jacobi methods is the existence of a strictly positive uniform lower bound for the cosines of the rotation angles (see~\cite{for+hen-60}). For the serial block Jacobi methods a sufficient condition for the global convergence is that the transformation matrices $U_k$ from relation~\eqref{jacobiblokagm} have a strictly positive uniform lower bound for the singular values of the diagonal blocks~\cite{drm-07}. That condition also appears in the global convergence analysis of more general serial Jacobi-type methods~\cite{har-15}. Unitary elementary block matrices which satisfy such property are called \emph{UBC} (uniformly bounded cosine) transformation matrices in~\cite{drm-07}. In the same paper it was shown that for every unitary matrix of order $n$ and every partition $\varsigma =(n_1,n_2)$ of $n$, there exists a permutation matrix $J$ such that for the leading $n_1\times n_1$ block of $\widetilde{U}=UJ$ one has
$$\sigma_{min}(\widetilde{U}_{11})\geq\gamma_{\varsigma}>\widetilde{\gamma}_{n}>0, \quad
\gamma_{\varsigma}=\frac{3}{\sqrt{(4^{n_1}+6n_2-1)(n_2+1)}}, \ \widetilde{\gamma}_{n}=\frac{3\sqrt{2}}{\sqrt{4^n+26}}.
$$
The second inequality, which involves $\widetilde{\gamma}_{n}$, has been proved in~\cite{har-15}. Hence, every unitary elementary block matrix can be made UBC by an appropriate permutation of its nontrivial columns.

In this paper we will use UBC transformation matrices. Therefore, for each $0<\varrho\leq 1$ we introduce the class $\ubce_{\pi} (\varrho)$ of elementary unitary UBC block matrices as follows. The unitary elementary block matrix $\mathbf{U}_{ij}$ from relation (\ref{Uij}) belongs to the class $\ubce_{\pi} (\varrho)$ if
\begin{equation}\label{uv_1}
\sigma_{min}(U_{ii})=\sigma_{min}(U_{jj}) \geq \varrho\gamma_{ij} >\varrho\widetilde{\gamma}_{n_i+n_j} \geq \varrho\widetilde{\gamma}_{n}
\end{equation}
holds, where
\begin{equation}\label{uv_1a}
\gamma_{ij}=\frac{3}{\sqrt{(4^{n_i}+6n_j-1)(n_j+1)}}.
\end{equation}
If $\pi$ is understood, we will write $\ubce (\varrho)$, while if $\varrho =1$, the $\varrho$ will also be omitted from the notation.
In definitions, statements and ordinary text, at every appearance of $\varrho$ we will automatically assume that $0<\varrho\leq 1$.
We will use the same notation when $\mathbf{U}_{ij}$ is real, i.e., orthogonal.

Note that $U_k$ from relation (\ref{jacobiblokagm}) is an orthogonal elementary matrix defined by the pivot pair $(i,j)$ where $i=i(k)$, $j=j(k)$. To make $U_k$ a $\ubce$ transformation, one has to find the permutation $J_k$ and then compute $U_k J_k$. This can be accomplished (see \cite{drm-07}) by performing the QR factorization with column pivoting of $[U_{ii}^{(k)}\ U_{ij}^{(k)}]$ from relation (\ref{blokpivotmatr}). That QR factorization yields $\widehat{J}_k$, which then defines $J_k$ as $J_k=\mathcal{E}(i,j,\widehat{J}_k)$. Then $\widehat{U}_k$ also belongs to the class $\ubce_{\varsigma_{ij}}$, where $\varsigma_{ij}=(n_i,n_j).$
If $(n_i,n_j)$ is understood, $\varsigma_{ij}$ will be omitted.
One ea\-si\-ly verifies that $U_k J_k$ diagonalizes the pivot submatrix $\widehat{A}^{(k)}$ and the similarity transformation with $J_k$ does not change the Frobenius norm of the affected blocks of $U_k^T A^{(k)}U_k$. In addition, one can show that once $S(A^{(k)})$ is sufficiently small, and the diagonal elements affiliated with the same eigenvalue occupy successive positions along the diagonal, the permutations $J_k$ are no longer needed (see~\cite{drm-07}), i.e., $J_k$ can be taken to be the identity.
If $\pi =(1,1,\ldots,1)$, then one can replace $\gamma_{ij}$ by $\sqrt{2}/2$. The uniform bound $\sqrt{2}/2\varrho$ is the one from the known Forsythe-Henrici condition \cite{for+hen-60}.

\begin{Remark}
The parameter $\varrho$ has been introduced for several reasons. First, it simplifies the convergence analysis of the more general iterative process described in Section 5. Second, as will be shown in Sections~3 and 4, the convergence proofs for the symmetric block Jacobi method hold for any $0<\varrho\leq 1$. Finally, for the case $\varrho =1$, the determination of the permutation $J_k$ requires the QR factorization with column pivoting of an $n_i\times (n_i+n_j)$ matrix. Possibly, for some smaller $\varrho$, an appropriate permutation matrix could be obtained at a smaller cost.
\end{Remark}

\subsection{Block Jacobi annihilators and operators}

The Jacobi annihilators and operators have been introduced in~\cite{hen+zim-68} as a tool for proving the global and quadratic convergence of the column-cyclic Jacobi method. Later they have been used for proving the global convergence of some norm-reducing Jacobi-type methods for general matrices~\cite{har-82,har-86}. In~\cite{har-09,har-15,B} they have been generalized to cope with the block Jacobi methods. Here we define a class of the Jacobi annihilators and operators designed precisely for the block Jacobi method for symmetric matrices. They will be referred to as block Jacobi annihilators and operators. This will move us to a more general point of view of the block Jacobi methods, which can be used in the convergence considerations.

First let us introduce some notation.
For an arbitrary $p\times q$ matrix $X$, we define the column vector comprising the columns of $X$,
\[
\textrm{col}(X)=\left[x_{11},x_{21},\ldots,x_{p1},\ldots,x_{1q},\ldots,x_{pq}\right]^T.
\]
Let $\pi=(n_1,\ldots,n_m)$ be a partition of $n$.
Let $\mathbf{S}_{n}$ denote the real vector space of symmetric matrices of order $n$.
Let $A=(A_{rs})\in\mathbf{S}_{n}$ be as in relation (\ref{blokmatrica}). Its block-matrix partition is determined by $\pi$. We define the vector-valued function $\textrm{vec}_{\pi}$ as follows (see~\cite{har-09,har-15}),
\begin{equation}\label{vec}
\textrm{vec}_{\pi}(A) = \left[\begin{array}{c}c_2\\ c_3\\ \vdots \\ c_m\end{array}\right], \qquad \text{where} \
c_j = \left[\begin{array}{c}
            \textrm{col}(A_{1j}) \\
            \textrm{col}(A_{2j}) \\
            \vdots \\
            \textrm{col}(A_{j-1,j}) \\
          \end{array}
        \right], \quad 2\leq j\leq m.
\end{equation}
Then
\begin{equation}\label{vec1}
\textrm{vec}_{\pi}:\mathbf{S}_{n}\rightarrow\mathbb{R}^{K}, \qquad K=N-\sum_{i=1}^m \frac{n_i(n_i-1)}{2}, \quad N=\frac{n(n-1)}{2},
\end{equation}
is a linear operator.

Note that $\textrm{vec}_{\pi}(A)$ contains all off-diagonal elements  from the upper block triangular part of $A$. They are arranged in $\textrm{vec}_{\pi}(A)$ using double column-wise ordering, one with respect to the blocks in $A$, the other with respect to the elements within each block.
The function $\textrm{vec}_{\pi}$ is a surjection, but not an injection. In order to make it bijection, we restrict it to the vector subspace $\mathbf{S}_{0,n}$ of $\mathbf{S}_{n}$, consisting of all matrices from $\mathbf{S}_{n}$ whose diagonal blocks (with respect to the block-matrix partition defined by $\pi$) are zero. Let $\textrm{vec}_{\pi,0}= \textrm{vec}_{\pi}\restrict{\mathbf{S}_{0,n}}$. Obviously, the function $\textrm{vec}_{\pi,0}$ is an invertible linear operator from $\mathbf{S}_{0,n}$ to $\mathbb{R}^{K}$.
In the following text we will often assume that the partition $\pi$ is known and it will be omitted from the notation. However, it will be denoted whenever an additional partition is also considered.

If $a\in\mathbb{R}^{K}$ and $A=\textrm{vec}_{0}^{-1} (a)$ then $A$ is obtained from $a$ using the block-matrix partition defined by $\pi$ and the double column-wise ordering, as is described in relation (\ref{vec}). The diagonal blocks are set to zero and the whole matrix is set to be symmetric. Obviously, $A$ is uniquely determined by $a$.

Beside the linear operators \textrm{vec} and $\textrm{vec}_{0}$, we will make use of
the linear operator $\mathcal{N}_{ij}:\mathbb{R}^{n\times n}\rightarrow\mathbb{R}^{n\times n}$, which also uses the block-matrix partition defined by $\pi$ and sets the pivot submatrix of the argument matrix to zero. When applied to $A\in  \mathbf{S}_{n}$,  $\mathcal{N}_{ij}(A)$ sets the blocks $A_{ij}$, $A_{ji}$, $A_{ii}$ and $A_{jj}$ to zero.

\begin{Definition}
Let $\pi =(n_1,\ldots ,n_m)$ be a partition of $n$, let
\[
\widehat{U}=\left[
    \begin{array}{cc}
    U_{11} & U_{12} \\
    U_{12}^T & U_{22} \\
    \end{array}
    \right] \begin{array}{c} n_i \\ n_j \end{array},
\]
be an orthogonal matrix of order $n_i+n_j$ and let $U=\mathcal{E}(i,j,\widehat{U})$ be the corresponding elementary block matrix. The transformation $\mathcal{R}_{ij}(\widehat{U})$ determined by
\[
\mathcal{R}_{ij}(\widehat{U})(vec(A)) = vec(\mathcal{N}_{ij}(U^TAU)), \quad A\in \mathbf{S}_{n},
\]
is called the $ij$-\emph{block Jacobi annihilator}. For each pair $1\leq i<j\leq m$,
\[
\Il_{ij} =\left\{\mathcal{R}_{ij}(\widehat{U}) \ \big{|} \ \widehat{U} \text{ is orthogonal matrix of order } n_i+n_j\right\}
\]
is the $ij$-class of the block Jacobi annihilators. If all $\widehat{U}$ are restricted to the class
$\ubce_{\varsigma_{ij}}(\varrho)$, then the resulting $ij$-class is denoted by $\Il_{ij}^{\ubce_{\pi}(\varrho)}$.
\end{Definition}

Given $i$, $j$, $\widehat{U}$, the following algorithm computes the vector $a'=\mathcal{R}_{ij}(\widehat{U})a$ for $a\in\mathbb{R}^K$. It is based on the formula $\mathcal{R}_{ij}(\widehat{U})a = \textrm{vec}(\mathcal{N}_{ij}(U^T\textrm{vec}_{0}^{-1}(a) U))$, which can be taken as an equivalent definition of $\mathcal{R}_{ij}(\widehat{U})$.

\begin{Algorithm}\label{bjanagm}
\vspace{0.5ex}\hrule\vspace{0.5ex}
\emph{Computing} $\mathcal{R}_{ij}(\widehat{U})a$
\vspace{0.5ex}\hrule
\vspace{-2ex}
\begin{multicols}{2}
\begin{algorithmic}
\State $a\in\mathbb{R}^K$
\State  $A=\textrm{vec}_0^{-1}(a)$
\For {$r=1,\ldots,m$}
\State $A'_{ri}=A_{ri}U_{ii}+A_{rj}U_{ji}$
\State $A'_{rj}=A_{ri}U_{ij}+A_{rj}U_{jj}$
\EndFor
\For {$r=1,\ldots,m$}
\State $A'_{ir}=U_{ii}^TA_{ir}+U_{ji}^TA_{jr}$
\State $A'_{jr}=U_{ij}^TA_{ir}+U_{jj}^TA_{jr}$
\EndFor
\State $A'_{ij}=0$, $A'_{ji}=0$, $A'_{ii}=0$, $A'_{jj}=0$
\State $a'=\textrm{vec}(A')$ \\
\textsf{\small \% \ an arbitrary vector} \\
\textsf{\small \% \ this invokes the module which computes} \\
\textsf{\small \% \ the symmetric matrix $A=\textrm{vec}_0^{-1}(a)$} \\
\\
\\
\textsf{\small \% \ this part of code computes} \\
\textsf{\small \% \ $A'=U^TAU$, $U=\mathcal{E}(i,j,\widehat{U})$} \\
\textsf{\small \% \ $\widehat{U}$ is partitioned as in relation (\ref{piv_mat_E})} \\
\\
\\
\textsf{\small \% \ this part updates $A'$, $A'\leftarrow\mathcal{N}_{ij}(A')$} \\
\textsf{\small \% \ the module which computes $\textrm{vec}(A')$}
\end{algorithmic}
\end{multicols}
\hrule
\end{Algorithm}

\bigskip

Matrix $A'$ from Algorithm~\ref{bjanagm} has the same partition as $A$. Note that the mapping $a\mapsto a'$ is a composition of linear transformations. Therefore, given a basis in $\mathbb{R}^{K}$, $\mathcal{R}_{ij}(\widehat{U})$ can be represented by some square matrix of order $K$. We will choose the canonic basis $(e)$, which consists of the columns of $I_{\!K}$, and denote the obtained matrix by the bold symbol. Hence,
\begin{equation}\label{Jannih}
\mathcal{R}_{ij}(\widehat{U})(a)=\bm{\mathcal{R}}_{ij}(\widehat{U})a, \quad a\in\mathbb{R}^{K}.
\end{equation}
We will call the matrix $\bm{\mathcal{R}}_{ij}(\widehat{U})$ by the same name, the $ij$-block Jacobi annihilator, and the appropriate class of matrices will be denoted by
\[
\Ilb_{ij} =\left\{\bm{\mathcal{R}}_{ij}(\widehat{U}) \ \big{|} \ \widehat{U} \ \text{is an orthogonal matrix of order} \ n_i+n_j \right\}.
\]
If all $\widehat{U}$ are restricted to $\ubce(\varrho)$, the obtained class of block Jacobi annihilators is denoted by $\Ilb_{ij}^{\ubce(\varrho)}$. In the sequel, every mention of the block Jacobi annihilator will refer to the matrix $\bm{\mathcal{R}}_{ij}(\widehat{U})$ from relation (\ref{Jannih}).

The following theorem reveals the structure of a block Jacobi annihilator. It is a  simplification of~\cite[Theorem~2.1]{har-09} and its proof can be found in~\cite{B}.
The theorem utilizes the function $\tau(i,j)=(j-1)(j-2)/2+i$, $1\leq i<j\leq m$,
and the Kronecker product of matrices. The vectors of length $K$ and the block Jacobi annihilators of order $K$ carry the block-partition determined by $(n_1n_2,n_1n_3,n_2n_3,\ldots,n_{m-1}n_m)$.
The spectral norm is denoted by $\|\cdot\|_2$.

\begin{Theorem}[\cite{har-09,B}]\label{tm:bl_jac_ann}
Let $\pi=(n_1,\ldots,n_m)$ be the partition of $n$ and let $K$, $N$ be integers defined in relation (\ref{vec1}). Let $(i,j)\in \mathcal{P}_m$, $\bm{\mathcal{R}}\in\Ilb_{ij}$, $\bm{\mathcal{R}}=\bm{\mathcal{R}}(\widehat{U})$, where $\widehat{U}$ is an orthogonal matrix of order $n_i+n_j$. Then $\bm{\mathcal{R}}$ differs from the identity matrix $I_{K}$ in exactly $m-1$ principal submatrices, which are given by the following relations:
\[
\mathcal{R}_{\tau(i,j),\tau(i,j)}=0,
\]
{\small
\begin{align*}
 \left[
    \begin{array}{cc}
      \mathcal{R}_{\tau(r,i),\tau(r,i)} & \mathcal{R}_{\tau(r,i),\tau(r,j)} \\
      \mathcal{R}_{\tau(r,j),\tau(r,i)} & \mathcal{R}_{\tau(r,j),\tau(r,j)}
    \end{array}
  \right]=& \left[
     \begin{array}{cc}
       U_{ii}^T\otimes I_{n_r} & U_{ji}^T\otimes I_{n_r} \\
       U_{ij}^T\otimes I_{n_r} & U_{jj}^T\otimes I_{n_r}
     \end{array}
   \right], \quad 1\leq r\leq i-1, \\
 \left[
    \begin{array}{cc}
      \mathcal{R}_{\tau(i,r),\tau(i,r)} & \mathcal{R}_{\tau(i,r),\tau(r,j)} \\
      \mathcal{R}_{\tau(r,j),\tau(i,r)} & \mathcal{R}_{\tau(r,j),\tau(r,j)}
    \end{array}
  \right]= & \left[
     \begin{array}{cc}
       I_{n_r}\otimes U_{ii}^T & S(U_{ji}^T\otimes I_{n_r}) \\
       \widetilde{S}(I_{n_r}\otimes U_{ij}^T) & U_{jj}^T\otimes I_{n_r}
     \end{array}
   \right], \ i\!+\!1\leq r\leq j\mbox{--}1,\\
   \left[
    \begin{array}{cc}
      \mathcal{R}_{\tau(i,r),\tau(i,r)} & \mathcal{R}_{\tau(i,r),\tau(j,r)} \\
      \mathcal{R}_{\tau(j,r),\tau(i,r)} & \mathcal{R}_{\tau(j,r),\tau(j,r)}
    \end{array}
  \right]=& \left[
     \begin{array}{cc}
       I_{n_r}\otimes U_{ii}^T & I_{n_r}\otimes U_{ji}^T \\
       I_{n_r}\otimes U_{ij}^T & I_{n_r}\otimes U_{jj}^T
     \end{array}
   \right], \quad j+1\leq r\leq m,
\end{align*}}
where
{\small
\begin{align*}
S & =\left[\!\!
      \begin{array}{c}
        I_{n_i}\otimes e_1^T \\
        \vdots \\
        I_{n_i}\otimes e_{n_r}^T \\
      \end{array}\!\! \right]=\left[\!\!
              \begin{array}{ccc}
                I_{n_r}\otimes \bar{e}_1 & \ldots & I_{n_r}\otimes \bar{e}_{n_i} \\
              \end{array}\!\!\right], \\
\widetilde{S} & =\left[\!\!
      \begin{array}{c}
        I_{n_r}\otimes \widetilde{e}_1^T \\
        \vdots \\
        I_{n_r}\otimes \widetilde{e}_{n_j}^T \\
      \end{array}\!\!\right]
      =\left[\!\!\begin{array}{ccc}
                I_{n_j}\otimes e_1 & \ldots & I_{n_j}\otimes e_{n_r} \\
              \end{array}\!\!\right].
\end{align*}}
Here, $e_i$, $\bar{e}_i$ and $\widetilde{e}_i$ denote the $i$th column of $I_{n_r}$, $I_{n_i}$ and $I_{n_j}$, respectively.

The matrix $\bm{\mathcal{R}}$ satisfies $\|\bm{\mathcal{R}}\|_2=1$, except in the case $m=2$, $(i,j)=(1,2)$, when $\bm{\mathcal{R}}=0$.
\end{Theorem}

\begin{Example}
Let $A\in\mathbb{R}^{8\times8}$, $\pi=(2,2,2,2)$, $i=1$, $j=2$. Then $K=24$ and
{\setlength{\arraycolsep}{0.6ex}
\[
\bm{\mathcal{R}} = \bm{\mathcal{R}}_{12}(\widehat{U})=
\mbox{\scriptsize $\displaystyle \left[
                        \begin{array}{cccc|cc|cc|cc|cc|cccc}
                          0 &  &  &  &  &  &  &  &  &  &  &  &  &  &  &  \\
                           & 0 &  &  &  &  &  &  &  &  &  &  &  &  &  &  \\
                           &  & 0 &  &  &  &  &  &  &  &  &  &  &  &  &  \\
                           &  &  & 0 &  &  &  &  &  &  &  &  &  &  &  &  \\ \hline
                           &  &  &  &  &  &  &  &  &  &  &  &  &  &  &  \\[-1ex]
                           &  &  &  & U_{11}^T &  & U_{21}^T &  &  &  &  &  &  &  &  &  \\[1ex]
                           &  &  &  &  & U_{11}^T &  & U_{21}^T &  &  &  &  &  &  &  &  \\[1ex] \hline
                           &  &  &  &  &  &  &  &  &  &  &  &  &  &  &  \\[-1ex]
                           &  &  &  & U_{12}^T &  & U_{22}^T &  &  &  &  &  &  &  &  &  \\[1ex]
                           &  &  &  &  & U_{12}^T &  & U_{22}^T &  &  &  &  &  &  &  &  \\[1ex] \hline
                           &  &  &  &  &  &  &  &  &  &  &  &  &  &  &  \\[-1ex]
                           &  &  &  &  &  &  &  & U_{11}^T &  & U_{21}^T &  &  &  &  &  \\[1ex]
                           &  &  &  &  &  &  &  &  & U_{11}^T &  & U_{21}^T &  &  &  &  \\[1ex] \hline
                           &  &  &  &  &  &  &  &  &  &  &  &  &  &  &  \\[-1ex]
                           &  &  &  &  &  &  &  & U_{12}^T &  & U_{22}^T &  &  &  &  &  \\[1ex]
                           &  &  &  &  &  &  &  &  & U_{12}^T &  & U_{22}^T &  &  &  &  \\[1ex]  \hline
                           &  &  &  &  &  &  &  &  &  &  &  & 1 &  &  & \\
                           &  &  &  &  &  &  &  &  &  &  &  &  & 1 &  & \\
                           &  &  &  &  &  &  &  &  &  &  &  &  &  & 1 &  \\
                           &  &  &  &  &  &  &  &  &  &  &  &  &  &  & 1 \\
                         \end{array}\right]$},
\]}
where $U_{11}$, $U_{12}$, $U_{21}$, $U_{22}$ are the blocks of order $2$ of $\widehat{U}\in \mathbb{R}^{4\times4}$ and $\widehat{U}$ is orthogonal.
\end{Example}

We see that $\bm{\mathcal{R}}$ is, up to similarity transformation with permutation, a direct sum of an orthogonal matrix and the zero matrix. Therefore, $\|\bm{\mathcal{R}}\|_2=1$, except in the case $m=2$, $(i,j)=(1,2)$ when it is the zero matrix.

\begin{Corollary}\label{tm:cor_bjann}
Let $\pi$, $(i,j)\in \mathcal{P}_m$ and $\bm{\mathcal{R}}\in\Ilb_{ij}$ be as in Theorem~\ref{tm:bl_jac_ann}. Then $\bm{\mathcal{R}}^T\in\Ilb_{ij}$.
Moreover, if $\bm{\mathcal{R}}\in\Ilb_{ij}^{\ubce(\varrho)}$, then $\bm{\mathcal{R}}^T\in\Ilb_{ij}^{\ubce(\varrho)}$.
\end{Corollary}

\begin{proof}
The proof has been moved to the Appendix.
$\qquad$ \end{proof}

The block Jacobi annihilators are used to define the block Jacobi operators, which make up our tool for proving the global convergence of the block Jacobi methods.

\begin{Definition}\label{def:jop}
Let $\pi=(n_1,\ldots,n_m)$ be a partition of $n$ and let
\[
\mathcal{O} = (i_0,j_0),(i_1,j_1),\ldots,(i_{T-1},j_{T-1})\in \Ol(\mathcal{P}_m), \quad T\geq M=\frac{m(m-1)}{2}.
\]
Then
\[
\Jl_{\!\!\! \mathcal{O}} = \{\mathcal{J} \ \big{|} \ \mathcal{J}=\bm{\mathcal{R}}_{i_{T-1}j_{T-1}}\ldots\bm{\mathcal{R}}_{i_1j_1}\bm{\mathcal{R}}_{i_0j_0},
\ \bm{\mathcal{R}}_{i_kj_k}\in\Ilb_{i_k j_k}, \ 0\leq k\leq T-1\}
\]
is called the class of block Jacobi operators associated with the sequence $\mathcal{O}$. The matrices $\mathcal{J}$ of order $K$ from $\Jl_{\!\!\! \mathcal{O}}$ are the block Jacobi operators.
If each $\Ilb_{ij}$ in $\Jl_{\!\!\! \mathcal{O}}$ is replaced by $\Ilb_{ij}^{\ubce(\varrho)}$, then the notation $\Jl_{\!\!\! \mathcal{O}}^{\ubce(\varrho)}$ will be used.
\end{Definition}

An element $\mathcal{J}\in\Jl_{\!\!\! \mathcal{O}}$ will sometimes be written $\mathcal{J}_{ \mathcal{O}}$.
The following lemma reveals some properties of the block Jacobi operators. The spectral radius of a square matrix $X$ is denoted by $\text{spr}(X)$.

\begin{Lemma}[\cite{har-15,B}]\label{tm:lemma_spran}
Let $\pi=(n_1,\ldots,n_m)$ be a partition of $n$, $\mathcal{O},\mathcal{O}'\in\Ol(\mathcal{P}_m)$ and $\mathcal{O}'\stackrel{\mathsf{w}}{\sim}\mathcal{O}$. Take
$\mathcal{J}_{\mathcal{O}}\in \Jl_{\!\!\! \mathcal{O}}$ and let $\mathcal{J}_{\mathcal{O}'}$ be comprised of the same block Jacobi annihilators as $\mathcal{J}_{\mathcal{O}}$.
Then $\text{spr}(\mathcal{J}_{\mathcal{O}})=\text{spr}(\mathcal{J}_{\mathcal{O}'})$.
If $\mathcal{O}\sim\mathcal{O}'$, then $\mathcal{J}_{\mathcal{O}}=\mathcal{J}_{\mathcal{O}'}$.
\end{Lemma}

\begin{proof}
The proof is the same as the proof of~\cite[Lemma~4.4]{har-15}.
$\qquad$ \end{proof}

If the spectral norm is used instead of the spectral radius, then we have the following result.

\begin{Proposition}\label{tm:prop_norman}
Let $\pi=(n_1,\ldots,n_m)$ be a partition of $n$, $\mathcal{O},\mathcal{O}'\in\Ol(\mathcal{P}_m)$ and $\mathcal{O}'\stackrel{\mathsf{w}}{\sim}\mathcal{O}$. Let $\mathcal{O}$ and $\mathcal{O}'$ be linked by the chain $\mathcal{O}=\mathcal{O}_0,\mathcal{O}_1,\ldots,\mathcal{O}_r=\mathcal{O}'$, as in relation (\ref{weak_eq}). Suppose that in the chain there are exactly $d$ pairs of neighboring terms that are shift equivalent. If
\[
\|\mathcal{J}\|_2\leq\mu_{\pi,\varrho} \quad \text{for all} \ \mathcal{J} \in \Jl_{\!\!\! \mathcal{O}}^{\ubce(\varrho)},
\]
then for any $d+1$ block Jacobi operators from $\Jl_{\!\!\! \mathcal{O}'}^{\ubce(\varrho)}$ one has
\[
\|\mathcal{J}_1' \mathcal{J}_2'\cdots \mathcal{J}_{d+1}'\|_2 \leq \mu_{\pi,\varrho},\qquad \mathcal{J}_1',\ldots,\mathcal{J}_{d+1}' \in \Jl_{\!\!\! \mathcal{O}'}^{\ubce(\varrho)}.
\]
The constant $\mu_{\pi,\varrho}$ may depend only on $\pi$ and $\varrho$.
\end{Proposition}

\begin{proof}
The proof is the same as the proof of \cite[Lemma~4.8(ii)]{har-15}. The role of the set $\Psi_{\pi}$ from~\cite[Lemma~4.8(ii)]{har-15} is played by the set $\bigcup_{i<j}\Ilb_{ij}^{\ubce(\varrho)}$.
$\qquad$ \end{proof}

\begin{Proposition}\label{tm:prop_invstr}
Let $\pi=(n_1,\ldots,n_m)$ be a partition of $n$, $\mathcal{O}\in\Ol(\mathcal{P}_m)$ and suppose
$\|\mathcal{J}\|_2\leq\mu_{\pi,\varrho}$ for all $\mathcal{J}\in\Jl_{\!\!\! \mathcal{O}}^{\ubce(\varrho)}$,
where $\mu_{\varrho}$ depends on $\pi$ and $\varrho$. Then
\[
\|\widetilde{\mathcal{J}}\|_2\leq\mu_{\pi,\varrho} \quad \text{for all} \ \widetilde{\mathcal{J}} \in \Jl_{\!\!\! \mathcal{O}^{\leftarrow}}^{\ubce(\varrho)}.
\]
The assertion holds provided that in both appearances the spectral norm is replaced by the spectral radius.
\end{Proposition}

\begin{proof}
Suppose $\mathcal{O} = (i_0,j_0),(i_1,j_1),\ldots,(i_{T-1},j_{T-1})\in\Ol(\mathcal{P}_m)$ and let
$\widetilde{\mathcal{J}}\in$ $\Jl_{\!\!\! \mathcal{O}^{\leftarrow}}^{\ubce(\varrho)}$ be arbitrary. Then
\[
\widetilde{\mathcal{J}}=\bm{\mathcal{R}}_{i(0),j(0)}(\widehat{U}_{0})\,\bm{\mathcal{R}}_{i(1),j(1)}(\widehat{U}_{1})\cdots \bm{\mathcal{R}}_{i(T-1),j(T-1)}(\widehat{U}_{T-1}),
\]
for some orthogonal $\ubce(\varrho)$ matrices $\widehat{U}_{k}$, $0\leq k\leq T-1$, of appropriate sizes.

If we show that ${\widetilde{\mathcal{J}}}^T\in \Jl_{\!\!\! \mathcal{O}}^{\ubce(\varrho)}$, first claim will follow from $\|\widetilde{\mathcal{J}}\|_2=\|\widetilde{\mathcal{J}}^T\|_2\leq\mu_{\pi,\varrho}$,
while the part about the spectral radius will be a consequence of $\text{spr}(\widetilde{\mathcal{J}})=\text{spr}(\widetilde{\mathcal{J}}^T) \leq \mu_{\pi,\varrho}$.
Note that
\begin{align*}
\widetilde{\mathcal{J}}^T & = \left[
\bm{\mathcal{R}}_{i(0),j(0)}(\widehat{U}_{0})\bm{\mathcal{R}}_{i(1),j(1)}(\widehat{U}_{1}) \cdots \bm{\mathcal{R}}_{i(T-1),j(T-1)}(\widehat{U}_{T-1})
\right]^T \\
& = [\bm{\mathcal{R}}_{i(T-1),j(T-1)}(\widehat{U}_{T-1})]^T \cdots [\bm{\mathcal{R}}_{i(1),j(1)}(\widehat{U}_1)]^T\cdot [\bm{\mathcal{R}}_{i(0),j(0)}(\widehat{U}_0)]^T.
\end{align*}
By Corollary~\ref{tm:cor_bjann} we know that
$[\bm{\mathcal{R}}_{ij}(\widehat{U}_k)]^T\in \Ilb_{ij}^{\ubce(\varrho)}$, hence $\widetilde{\mathcal{J}}^T\in \Jl_{\!\!\! \mathcal{O}}^{\ubce(\varrho)}$.
$\qquad$ \end{proof}

\subsubsection{Permutation equivalence and the block Jacobi operators}

Here we derive a similar result for the block Jacobi operators $\mathcal{J}_{\mathcal{O}}$ and $\mathcal{J}_{\mathcal{O}(\mathsf{p})}$.

\begin{Theorem}\label{tm:prop_permstr}
Let  $\pi=(n_1,\ldots,n_m)$ be a partition of $n$ and take $\mathcal{O}\in\Ol(\mathcal{P}_m)$.
Let $\mathsf{p}$ be a permutation of $\mathcal{S}_m$ and set $\widetilde{\mathcal{O}}=\mathcal{O}(\mathsf{p})$.
\begin{itemize}
\item[(i)] If $\displaystyle \|\mathcal{J}\|_2\leq\mu_{\pi,\varrho}$ for all $\displaystyle \mathcal{J}\in\Jl_{\!\!\! \mathcal{O}}^{\ubce_{\pi}(\varrho)}$,
where $\mu_{\pi,\varrho}$ only depends on $\pi$ and $\varrho$, then
\[
\|\widetilde{\mathcal{J}}\|_2\leq\mu_{\pi,\varrho} \quad \text{for any} \ \widetilde{\mathcal{J}}\in\Jl_{\!\!\! \widetilde{\mathcal{O}}}^{\ubce_{\pi_{\mathsf{p}}} (\varrho)}.
\]
\item[(ii)] If $\displaystyle \|\mathcal{J}_1 \mathcal{J}_2\cdots\mathcal{J}_{d+1}\|_2\leq\mu_{\pi,\varrho}$, for all
$\displaystyle \mathcal{J}_1,\ldots,\mathcal{J}_{d+1}\in\Jl_{\!\!\! \mathcal{O}}^{\ubce_{\pi} (\varrho)}$, where $\mu_{\pi,\varrho}$ only depends on $\pi$ and $\varrho$, then
\[
\|\widetilde{\mathcal{J}}_1 \widetilde{\mathcal{J}}_2\cdots \widetilde{\mathcal{J}}_{d+1}\|_2 \leq \mu_{\pi,\varrho} \quad \text{for any} \ \widetilde{\mathcal{J}}_1,\ldots,\widetilde{\mathcal{J}}_{d+1}\in \Jl_{\!\!\! \widetilde{\mathcal{O}}}^{\ubce_{\pi_{\mathsf{p}}} (\varrho)}.
\]
\end{itemize}
The assertions also hold provided that every appearance of the spectral norm is replaced by the spectral radius.
\end{Theorem}

\begin{proof} $(i)$ Let $\mathcal{O} = (i_0,j_0),(i_1,j_1),\ldots,(i_{T-1},j_{T-1})$, so that $\widetilde{\mathcal{O}}=(\mathsf{p}(i_0),\mathsf{p}(j_0)),(\mathsf{p}(i_1),\mathsf{p}(j_1)),\ldots,$ \linebreak $(\mathsf{p}(i_{T-1}),\mathsf{p}(j_{T-1}))$. Let $\widetilde{\mathcal{J}}\in\Jl_{\!\!\! \widetilde{\mathcal{O}}}^{\ubce_{\pi_{\mathsf{p}}} (\varrho)}$ be arbitrary. Then
\begin{equation}\label{2.19aa}
\widetilde{\mathcal{J}} = \bm{\mathcal{R}}_{\mathsf{p}(i_{T-1}),\mathsf{p}(j_{T-1})}(\widehat{U}_{T-1})\cdots \bm{\mathcal{R}}_{\mathsf{p}(i_{1}),\mathsf{p}(j_{1})}(\widehat{U}_{1})
\bm{\mathcal{R}}_{\mathsf{p}(i_{0}),\mathsf{p}(j_{0})}(\widehat{U}_{0}),
\end{equation}
where $\widehat{U}_{k}$ is an orthogonal $\ubce_{\pi_{\mathsf{p}}} (\varrho)$ matrix of order $n_{\mathsf{p}(i_{k})}+n_{\mathsf{p}(j_{k})}$ for each $0\leq k\leq T-1$.

Let $a\in \mathbb{R}^{K}$ be an arbitrary nonzero vector, and consider the computation of $a' = \widetilde{\mathcal{J}}a$.
Using Algorithm~\ref{bjanagm}, the vector $a'$ can be obtained by the following procedure:
\begin{itemize}
\item Compute the symmetric matrix $A^{(0)}\equiv A=\ve_{\pi_{\mathsf{p}},0}^{-1}(a)$.
\item Recursively compute: $A^{(k+1)}=\mathcal{N}_{\mathsf{p}(i)\mathsf{p}(j)}(U_{k}^T A^{(k)} U_{k})$, $k=0,1,\ldots,T-1$.
\item Compute the vector $a'= \ve_{\pi_{\mathsf{p}}} (A^{(T)})$.
\end{itemize}
Here $U_{k}=\mathcal{E}(\mathsf{p}(i_k),\mathsf{p}(j_k),\widehat{U}_k)$, $0\leq k\leq T-1$, and the matrices $A$ and $ A^{(k)}$ carry the matrix block-partition defined by $\pi_{\mathsf{p}}$.
Let $P$ be the matrix from relation (\ref{UTAU}), which is defined by $Pe_t=e_{p(t)}$, $1\leq t\leq n$, where $p$ is defined as in Section~\ref{sec:2.2.1}, and $\mathsf{p}$ is from the statement of this lemma. Let $X$ be a square matrix of order $n$, partitioned in accordance with $\pi_{\mathsf{p}}$. Then for any $1\leq s,t\leq m$ the transformation $X\mapsto P^T X P$ changes the partition from $\pi_{\mathsf{p}}$ to $\pi$ and moves the block $X_{\mathsf{p}(s)\mathsf{p}(t)}$ to the $(s,t)$ position.
Therefore, we have
\begin{align}
P^T A^{(k+1)}P & = P^T \mathcal{N}_{\mathsf{p}(i)\mathsf{p}(j)}\left(U_{k}^T A^{(k)} U_{k}\right)P
= \mathcal{N}_{ij}\left( P^T U_{k}^T A^{(k)} U_{k} P\right) \label{2.19a} \\
& = \mathcal{N}_{ij}\left( (P^T U_{k} P)^T (P^T A^{(k)}P) (P^TU_{k} P)\right), \quad
k=0,1,\ldots,T-1. \nonumber
\end{align}
If we set $\bar{A}^{(k)} = P^T A^{(k)}P$, $\bar{A} = P^T AP$ and $\bar{U}^{(k)} = P^T U^{(k)}P$, then recurrence (\ref{2.19a}) takes the form
\begin{equation}\label{2.19a1}
\bar{A}^{(k+1)} = \mathcal{N}_{ij}\left([\bar{U}^{(k)}]^T \bar{A}^{(k)}\bar{U}^{(k)}\right), \quad\ k=0,1,\ldots,T-1; \quad \bar{A}^{(0)} = \bar{A}.
\end{equation}
Obviously, we have $\ve_{\pi} (\bar{A}) = \ve_{\pi} (P^T AP) =\mathbf{P}a$ for some permutation matrix $\mathbf{P}$ of order $K$. Applying the $\ve$ function to relation (\ref{2.19a1}), one obtains
\begin{equation}\label{2.19a2}
\bar{a}^{(k+1)} = \bm{\mathcal{R}}_{i(k),j(k)}(\widehat{\bar{U}}_{k}) \bar{a}^{(k)}, \quad k=0,1,\ldots,T-1; \quad \bar{a}^{(0)} = \mathbf{P}a,
\end{equation}
where $\bm{\mathcal{R}}_{ij}(\widehat{\bar{U}}_{k})\in\Ilb_{ij}^{\ubce_{\pi} (\varrho)}$.
This process is associated with the sequence $\mathcal{O}$ and results in the final form
\[
\mathbf{P}a' = \bar{\mathcal{J}}_{\mathcal{O}}\mathbf{P}a,\quad \bar{\mathcal{J}}_{\mathcal{O}}\in \Jl_{\!\!\! \mathcal{O}}^{\ubce_{\pi} (\varrho)}.
\]
Because $a$ is an arbitrary vector, we have $\widetilde{\mathcal{J}} = \mathbf{P}^T \bar{\mathcal{J}}_{\mathcal{O}}\mathbf{P}$. This implies $\|\widetilde{\mathcal{J}}\|_2$ $=$
$\|\bar{\mathcal{J}}_{\mathcal{O}}\|_2 \leq \mu_{\pi,\varrho}$.
Since, $\widetilde{\mathcal{J}}$ and $\bar{\mathcal{J}}_{\mathcal{O}}$ are similar, their spectral radius is the same.

$(ii)$ The proof is similar to the proof of $(i)$. We start our consideration with an arbitrary nonzero vector $a\in\mathbb{R}^{K}$ and consider the computation of
$a'=\widetilde{\mathcal{J}}_1 \widetilde{\mathcal{J}}_2\cdots \widetilde{\mathcal{J}}_{d+1}a$. Since we have $d+1$ block Jacobi operators, we will use altogether $(d+1)T$ block Jacobi annihilators. All we have to change in the proof of (i) is the the range of the index $k$ in the relations (\ref{2.19aa}) -- (\ref{2.19a2}): instead of $T-1$ its largest value will be $(d+1)T-1$. At the end we will have
$\displaystyle \widetilde{\mathcal{J}}_1 \widetilde{\mathcal{J}}_2\cdots \widetilde{\mathcal{J}}_{d+1}=\mathbf{P}^T\mathcal{J}_1 \mathcal{J}_2\cdots \mathcal{J}_{d+1}\mathbf{P}$
and the conclusion will follow.
$\qquad$ \end{proof}

\section{Generalized serial strategies}\label{cyclic}

The aim of this section is to significantly enlarge the class of the known ``convergent'' cyclic pivot strategies, namely the serial ones and those that are weak equivalent to them (the so-called weak wavefront strategies~\cite{shr+sch-89}).
We study several classes of cyclic pivot strategies, which are generalizations of the serial ones. The first (resp.\@ second) of those classes is defined by the set $\mathcal{B}_c^{(m)}$ (resp.\@ $\mathcal{B}_r^{(m)}$) of pivot orderings which arise from column-wise (resp.\@ row-wise) orderings of $\mathcal{P}_m$. The other two are defined by the first two using reverse orderings.
Once the global convergence of the block Jacobi method under these strategies is proved, one can easily expand the obtained set of pivot strategies using the theory of equivalent strategies.

\subsection{The class $\mathcal{B}_c^{(m)}$}

We start with the class of cyclic strategies that choose $(1,2)$-block as the first pivot block, then choose all blocks from the second block-column in some order, etc. At the last stage they choose all blocks from the last block-column in some order. For the precise definition of that class, we denote the set of all permutations of the set $\{l_1,l_1+1,\ldots,l_2\}$ by $\Pi^{(l_1,l_2)}$. Let
\begin{align}
\label{b1def} \mathcal{B}_c^{(m)} = \big{\{} \mathcal{O}\in\mathcal{\Ol}(\mathcal{P}_m) \ \big{|} & \ \mathcal{O}= (1,2),(\tau_{3}(1),3),(\tau_{3}(2),3),\ldots,(\tau_{m}(1),m),\ldots, \\
& (\tau_{m}(m-1),m), \quad \tau_{j}\in\Pi^{(1,j-1)}, \ 3\leq j\leq m \big{\}}. \nonumber
\end{align}
The set $\mathcal{B}_c^{(m)}$ is a part of the class of \emph{column-wise orderings with permutations} of the set $\mathcal{P}_m$, which will be described in Definition~\ref{tm:def3.5}.
A typical ordering $\mathcal{O}\in\mathcal{B}_c^{(6)}$ is represented by $\mathsf{M}_{\mathcal{O}}$ below. The second matrix
$\mathsf{M}_{\widetilde{\mathcal{O}}}$ is defined by  some $\widetilde{\mathcal{O}}\stackrel{\mathsf{w}}{\sim}\mathcal{O}$. Its purpose is to see how far from the ``serial structure'' this equivalence can push $\mathcal{O}$.
{\small
\[\mbox{\normalsize $\mathsf{M}_{\mathcal{O}}$} =\left[
              \begin{array}{cccccc}
                * & 0 & 2 & 4 & 9 & 12 \\
                0 & * & 1 & 5 & 8 & 10 \\
                2 & 1 & * & 3 & 7 & 13 \\
                4 & 5 & 3 & * & 6 & 11 \\
                9 & 8 & 7 & 6 & * & 14 \\
                12 & 10 & 13 & 11 & 14 & * \\
              \end{array}\right], \qquad
\mbox{\normalsize $\mathsf{M}_{\widetilde{\mathcal{O}}}$}= \left[\begin{array}{cccccc}
                * & 7 & 9 & 0 & 2 & 5 \\
                7 & * & 10 & 13 & 14 & 6 \\
                9 & 10 & * & 11 & 12 & 8 \\
                0 & 13 & 11 & * & 1 & 4 \\
                2 & 14 & 12 & 1 & * & 3 \\
                5 & 6 & 8 & 4 & 3 & * \\
              \end{array}\right].
\]}
From $\mathsf{M}_{\mathcal{O}}$ we can see that the permutation $\tau_j$ from~\eqref{b1def} is linked with the $j$th block-column of the matrix.
The next theorem proves the global convergence of the block Jacobi method under the cyclic strategies $I_{\mathcal{O}}$ defined by the orderings $\mathcal{O}\in\mathcal{B}_c^{(m)}$.

\begin{Theorem}\label{tm:b1}
Let $\pi=(n_1,\ldots ,n_m)$ be a partition of $n$ and let $\mathcal{O}\in\mathcal{B}_c^{(m)}$. Let $A\in\mathbf{S}_{n}$ be partitioned as in relation \eqref{blokmatrica}. Suppose that $A'$ is obtained from $A$ by applying one sweep
of the cyclic block Jacobi method defined by the strategy $I_{\mathcal{O}}$.  If all transformation matrices are from the class $\ubce(\varrho)$, then there are constants $\eta_{\pi,\varrho}$ (depending only on $\pi$ and $\varrho$) and $\widetilde{\eta}_{n,\varrho}$ (depending only on $n$ and $\varrho$) such that
\[
S^2(A')\leq \eta_{\pi,\varrho} S^2(A), \quad 0\leq \eta_{\pi,\varrho} <\widetilde{\eta}_{n,\varrho}<1.
\]
\end{Theorem}

\begin{proof}
The proof is lengthy and has been moved to the Appendix.
$\qquad$ \end{proof}

We have to explain why we use the two bounds satisfying $\mu_{\pi,\varrho}<\widetilde{\mu}_{n,\varrho}<1$.
Recall that each block Jacobi method is defined by some partition $\pi$ of $n$. Different partitions define different block Jacobi methods, even in the case when the pivot orderings are the same. The second bound $\widetilde{\mu}_{n,\varrho}$ can be used in the global convergence statements for the block Jacobi method, when the order $n$ of the initial matrix is known, while about the pivot ordering it is only known that it belongs to the set $\cup_{3\leq m\leq n}\mathcal{B}_c^{(m)}$. It means that, for a given $m$, the convergence result holds for the block Jacobi method defined by any $\pi$ such that $n_1+\cdots +n_m=n$.

Combining Theorem~\ref{tm:b1} with Theorem~\ref{tm:SSekv}, we see that we can enlarge the class of
``convergent orderings'' from $\mathcal{B}_c^{(m)}$ to the class of all orderings that are weak equivalent to orderings from $\mathcal{B}_c^{(m)}$. Thus, the ordering $\widetilde{\mathcal{O}}$ linked with the above matrix $\mathsf{M}_{\widetilde{\mathcal{O}}}$ is also convergent.

The next result is a slight generalization of Theorem~\ref{tm:b1} and it deals with the block Jacobi operators. The role of the block Jacobi operators will be explained in Section~5, especially by inspecting the proof of Theorem~\ref{opcitm}.

\begin{Theorem}\label{tm:b1jop}
Let $\pi=(n_1,\ldots,n_m)$ be a partition of $n$. Take $\mathcal{O}\in\mathcal{B}_c^{(m)}$ and let $\mathcal{J}\in\Jl_{\!\!\! \mathcal{O}}^{\ubce(\varrho )}$ be a block Jacobi operator. Then there are constants $\mu_{\pi,\varrho}$ and $\widetilde{\mu}_{n,\varrho}$ depending only on $\pi$, $\varrho$ and $n$, $\varrho$, respectively, such that
\[
\|\mathcal{J}\|_2\leq\mu_{\pi,\varrho}, \quad 0\leq\mu_{\pi,\varrho}<\widetilde{\mu}_{n,\varrho}<1.
\]
\end{Theorem}

\begin{proof}
Let $a\in\mathbb{R}^K$ be an arbitrary nonzero vector and let $a'=\mathcal{J} a$. To track how $a'$ is obtained from $a$, we can assume $\mathcal{J}=\bm{\mathcal{R}}_{i_{M-1}j_{M-1}}\bm{\mathcal{R}}_{i_{M-2}j_{M-2}}\ldots\bm{\mathcal{R}}_{i_0j_0}$, where $\mathcal{O}= (i_0,j_0),\ldots, (i_{M-1},j_{M-1})$. If we define
\begin{equation}\label{vec-ak}
a^{(k+1)}=\bm{\mathcal{R}}_{i_kj_k}a^{(k)}, \quad 0\leq k\leq M-1; \quad a^{(0)}=a,
\end{equation}
we obtain $a'=a^{(M)}$. Recall that
Algorithm~\ref{bjanagm} describes the $k$th step of the process~\eqref{vec-ak}, i.e., how the vector $a^{(k+1)}$ is obtained from $a^{(k)}$. That algorithm computes the matrix $A^{(k+1)}=\ve_0^{-1}(a^{(k+1)})$ from the matrix $A^{(k)}=\ve_0^{-1}(a^{(k)})$. Note that
\[
S^2(A)=\|a\|_2^2 \qquad \text{and} \qquad S^2(A')=\|a'\|_2^2=\|\mathcal{J} a\|_2^2.
\]
If we prove
\begin{equation}\label{haha1}
S^2(A')\leq\eta_{\pi,\varrho}S^2(A), \quad \eta_{\pi,\varrho}<\widetilde{\eta}_{n,\varrho}<1,
\end{equation}
and take into account that $a$ is an arbitrary nonzero vector, we will straightforwardly obtain
\[
\|\mathcal{J}_{\mathcal{O}}\|_2=\max_{a\neq0}\frac{\|\mathcal{J}_{\mathcal{O}} \ a\|_2}{\|a\|_2}\leq\mu_{\pi,\varrho}<
\widetilde{\mu}_{n,\varrho}, \qquad \mu_{\pi,\varrho}=\sqrt{\eta_{\pi,\varrho}}, \quad \widetilde{\mu}_{n,\varrho} = \sqrt{\widetilde{\eta}_{n,\varrho}}.
\]
To prove \eqref{haha1}, we can rely on the proof of Theorem~\ref{tm:b1}. Indeed, let us compare the computation of the matrix
$A^{(k+1)}$ from $A^{(k)}$ using Algorithm~\ref{bjanagm} with the $k$th step of the block Jacobi method.
If we neglect the diagonal blocks, both amount to the same procedure, except for the fact that the block Jacobi method actually computes the orthogonal elementary matrix which diagonalizes the pivot submatrix, while in process~\eqref{vec-ak} that transformation is given via the matrix $\bm{\mathcal{R}}_{i_kj_k}=\bm{\mathcal{R}}_{i_kj_k}(\widehat{U}_{k})$.
The two procedures will naturally generate different iteration matrices, but all estimates and the whole proof will be the same. The quantity $\zeta_l$ from relation~\eqref{varsigma} will be different for the two procedures, but all that is needed for the proof is that $\zeta_l$ is uniformly bounded from below by some positive constant, which is certainly satisfied.
$\qquad$ \end{proof}

In the special case when $\pi=(1,1,\ldots,1)$, all blocks are $1\times1$ matrices, i.e., the single elements, so the block method reduces to the standard Jacobi method. In this case we will denote the class $\mathcal{B}_c^{(n)}$ by $\mathcal{C}_c^{(n)}$. Theorem~\ref{tm:b1} then reduces to the following corollary.

\begin{Corollary}\label{tm:c1}
Let $A\in\mathbf{S}_{n}$, $\mathcal{O}\in\mathcal{C}_c^{(n)}$ and let $A'$ be obtained from $A$ by applying one sweep of the cyclic Jacobi method defined by the strategy $I_{\mathcal{O}}$, with rotation angles from the interval $[-\frac{\pi}{4},\frac{\pi}{4}]$.
Then there is a constant $\eta_n$ depending only on $n$, such that
\[
S^2(A')\leq \eta_n S^2(A), \quad 0\leq \eta_n<1.
\]
\end{Corollary}

\begin{proof}
The proof has been moved to the Appendix.
$\qquad$ \end{proof}

In this special case, the class $\mathcal{C}_c^{(n)}$ is a subset of the set of Nazareth's orderings from~\cite{naz-75}. However, note that the bounds obtained here are much better than those in~\cite{naz-75}. To illustrate that, observe that for $n=3$ (resp.\@ $n=4$) the value of $\eta_n$ is equal to $\max\{\frac{3}{4},\frac{3}{4}\}=\frac{3}{4}$
(resp.\@ $\max\{\frac{7}{8},\frac{27}{28}\}=\frac{27}{28}$). In~\cite{naz-75} the corresponding bounds are $1-1/(3\cdot 2^{104})$ and $1-1/(6\cdot 2^{294})$. This comparison has also been studied in~\cite[page 57]{B}.

\subsection{The classes $\mathcal{B}_{cp}^{(m)}$, $\mathcal{B}_{rp}^{(m)}$ and $\mathcal{B}_{sp}^{(m)}$}\label{sec:3.2}

The same results hold for the class of cyclic pivot strategies which take the pivot blocks from the block-rows.
Let
\begin{align*}
\mathcal{B}_r^{(m)} & = \big{\{}\mathcal{O}\in\mathcal{\Ol}(\mathcal{P}_m) \ \big{|} \ \mathcal{O}=(m-1,m),(m-2,\widetilde{\tau}_{m-2}(m-1)),(m-2,\widetilde{\tau}_{m-2}(m)), \\
& \qquad \ldots,(1,\widetilde{\tau}_{1}(2)),\ldots,(1,\widetilde{\tau}_{1}(m)), \quad \widetilde{\tau}_{i}\in\Pi^{(i+1,m)}, \ 1\leq i\leq m-2 \big{\}}.
\end{align*}
The set $\mathcal{B}_r^{(m)}$ is a part of the class of \emph{row-wise orderings with permutations} of the set $\mathcal{P}_m$ (see Definition~\ref{tm:def3.5} below).
A typical ordering from $\mathcal{O}\in\mathcal{B}_r^{(6)}$ is represented by $\mathsf{M}_{\mathcal{O}}$ below. The second matrix $\mathsf{M}_{\widetilde{\mathcal{O}}}$ is defined by  $\widetilde{\mathcal{O}}\stackrel{\mathsf{w}}{\sim}\mathcal{O}$.
{\small
\[\mbox{\normalsize $\mathsf{M}_{\mathcal{O}}$} =\left[
              \begin{array}{cccccc}
                * & 11 & 13 & 12 & 10 & 14 \\
                10 & * & 9 & 7 & 6 & 8 \\
                11 & 9 & * & 5 & 3 & 4 \\
                12 & 6 & 5 & * & 1 & 2 \\
                13 & 7 & 3 & 1 & * & 0 \\
                14 & 8 & 4 & 2 & 0 & * \\
              \end{array}\right], \qquad
\mbox{\normalsize $\mathsf{M}_{\widetilde{\mathcal{O}}}$}=\left[
              \begin{array}{cccccc}
                * & 14 & 1 & 0 & 11 & 2 \\
                14& * & 13 & 10 & 7 & 12 \\
                1 & 13 & * & 9 & 6 & 8 \\
                0 & 10 & 9 & * & 4 & 5 \\
                11& 7 & 6 & 4 & * & 3 \\
                2 & 12 & 8 & 5 & 3 & * \\
              \end{array}\right].
\]}
From the matrix $\mathsf{M}_{\mathcal{O}}$ we can see that the permutation $\widetilde{\tau}_i$ from~\eqref{b1def} is linked with the $i$th block-row of the matrix.
It is immediately clear that
\begin{equation}\label{brcm}
\mathcal{B}_r^{(m)} = \big{\{}\mathcal{O}(\widetilde{\mathsf{e}})\mid \mathcal{O}\in\mathcal{B}_c^{(m)}\big{\}},
\end{equation}
where $\widetilde{\mathsf{e}}$ is defined by relation (\ref{etilde}). Theorems~\ref{tm:b1} and~\ref{tm:b1jop} remain to hold for $\mathcal{O}\in\mathcal{B}_r^{(m)}$. The proofs are almost identical to those for the case $\mathcal{O}\in\mathcal{B}_c^{(m)}$.

The version of theorem~\ref{tm:b1jop} for $\mathcal{O}\in\mathcal{B}_r^{(m)}$ follows directly from the original Theorem~\ref{tm:b1jop} combined with Theorem~\ref{tm:prop_permstr}$(i)$ (when the permutation $\mathsf{p}$ is specified to be  $\widetilde{\mathsf{e}}$).

Another way how Theorem~\ref{tm:b1} can be established for $\mathcal{O}\in\mathcal{B}_r^{(m)}$ is to follow the lines of the proof of Theorem~\ref{tm:bsjop} (below), but then a version of Theorem~\ref{tm:b1jop} for $\mathcal{O}\in\mathcal{B}_r^{(m)}$ should be proved first.

If $\pi=(1,1,\ldots,1)$, then the class of orderings $\mathcal{B}_r^{(n)}$ is denoted by $\mathcal{C}_r^{(n)}$.
Corollary~\ref{tm:c1} holds with the same constant $\eta_n$ provided that $\mathcal{C}_c^{(n)}$ is replaced by $\mathcal{C}_r^{(n)}$.

We are interested in two more classes of pivot strategies for the block methods. The first (resp.\@ second) one selects the pivot blocks by block-columns (resp.\@ block-rows), but now from the last one to the second one (resp.\@ from the first one to the next-to-last one). They are defined as
\[
\overleftarrow{\mathcal{B}}_c^{(m)} = \big{\{} \mathcal{O}\in\mathcal{\Ol}(\mathcal{P}_m) \ \big{|} \ \mathcal{O}^{\leftarrow }\in \mathcal{B}_c^{(m)} \big{\}}, \qquad
\overleftarrow{\mathcal{B}}_r^{(m)} = \big{\{} \mathcal{O}\in\mathcal{\Ol}(\mathcal{P}_m) \ \big{|} \ \mathcal{O}^{\leftarrow }\in \mathcal{B}_r^{(m)} \big{\}}.
\]
Typical orderings from $\overleftarrow{\mathcal{B}}_c^{(6)}$ and $\overleftarrow{\mathcal{B}}_r^{(6)}$ are represented by $\mathsf{M}_{\overleftarrow{\mathcal{O}}_c}$ and $\mathsf{M}_{\overleftarrow{\mathcal{O}}_r}$ below.
{\small
\[\mathsf{M}_{\overleftarrow{\mathcal{O}}_c}=
              \left[
              \begin{array}{cccccc}
                * & 14 & 12 & 10 & 5 & 2 \\
                14 & * & 13 & 9 & 6 & 4 \\
                12 & 13 & * & 11 & 7 & 1 \\
                10 & 9 & 11 & * & 8 & 3 \\
                5 & 6 & 7 & 8 & * & 0 \\
                2 & 4 & 1 & 3 & 0 & *
              \end{array}
            \right], \qquad
\mathsf{M}_{\overleftarrow{\mathcal{O}}_r}=
              \left[
              \begin{array}{cccccc}
                * & 4 & 3 & 2 & 1 & 0 \\
                4 & * & 5 & 8 & 7 & 6 \\
                3 & 5 & * & 9 & 11 & 10 \\
                2 & 8 & 9 & * & 13 & 12 \\
                1 & 7 & 11 & 13 & * & 14 \\
                0 & 6 & 10 & 12 & 14 & *
              \end{array}
            \right].
\]}
As it has already been noticed at the end of Section~\ref{sec:2.2}, Lemma~\ref{tm: lema_2.3} remains to hold if $\stackrel{w}{\sim}$ is replaced with $\stackrel{\mathsf{p}}{\sim}$. Therefore, relation (\ref{brcm}) implies
$\overleftarrow{\mathcal{B}}_r^{(m)} = \big{\{} \mathcal{O}(\widetilde{e})\mid \mathcal{O}\in \overleftarrow{\mathcal{B}}_c^{(m)} \big{\}}$.

\begin{Definition}\label{tm:def3.5}
Let $\mathcal{B}_{cp}^{(m)}=\mathcal{B}_c^{(m)} \cup \overleftarrow{\mathcal{B}}_c^{(m)}$,
$\mathcal{B}_{rp}^{(m)}=\mathcal{B}_r^{(m)} \cup \overleftarrow{\mathcal{B}}_r^{(m)}$. The set $\mathcal{B}_{cp}^{(m)}$ (resp.\@ $\mathcal{B}_{rp}^{(m)}$) is the class of the
column-wise (resp.\@ row-wise) orderings with permutations of $\mathcal{P}_m$ and $\{I_{\mathcal{O}}\mid \mathcal{O}\in\mathcal{B}_{cp}^{(m)}\}$
(resp.\@ $\{I_{\mathcal{O}}\mid \mathcal{O}\in\mathcal{B}_{rp}^{(m)}\}$) is the class of  the column-cyclic (resp.\@ row-cyclic) strategies with permutations.

The set $\displaystyle \mathcal{B}_{sp}^{(m)}= \mathcal{B}_{cp}^{(m)}\cup \mathcal{B}_{rp}^{(m)}$ is the class of the serial orderings with permutations of $\mathcal{P}_m$ and  $\{I_{\mathcal{O}}\mid \mathcal{O}\in\mathcal{B}_{sp}^{(m)}\}$ is the class of the serial strategies with permutations.
\end{Definition}

\begin{Theorem}\label{tm:bsjop}
Let $\pi=(n_1,\ldots ,n_m)$ be a partition of $n$, $\mathcal{O}\in\mathcal{B}_{sp}^{(m)}$ and let $\mathcal{J}\in\Jl_{\!\!\! \mathcal{O}}^{\ubce(\varrho)}$ be a block Jacobi operator. Then there are constants $\mu_{\pi,\varrho}$ and $\widetilde{\mu}_{n,\varrho}$ depending only on $\pi$, $\varrho$ and $n$, $\varrho$, respectively, such that
\[
\|\mathcal{J}\|_2\leq\mu_{\pi,\varrho}, \quad 0\leq\mu_{\pi,\varrho}<\widetilde{\mu}_{n,\varrho}<1.
\]
\end{Theorem}

\begin{proof}
If $\mathcal{O}\in\mathcal{B}_{c}^{(m)}$ (resp.\@ $\mathcal{O}\in\mathcal{B}_{r}^{(m)}$), then the theorem reduces to Theorem~\ref{tm:b1jop} (resp.\@ Theorem~\ref{tm:b1jop} combined with Theorem~\ref{tm:prop_permstr}$(i)$). If $\mathcal{O}\in\overleftarrow{\mathcal{B}}_c^{(m)}$ (resp.\@ $\mathcal{O}\in\overleftarrow{\mathcal{B}}_r^{(m)}$), then additionally Proposition~\ref{tm:prop_invstr} is used.
$\qquad$ \end{proof}

The following result is a corollary of Theorem~\ref{tm:bsjop}, but because of its importance, it is stated as a stand-alone result.

\begin{Theorem}\label{tm:bs}
Let $\pi=(n_1,\ldots,n_m)$ be a partition of $n$, $\mathcal{O}\in\mathcal{B}_{sp}^{(m)}$, and let $A\in\mathbf{S}_{n}$ be partitioned as in relation~\eqref{blokmatrica}. Let $A'$ be obtained from $A$ by applying one sweep of the cyclic block Jacobi method defined by the strategy $I_{\mathcal{O}}$. If all transformation matrices are from the class $\ubce(\varrho)$, then there are constants $\eta_{\pi,\varrho}$ (depending only on $\pi$, $\varrho$) and $\widetilde{\eta}_{n,\varrho}$ (depending only on $n$, $\varrho$) such that
\[
S^2(A')\leq \eta_{\pi,\varrho} S^2(A), \quad 0\leq \eta_{\pi,\varrho}<\widetilde{\eta}_{n,\varrho}<1.
\]
\end{Theorem}

\begin{proof}
Let $\mathcal{O}=(i_0,j_0),(i_1,j_1),\ldots,(i_{M-1},j_{M-1})$.
The method in the statement of the theorem generates the recurrence relation of the form (\ref{jacobiblokagm}). If we observe how the elements in the block upper-triangle are being updated, we arrive at the recursion
\[
a^{(k+1)} = \bm{\mathcal{R}}_{ij}a^{(k)}, \quad k\geq 0; \quad  a^{(0)}=a=\textrm{vec}(A).
\]
Here, for each $k$, $a^{(k)}= \textrm{vec}(A^{(k)})\in\mathbb{R}^{K}$ and $\bm{\mathcal{R}}_{ij}=\bm{\mathcal{R}}_{ij}(\widehat{U}_k)$ is the block Jacobi annihilator associated with step $k$ of the method. We have $(i,j)=(i(k),j(k))=(i_k,j_k)$ for $0\leq k\leq M-1$. After the first sweep is completed, one obtains
\begin{equation}\label{in-thmbs-1}
a^{(M)} = \mathcal{J}_{\mathcal{O}}a, \quad \mathcal{J}_{\mathcal{O}}=\bm{\mathcal{R}}_{i(M-1)j(M-1)}\cdots \bm{\mathcal{R}}_{i(0)j(0)}.
\end{equation}
Since all transformation matrices in the block Jacobi method are from the class $\ubce(\varrho)$, we have $\bm{\mathcal{R}}_{ij}(\widehat{U}_k)$ $\in$ $\Ilb_{ij}^{\ubce(\varrho)}$. Therefore, $\mathcal{J}_{\mathcal{O}}\in\Jl_{\!\!\! \mathcal{O}}^{\ubce(\varrho )}$
and by applying Theorem~\ref{tm:bsjop} one obtains $\|\mathcal{J}\|_2\leq\mu_{\pi,\varrho}$,
$\mu_{\pi,\varrho}<\widetilde{\mu}_{n,\varrho}<1$. Hence, if one takes the Euclidean vector norm of both sides of the left equation in (\ref{in-thmbs-1}), it follows that
\[
S^2(A^{(M)})= \|a^{(M)}\|_2^2 \leq \mu_{\pi,\varrho}^2 \|a\|_2^2 = \mu_{\pi,\varrho}^2 S^2(A).
\]
It remains to set $\eta_{\pi,\varrho}=\mu_{\pi,\varrho}^2$ and $\widetilde{\eta}_{n,\varrho}=\widetilde{\mu}_{n,\varrho}^2$.
$\qquad$ \end{proof}

Obviously, since the assertion of Theorem~\ref{tm:bs} holds for any single sweep, we conclude that $S(A^{(tM)})\rightarrow 0$  as $t\rightarrow\infty$. Since the sequence $S(A^{(k)})$, $k\geq 0$, is nonincreasing, one obtains $S(A^{(k)})\rightarrow 0$ as $k\rightarrow\infty$. Together with Theorem~\ref{tm:gl_conv-opp} this implies the global convergence.

In the case $\pi =(1,1,\ldots\ ,1)$ we have $m=n$. We can write
$\mathcal{C}_{cp}^{(n)}$, $\mathcal{C}_{rp}^{(n)}$ and $\mathcal{C}_{sp}^{(n)}$ in the places of $\mathcal{B}_{cp}^{(m)}$, $\mathcal{B}_{rp}^{(m)}$ and $\mathcal{B}_{sp}^{(m)}$, respectively.
Once again, Corollary~\ref{tm:c1} holds with the same constant $\eta_n$ if  $\mathcal{C}_c^{(n)}$ is replaced by $\mathcal{C}_{sp}^{(n)}$. Theorem~\ref{tm:bsjop} also holds with
$\mu_{n,\varrho}$ replacing $\mu_{\pi,\varrho}$.

\subsection{Generalized serial strategies}

To further enlarge the class of convergent strategies one can start with elements from $\mathcal{B}_{sp}^{(m)}$ and use all conceivable chains which comprise the equivalence relations  $\sim$, $\stackrel{\mathsf{s}}{\sim}$, $\stackrel{\mathsf{w}}{\sim}$ and $\stackrel{\mathsf{p}}{\sim}$. Fortunately, using Proposition~\ref{tm:peq1} and Definition~\ref{tm:can_chain}, we know that the set obtained that way can actually be obtained by using just one $\stackrel{\mathsf{w}}{\sim}$ and one $\stackrel{\mathsf{p}}{\sim}$.

\begin{Definition}\label{tm:gss}
Let
\begin{align*}
\mathcal{B}_{spg}^{(m)} & = \big{\{}\mathcal{O}\in\Ol(\mathcal{P}_m) \ \big{|} \ \mathcal{O}\stackrel{\mathsf{p}}{\sim}\mathcal{O'}\sim\mathcal{O''} \ \text{or} \
\mathcal{O}\sim\mathcal{O'}\stackrel{\mathsf{p}}{\sim}\mathcal{O''}, \ \mathcal{O''}\in\mathcal{B}_{sp}^{(m)}\big{\}}, \\
\mathcal{B}_{sg}^{(m)} & = \big{\{}\mathcal{O}\in\Ol(\mathcal{P}_m) \ \big{|} \ \mathcal{O}\stackrel{\mathsf{p}}{\sim}\mathcal{O'}\stackrel{\mathsf{w}}{\sim}\mathcal{O''} \ \text{or} \ \mathcal{O}\stackrel{\mathsf{w}}{\sim}\mathcal{O'}\stackrel{\mathsf{p}}{\sim}\mathcal{O''}, \ \mathcal{O''}\in\mathcal{B}_{sp}^{(m)}\big{\}},
\end{align*}
where the chains are in the canonical form and $\mathcal{O'}\in\Ol(\mathcal{P}_m)$.
The set $\mathcal{B}_{sg}^{(m)}$ (resp.\@ $\{I_{\mathcal{O}} \mid \mathcal{O}\in \mathcal{B}_{sg}^{(m)}\}$) is the class of \emph{generalized serial pivot orderings} of $\mathcal{P}_m$ (resp.\@ \emph{generalized serial pivot strategies}). The set $\mathcal{B}_{spg}^{(m)}$ is a subclass of $\mathcal{B}_{sg}^{(m)}$ whose elements are linked by chains that do not use shifts.
\end{Definition}

\begin{Theorem}\label{tm:corr3-8}
Let $\pi=(n_1,\ldots,n_m)$ be a partition of $n$, $\mathcal{O}\in\mathcal{B}_{spg}^{(m)}$ and let $\mathcal{J}\in\Jl_{\!\!\! \mathcal{O}}^{\ubce_{\pi}(\varrho )}$ be the block Jacobi operator. Suppose that either
$\mathcal{O}\stackrel{\mathsf{p}}{\sim}\mathcal{O'}\sim\mathcal{O''}\in\mathcal{B}_{sp}^{(m)}$ with $\mathcal{O'}=\mathcal{O}(\mathsf{q})$, or $\mathcal{O}\sim\mathcal{O'}\stackrel{\mathsf{p}}{\sim}\mathcal{O''}\in\mathcal{B}_{sp}^{(m)}$ with $\mathcal{O''}=\mathcal{O'}(\mathsf{q})$, for some permutation $\mathsf{q}$ of the set $\mathcal{S}_m$. Then there exist constants $\mu_{\pi_{\mathsf{q}},\varrho}$ and $\widetilde{\mu}_{n,\varrho}$, depending only on $\pi_{\mathsf{q}}$, $\varrho$ and $n$, $\varrho$, respectively, such that
\[
\|\mathcal{J}\|_2\leq\mu_{\pi_{\mathsf{q}}}, \quad 0\leq\mu_{\pi_{\mathsf{q}},\varrho}<\widetilde{\mu}_{n,\varrho}<1.
\]
\end{Theorem}

\begin{proof}
Let us first consider the case $\mathcal{O}\stackrel{\mathsf{p}}{\sim}\mathcal{O'}\sim\mathcal{O''}\in\mathcal{B}_{sp}^{(m)}$ with $\mathcal{O'}=\mathcal{O}(\mathsf{q})$. Let
\[
\widetilde{\mathcal{B}}_{sp}^{(m)} = \big{\{}\mathcal{O'}\in\Ol(\mathcal{P}_m) \ \big{|} \ \mathcal{O'}\sim\mathcal{O''}, \mathcal{O''}\in\mathcal{B}_{sp}^{(m)}\big{\}}.
 \]
Theorem~\ref{tm:bsjop} and Lemma~\ref{tm:lemma_spran} imply
\begin{equation}\label{rel3.7}
\|\mathcal{J}_{\mathcal{O'}}\|_2\leq\mu_{\pi_{\mathsf{q}},\varrho}, \quad
0\leq \mu_{\pi_{\mathsf{q}},\varrho}<\widetilde{\mu}_{n,\varrho}<1, \ \mathcal{J}_{\mathcal{O'}}\in\Jl_{\!\!\! \mathcal{O'}}^{\ubce_{\pi_{\mathsf{q}}}(\varrho )}\!\!, \
\mathcal{O'}\in \widetilde{\mathcal{B}}_{sp}^{(m)}.
\end{equation}
Since $\mathcal{O}=\mathcal{O'}(\mathsf{q^{-1}})$ and $\pi_{\mathsf{q}^{-1}\circ\mathsf{q}}=\pi$, formula (\ref{rel3.7}) and Theorem~\ref{tm:prop_permstr}$(i)$ imply
\begin{equation}\label{rel3.7a}
\|\mathcal{J}_{\mathcal{O}}\|_2\leq\mu_{\pi_{\mathsf{q}},\varrho}, \quad 0\leq\mu_{\pi_{\mathsf{q}},\varrho}<\widetilde{\mu}_{n,\varrho}<1,
\quad \text{for any} \ \mathcal{J}_{\mathcal{O}}\in \Jl_{\!\!\! \mathcal{O}}^{\ubce_{\pi}(\varrho )}.
\end{equation}
Hence, relation (\ref{rel3.7a}) holds for $\mathcal{J}$, which in turn proves the theorem.

Now suppose that $\mathcal{O}\sim\mathcal{O'}\stackrel{\mathsf{p}}{\sim}\mathcal{O''}\in\mathcal{B}_{sp}^{(m)}$ with $\mathcal{O''}=\mathcal{O'}(\mathsf{q})$. Let
\[
\widehat{\mathcal{B}}_{sp}^{(m)} = \big{\{}\mathcal{O'}\in\Ol(\mathcal{P}_m) \ | \ \mathcal{O'}=\mathcal{O''}(\mathsf{q}^{-1}), \ \mathcal{O''}\in\mathcal{B}_{sp}^{(m)}\big{\}}.
\]
By Theorem~\ref{tm:bsjop} and Theorem~\ref{tm:prop_permstr}$(i)$ we have
\[
\|\mathcal{J}_{\mathcal{O'}}\|_2\leq\mu_{\pi_{\mathsf{q}},\varrho}, \quad 0\leq \mu_{\pi_{\mathsf{q}},\varrho}<\widetilde{\mu}_{n,\varrho}<1, \ \mathcal{J}_{\mathcal{O'}}\in\Jl_{\!\!\! \mathcal{O}}^{\ubce_{\pi (\varrho )}}, \ \mathcal{O}'\in \widehat{\mathcal{B}}_{sp}^{(m)}.
\]
Here we have used $\pi_{\mathsf{q}^{-1}\circ\mathsf{q}}=\pi$. The last formula holds for any
$\mathcal{J}_{\mathcal{O'}}\in\Jl_{\!\!\! \mathcal{O}}^{\ubce_{\pi (\varrho )}}$.
Since $\mathcal{O}\sim\mathcal{O'}$, Lemma~\ref{tm:lemma_spran} implies that the same formula holds for any
$\mathcal{J}_{\mathcal{O}}\in\Jl_{\!\!\! \mathcal{O}}^{\ubce_{\pi (\varrho )}}$.
This completes the proof since $\mathcal{J}$ is just one of those $\mathcal{J}_{\mathcal{O}}$.
$\qquad$ \end{proof}

\begin{Theorem}\label{tm:tm3-9}
Let $\pi=(n_1,\ldots ,n_m)$ be a partition of $n$ and $\mathcal{O}\in\mathcal{B}_{sg}^{(m)}$. Suppose that the chain connecting $\mathcal{O}$ to $\mathcal{O''}\in\mathcal{B}_{sp}^{(m)}$ contains $d$ shift equivalences. Moreover, suppose that either $\mathcal{O}\stackrel{\mathsf{p}}{\sim}\mathcal{O'}\stackrel{\mathsf{w}}{\sim}\mathcal{O''}\in
\mathcal{B}_{sp}^{(m)}$ with $\mathcal{O'}=\mathcal{O}(\mathsf{q})$ holds, or $\mathcal{O}\stackrel{\mathsf{w}}{\sim}\mathcal{O'}\stackrel{\mathsf{p}}{\sim}\mathcal{O''}\in
\mathcal{B}_{sp}^{(m)}$ with $\mathcal{O''}=\mathcal{O'}(\mathsf{q})$ holds, for some permutation $\mathsf{q}$ of the set $\mathcal{S}_m$.
Then there exist constants $\mu_{\pi,\varrho}$ and $\widetilde{\mu}_{n,\varrho}$ depending only on $\pi$, $\varrho$ and $n$, $\varrho$, respectively, such that for any $d+1$ block Jacobi operators $\mathcal{J}_1, \mathcal{J}_2,\ldots ,\mathcal{J}_{d+1}\in\Jl_{\!\!\! \mathcal{O}}^{\ubce_{\pi}(\varrho)}$ one has
\[
\|\mathcal{J}_1 \mathcal{J}_2\cdots \mathcal{J}_{d+1}\|_2 \leq \mu_{\pi_{\mathsf{q}},\varrho}, \quad 0\leq\mu_{\pi_{\mathsf{q}},\varrho}<\widetilde{\mu}_{n,\varrho}<1.
\]
\end{Theorem}

\begin{proof}
Let us first consider the case $\mathcal{O}\stackrel{\mathsf{p}}{\sim}\mathcal{O'}\stackrel{\mathsf{w}}{\sim}\mathcal{O''}\in
\mathcal{B}_{sp}^{(m)}$ with $\mathcal{O'}=\mathcal{O}(\mathsf{q})$. Denote
\[
\widetilde{\mathcal{B}}_{sp}^{(m)} = \big{\{}\mathcal{O'}\in\Ol(\mathcal{P}_m) \ \big{|} \ \mathcal{O'}\stackrel{\mathsf{w}}{\sim}\mathcal{O''}, \ \mathcal{O''}\in\mathcal{B}_{sp}^{(m)}\big{\}}.
\]
Theorem~\ref{tm:bsjop} and Proposition~\ref{tm:prop_norman} imply
\begin{equation}\label{rel3.7b}
\|\mathcal{J}_1'\mathcal{J}_2'\cdots \mathcal{J}_{d+1}'\|_2\leq\mu_{\pi_{\mathsf{q}},\varrho}, \quad
\mathcal{J}_1',\mathcal{J}_2',\ldots,\mathcal{J}_{d+1}'\in\Jl_{\!\!\! \mathcal{O'}}^{\ubce_{\pi_{\mathsf{q}}}(\varrho )}\!\!, \ \mathcal{O'}\in \widetilde{\mathcal{B}}_{sp}^{(m)},
\end{equation}
where $0\leq \mu_{\pi_{\mathsf{q}},\varrho}<\widetilde{\mu}_{n,\varrho}<1$.
Since $\mathcal{O}=\mathcal{O'}(\mathsf{q^{-1}})$ and $\pi_{\mathsf{q}^{-1}\circ\,\mathsf{q}}=\pi$, relation (\ref{rel3.7b}) and Theorem~\ref{tm:prop_permstr}$(ii)$ give
\[
\|\mathcal{J}_1\mathcal{J}_2\cdots \mathcal{J}_{d+1}\|_2\leq\mu_{\pi_{\mathsf{q}},\varrho} \quad \text{for any} \
\mathcal{J}_1,\mathcal{J}_2,\ldots,\mathcal{J}_{d+1}\in\Jl_{\!\!\! \mathcal{O}}^{\ubce_{\pi}(\varrho )}\!\!,
\]
which proves the theorem.

Now suppose that $\mathcal{O}\stackrel{\mathsf{w}}{\sim}\mathcal{O'}\stackrel{\mathsf{p}}{\sim}\mathcal{O''}\in\mathcal{B}_{sp}^{(m)}$ with $\mathcal{O''}=\mathcal{O'}(\mathsf{q})$. Let
\[
\widehat{\mathcal{B}}_{sp}^{(m)} = \big{\{}\mathcal{O'}\in\Ol(\mathcal{P}_m) \ \big{|} \ \mathcal{O'}=\mathcal{O''}(\mathsf{q}^{-1}),\ \mathcal{O''}\in\mathcal{B}_{sp}^{(m)}\big{\}}.
 \]
By Theorem~\ref{tm:bsjop} and Theorem~\ref{tm:prop_permstr}$(i)$ we have
\[
\|\mathcal{J}_{\mathcal{O'}}\|_2\leq\mu_{\pi_{\mathsf{q}},\varrho}, \quad 0\leq \mu_{\pi_{\mathsf{q}},\varrho}<\widetilde{\mu}_{n,\varrho}<1, \ \mathcal{J}_{\mathcal{O'}}\in\Jl_{\!\!\! \mathcal{O}}^{\ubce_{\pi (\varrho )}}, \ \mathcal{O}'\in \widehat{\mathcal{B}}_{sp}^{(m)}.
\]
Here we used $\pi_{\mathsf{q}^{-1}\circ\,\mathsf{q}}=\pi$ once again. Since $\mathcal{O}\stackrel{\mathsf{w}}{\sim}\mathcal{O'}$, Proposition~\ref{tm:prop_norman} completes the proof.
$\qquad$ \end{proof}

We end this section by shifting our attention from block Jacobi operators to cyclic block Jacobi methods, defined by the generalized serial strategies.

\begin{Theorem}\label{tm:tm-spg}
Let $\pi=(n_1,\ldots ,n_m)$ be a partition of $n$, $\mathcal{O}\in\mathcal{B}_{spg}^{(m)}$ and let $A\in\mathbf{S}_{n}$ be partitioned as in relation \eqref{blokmatrica}. Let $A'$ be obtained from $A$ by applying one sweep of the cyclic block Jacobi method defined by the strategy $I_{\mathcal{O}}$.  If all transformation matrices are from the class $\ubce_{\pi_{\mathsf{q}}}(\varrho )$ for an appropriate permutation $\mathsf{q}$ of the set $\mathcal{S}_m$, then there are constants $\eta_{\pi_{\mathsf{q}},\varrho}$ and $\widetilde{\eta}_{n,\varrho}$ depending only on $\pi_{\mathsf{q}}$, $\varrho$ and $n$,  $\varrho$, respectively, such that
\[
S^2(A')\leq \eta_{\pi_{\mathsf{q}},\varrho} S^2(A), \quad 0\leq \eta_{\pi_{\mathsf{q}},\varrho}<\widetilde{\eta}_{n,\varrho}<1.
\]
\end{Theorem}

\begin{proof}
The proof is almost identical to the proof of Theorem~\ref{tm:bs}.
The difference is that $\mathcal{O}\in\mathcal{B}_{spg}^{(m)}$ and we use Theorem~\ref{tm:corr3-8} instead of Theorem~\ref{tm:bsjop}.
$\qquad$ \end{proof}

\begin{Theorem}\label{tm:tm-spg-d}
Let $\pi=(n_1,\ldots,n_m)$ be a partition of $n$, $\mathcal{O}\in\mathcal{B}_{sg}^{(m)}$ and let $A\in\mathbf{S}_{n}$ be partitioned as in relation \eqref{blokmatrica}. Suppose that the chain connecting $\mathcal{O}$ and $\mathcal{O''}\in\mathcal{B}_{sp}^{(m)}$ is in the canonical form and contains $d$ shift equivalences.
Let $A'$ be obtained from $A$ by applying $d+1$ sweeps of the cyclic block Jacobi method defined by the strategy $I_{\mathcal{O}}$. If all transformation matrices are from the class $\ubce_{\pi_{\mathsf{q}}}(\varrho )$ for an appropriate permutation $\mathsf{q}$ of the set $\mathcal{P}_m$, then there are constants $\eta_{\pi_{\mathsf{q}},\varrho}$ and $\widetilde{\eta}_{n,\varrho}$ depending only on $\pi_{\mathsf{q}}$, $\varrho$ and $n$, $\varrho$, respectively, such that
\[
S^2(A')\leq \eta_{\pi_{\mathsf{q}},\varrho} S^2(A), \quad 0\leq \eta_{\pi_{\mathsf{q}},\varrho}<\widetilde{\eta}_{n,\varrho}<1.
\]
Here $\mathsf{q}$ is an appropriate permutation of the set $\mathcal{P}_m$.
\end{Theorem}

\begin{proof}
The proof follows the lines of the proof of Theorem~\ref{tm:bs}. Since we consider $d+1$ sweeps, instead of relation (\ref{in-thmbs-1}), we will obtain
$$a'=a^{((d+1)M)} = \mathcal{J}_{\mathcal{O}}^{[d+1]}\mathcal{J}_{\mathcal{O}}^{[d]}\cdots
\mathcal{J}_{\mathcal{O}}^{[1]}a,\quad \mathcal{O}\in\mathcal{B}_{sg}^{(m)}.
$$
Here $\mathcal{J}_{\mathcal{O}}^{[s]\,}$ is the block Jacobi operator associated with cycle $s$ of the block Jacobi method and $a=\ve (A)$, $a'=\ve (A')$.
From Theorem~\ref{tm:tm3-9} we know that
\[
\|\mathcal{J}_{\mathcal{O}}^{[d+1]}\mathcal{J}_{\mathcal{O}}^{[d]}\cdots
\mathcal{J}_{\mathcal{O}}^{[1]}\|_2 \leq \mu_{\pi_{\mathsf{q}},\varrho}, \quad 0\leq\mu_{\pi_{\mathsf{q}},\varrho}<\widetilde{\mu}_{n,\varrho}<1,
\]
so
\[
S^2(A')= \|a'\|_2^2 \leq \mu_{\pi_{\mathsf{q}},\varrho}^2 \|a\|_2^2 \leq \mu_{\pi,\varrho}^2 S^2(A).
\]
It remains to set $\eta_{\pi_{\mathsf{q}},\varrho}=\mu_{\pi_{\mathsf{q}},\varrho}^2$ and $\widetilde{\eta}_{n\varrho}=\widetilde{\mu}_{n,\varrho}^2$.
$\qquad$ \end{proof}

In the case $\pi =(1,1,\ldots,1)$, we have $m=n$ and we can use notation $\mathcal{C}_{spg}^{(m)}$ and $\mathcal{C}_{sg}^{(m)}$ instead of $\mathcal{B}_{spg}^{(m)}$ and $\mathcal{B}_{sg}^{(m)}$, respectively. Corollary~\ref{tm:c1} holds with the same constant $\eta_n$ if $\mathcal{C}_c^{(n)}$ is replaced with $\mathcal{C}_{spg}^{(n)}$. Theorems \ref{tm:tm3-9} -- \ref{tm:tm-spg-d} also hold with $\mu_{n,\varrho}$ instead of $\mu_{\pi_{\mathsf{q}},\varrho}$.

It is not easy to count how many pivot orderings are contained in $\mathcal{B}_{sg}^{(m)}$. In $\mathcal{B}_{c}^{(m)}$ we have $2!\cdot 3!\cdots (m-1)!$ elements.
The sets $\overleftarrow{\mathcal{B}}_c^{(m)}$, $\mathcal{B}_{r}^{(m)}$ and $\overleftarrow{\mathcal{B}}_r^{(m)}$ have the same number of elements. Hence, for large $m$ we expect that $\mathcal{B}_{sp}^{(m)}$ contains $4\cdot 2!\cdot 3!\cdots (m-1)!$ elements. For each ordering $\mathcal{O}\in  \mathcal{B}_{sp}^{(m)}$ there are $m!$ orderings of the form $\mathcal{O}(\mathsf{p})$, so for large $m$ we expect at least $4\cdot 2!\cdot 3!\cdots m!$ elements in $\mathcal{B}_{spg}^{(m)}$ (where we have not taken the equivalences $\sim$ and $\stackrel{\mathsf{s}}{\sim}$ into account). Obviously, for small $m$ (like $m=3,4,5$) this count is not realistic.

Nevertheless, the results obtained here have been used in~\cite{beg+har15a} to prove that every cyclic Jacobi method for symmetric matrices of order $4$ is globally convergent. Note that there are altogether $720$ cyclic strategies when $n=4$.

\section{Quasi-cyclic pivot strategies}

Our next step is to enlarge the scope of generalized serial strategies by allowing repetition of some Jacobi steps within one sweep. This leads us to special quasi-cyclic pivot sequences, which are closely related to the orderings from Section~\ref{cyclic}. This change often leads to faster convergence of the Jacobi method~\cite{drm+ves-04a,drm+ves-04b,har-07}. To keep our consideration within reasonable framework, we can assume that the length of each quasi-cyclic pivot sequence is smaller than $2M$, where $M=\frac{m(m-1)}{2}$.

Our basic class of quasi-cyclic pivot sequences is derived from the class $\mathcal{B}_c^{(m)}$. More precisely,
\begin{align*}
\bar{\mathcal{B}}_c^{(m)} & = \big{\{}\mathcal{O}\in\mathcal{\Ol}(\mathcal{P}_m) \ \big{|} \ \mathcal{O}=(1,2),(\pi_{3}(1),3),(\pi_{3}(2),3),\Os_3,\ldots,(\pi_{m}(1),m), \\
& \quad \ldots,(\pi_{m}(m-1),m),\Os_m, \ \pi_{j}\in\Pi^{(1,j-1)}, \ \Os_j\in\Ol(\mathcal{S}_j), \ \mathcal{S}_j\subseteq \mathcal{P}_j, \ 3\leq j\leq m\big{\}}.
\end{align*}
Thus, the quasi-cyclic pivot strategy $I_{\mathcal{O}}$ defined by some $\mathcal{O}\in\bar{\mathcal{B}}_c^{(m)}$ selects pivot blocks by block-columns. After the pivot blocks within the $j$th block-column have all been annihilated once, it is allowed to annihilate again any block that lies within the first $j$ block-columns. Here, $3\leq j\leq m$.

\begin{Theorem}\label{kvazib1}
Let $\pi=(n_1,\ldots,n_m)$ be a partition of $n$, $\mathcal{O}\in\bar{\mathcal{B}}_c^{(m)}$ and let $A\in\mathbf{S}_{n}$ be partitioned as in relation \eqref{blokmatrica}. Let $A'$ be obtained from $A$ by applying one sweep of the quasi-cyclic block Jacobi method defined by the strategy $I_{\mathcal{O}}$. If all transformation matrices are from the class $\ubce(\varrho)$, then there are constants $\eta_{\pi,\varrho}$ and $\widetilde{\eta}_{n,\varrho}$, depending only on $\pi$, $\varrho$ and $n$, $\varrho$, respectively, such that
\[
S^2(A')\leq \eta_{\pi,\varrho} S^2(A), \quad 0\leq \eta_{\pi,\varrho}<\widetilde{\eta}_{n,\varrho}<1.
\]
\end{Theorem}

\begin{proof}
The proof has been moved to the Appendix.
$\qquad$ \end{proof}

Although the quasi-cyclic strategy from~\cite{drm+ves-04a,drm+ves-04b,har-07} belongs to the class of block-oriented strategies, its ``full block'' analogue is $I_{\mathcal{O}}$ for some special $\mathcal{O}\in\bar{\mathcal{B}}_c^{(m)}$. As has been shown in~\cite{har+sin-10,har+sin-11}, on large matrices, full-block Jacobi-type methods are generally more efficient than the block-oriented ones. This implies that the Jacobi method from LAPACK can be upgraded to full block version and Theorem~\ref{kvazib1} ensures its convergence.

Now, it is easy to prove Theorem~\ref{tm:b1jop} for $\mathcal{O}\in\bar{\mathcal{B}}_c^{(m)}$. The proof remains the same except that the word ``cyclic'' should be replaced with ``quasi-cyclic'' and the notation $\mathcal{B}_c^{(m)}$ should be replaced with $\bar{\mathcal{B}}_c^{(m)}$. The case $m=n$ is treated in the same way.

Following the ideas from Section~\ref{cyclic}, we define
$\displaystyle \bar{\mathcal{B}}_r^{(m)} = \big{\{}\mathcal{O}(\widetilde{e})\mid \mathcal{O}\in\bar{\mathcal{B}}_c^{(m)}\big{\}}$,
\[
\overleftarrow{\bar{\mathcal{B}}}_c^{(m)} = \big{\{} \mathcal{O}\in\mathcal{\Ol}(\mathcal{P}_m) \mid \mathcal{O}^{\leftarrow }\in \bar{\mathcal{B}}_c^{(m)} \big{\}}, \qquad
\overleftarrow{\bar{\mathcal{B}}}_r^{(m)} = \big{\{} \mathcal{O}\in\mathcal{\Ol}(\mathcal{P}_m) \mid \mathcal{O}^{\leftarrow }\in \bar{\mathcal{B}}_r^{(m)} \big{\}}
\]
and
\[
\bar{\mathcal{B}}_{cp}^{(m)}= \bar{\mathcal{B}}_c^{(m)} \cup \overleftarrow{\bar{\mathcal{B}}}_c^{(m)}, \qquad
\bar{\mathcal{B}}_{rp}^{(m)}= \bar{\mathcal{B}}_r^{(m)} \cup \overleftarrow{\bar{\mathcal{B}}}_r^{(m)}, \qquad
\bar{\mathcal{B}}_{sp}^{(m)}= \bar{\mathcal{B}}_{cp}^{(m)}\cup \bar{\mathcal{B}}_{rp}^{(m)}.
\]
It is easy to check that both Theorem~\ref{tm:bsjop} and Theorem~\ref{tm:bs} hold with
$\bar{\mathcal{B}}_{sp}^{(m)}$ in the place of $\mathcal{B}_{sp}^{(m)}$. Finally, we can define
\begin{align*}
\bar{\mathcal{B}}_{spg}^{(m)} & = \big{\{}\mathcal{O}\in\Ol(\mathcal{P}_m) \ \big{|} \ \mathcal{O}\stackrel{\mathsf{p}}{\sim}\mathcal{O''}\sim\mathcal{O''} \ \text{or} \
\mathcal{O}\sim\mathcal{O'}\stackrel{\mathsf{p}}{\sim}\mathcal{O''}, \ \mathcal{O''}\in\bar{\mathcal{B}}_{sp}^{(m)} \big{\}}, \\
\bar{\mathcal{B}}_{sg}^{(m)} & = \big{\{}\mathcal{O}\in\Ol(\mathcal{P}_m) \ \big{|} \ \mathcal{O}\stackrel{\mathsf{p}}{\sim}\mathcal{O'}\stackrel{\mathsf{w}}{\sim}\mathcal{O''}
\ \text{or} \ \mathcal{O}\stackrel{\mathsf{w}}{\sim}\mathcal{O'}\stackrel{\mathsf{p}}{\sim}\mathcal{O''}, \ \mathcal{O''}\in\bar{\mathcal{B}}_{sp}^{(m)}\big{\}}.
\end{align*}
It is easy to check that all four theorems, Theorem~\ref{tm:corr3-8} -- Theorem~\ref{tm:tm-spg-d}, hold with $\bar{\mathcal{B}}_{spg}^{(m)}$ and $\bar{\mathcal{B}}_{sg}^{(m)}$.

In the case $m=n$, one can reestablish the corresponding results for the standard Jacobi method and the associated Jacobi operators.

\section{Convergence of more general block Jacobi-type methods}

The obtained results for the block Jacobi operators and annihilators can be used to prove convergence of more general block Jacobi-type methods. This section is similar to~\cite[Section~5]{har-15}, and we will refer to some results from there. First, we prove the main result and then we apply it to the block J-Jacobi method from~\cite{har+sin-11}.

Consider the block Jacobi-type process
\begin{equation}\label{opciproces}
A^{(k+1)}=F_k^{T}A^{(k)}F_k, \quad k\geq0; \quad A^{(0)}=A,
\end{equation}
where $A$ is a symmetric matrix of order $n$, partitioned as in relation (\ref{blokmatrica}), and $F_k$, $k\geq0$, are elementary block matrices. Their pivot submatrices are only required to be nonsingular. Since all $F_k$ are nonsingular, $A\ne0$ implies $A^{(k)}\ne0$ for all $k$. The process is said to be Jacobi-type since it is generally not required that the pivot submatrices are diagonalized. We assume that for the process $(\ref{opciproces})$ the following assumptions hold.
\begin{itemize}
\item[\textbf{A1}] $\displaystyle \mathcal{O}\in\mathcal{B}_{sg}^{(m)}$, i.e., the pivot strategy $I_{\mathcal{O}}$ of the process is a generalized serial one.
\item[\textbf{A2}] There is a sequence of orthogonal elementary block matrices $U_k$, $k\geq0$, such that
$$\lim_{k\rightarrow\infty}(F_k-U_k)=0.$$
\item[\textbf{A3}] For the diagonal block $F_{ii}^{(k)}$ of $F_k$ one has
\[
\sigma =\liminf_{k\rightarrow\infty}\sigma^{(k)}>0, \quad\text{where} \ \sigma^{(k)}=\sigma_{\min}\big{(}F_{ii}^{(k)}\big{)}, \ k\geq 0.
\]
\end{itemize}
The first assumption \textbf{A1} deserves a comment. By Definition~\ref{tm:gss}, the set $\mathcal{B}_{sg}^{(m)}$ is defined using a single permutation equivalence. In order to make use of Theorems~\ref{tm:tm3-9} and~\ref{tm:tm-spg-d}, we can presume that in Definition~\ref{tm:gss} either $\mathcal{O'}=\mathcal{O}(\mathsf{q})$ or $\mathcal{O''}=\mathcal{O'}(\mathsf{q})$ holds for some permutation $\mathsf{q}$.

Because of the condition \textbf{A2}, in the assumption \textbf{A3} one can replace $F_{ii}^{(k)}$ by $U_{ii}^{(k)}$. From the CS decomposition of the orthogonal $(n_i+n_j)\times (n_i+n_j)$ matrix $\widehat{U}_{ij}$, we have $\sigma_{\min}\big{(}U_{ii}^{(k)}\big{)}=\sigma_{\min}\big{(}U_{jj}^{(k)}\big{)}$.
Therefore, in the definition of $\sigma^{(k)}$, instead of $F_{ii}^{(k)}$ one can use $F_{jj}^{(k)}$, $U_{ii}^{(k)}$ or $U_{jj}^{(k)}$. Recall that $F_{ii}^{(k)}$, $F_{jj}^{(k)}$, stands for $F_{i(k)i(k)}^{(k)}$, $F_{j(k)j(k)}^{(k)}$, respectively.

Since for each $U_k$, there is a permutation matrix $P_k$ that makes $U_k P_k$ a $\ubce$ matrix, the condition \textbf{A2} shows that for large enough $k$ each $F_k P_k$ will be arbitrary close to some $\ubce$ matrix. However, $U_k$, and therefore also $P_k$, is generally not known, while $F_k$ is available. Thus, one can perform the QR factorization with column pivoting of $[F_{ii}^{(k)} F_{ij}^{(k)}]$ to obtain $P_k$ and then replace $F_k$ with $F_k P_k$. The corresponding matrix $\widetilde{U}_k = U_k P_k$ may not be from $\ubce(1)$, but it is certainly from $\ubce(\varrho)$ for some $0<\varrho<1$ and when $k$ is large enough.

\begin{Theorem}\label{opcitm}
Let $\pi=(n_1,\ldots,n_m)$ be a partition of $n$, $\mathcal{O}\in\mathcal{B}_{sg}^{(m)}$ and let $A\in\mathbf{S}_{n}$, $A\ne 0$, be partitioned as in relation \eqref{blokmatrica}. Let the sequence of matrices $(A^{(k)};k\geq0)$ be generated by the block Jacobi-type process (\ref{opciproces}). If the assumptions $\mathbf{A1}-\mathbf{A3}$ are met, then the following two assertions are equivalent:
\begin{itemize}
\item[(i)] $\displaystyle \lim_{k\rightarrow\infty} \frac{S\big{(}\widehat{A}_{ij}^{(k+1)}\big{)}}{\|A^{(k)}\|_F}=0$, \quad where
    $\widehat{A}_{ij}^{(k+1)}=\widehat{F}_k^T\widehat{A}_{ij}^{(k)}\widehat{F}_k$;
\item[(ii)] $\displaystyle \lim_{k\rightarrow\infty} \frac{S(A^{(k)})}{\|A^{(k)}\|_F}=0$.
\end{itemize}
\end{Theorem}

\begin{proof}
The proof has been moved to the Appendix.
$\qquad$ \end{proof}

Thus, the condition (i) is sufficient for the convergence of $S(A^{(k)})$ to zero.
In the case of block-wise or element-wise Jacobi methods (i.e., those which diagonalize the pivot submatrix at each step) the condition (i) is trivially fulfilled. Note that
$S\big{(}\widehat{A}_{ij}^{(k+1)}\big{)}$ and $S(A^{(k)})$ are being divided by $\|A^{(k)}\|_F$, which is appropriate, since the theorem deals with nonorthogonal transformations.

In some applications the following corollaries can be used.

\begin{Corollary}
Theorem~\ref{opcitm} holds provided that the recurrence relation (\ref{opciproces}) is replaced with
\begin{equation}\label{eq2.38}
A^{(k+1)}= F_k^TA^{(k)}F_k + E^{(k)}, \quad k\geq 0,
\end{equation}
where $\lim_{k\rightarrow\infty} S(E^{(k)})/\|A^{(k)}\|_F=0$
and $E^{(k)} \not= -F_k^TA^{(k)}F_k$, $k\geq 0$.
The last condition on $E^{(k)}$ can be replaced by the requirement that $E^{(k)}=0$ whenever  $A^{(k)}=0$ for some $k$.
\end{Corollary}

\begin{proof}
Comparing with the proof of Theorem~\ref{opcitm}, the only difference appears in the definition of each vector $g^{(k)}$, now including the vector $e^{(k)}$, which in turn results from the matrix $E^{(k)}$.
$\qquad$ \end{proof}

\begin{Corollary} \label{tm:corr5.3}
Let $A\not=O$ be a matrix of order $n$ and let the sequence $A^{(0)}=A$, $A^{(1)},\ldots$ be generated by a block Jacobi-type process defined by relation (\ref{eq2.38}). Assume that the assumptions \textbf{A1} -- \textbf{A3} hold. Suppose that the sequence $(A^{(k)};k\geq 0)$ is bounded and
\begin{equation}\label{eq2.39}
\lim_{t\rightarrow\infty} S(E^{(k)}) = 0.
\end{equation}
Then the following two conditions are equivalent:
\begin{itemize}
\item[(iii)] ${\displaystyle \lim_{k\rightarrow\infty} S(\widehat{A}_{ij}^{(k+1)})=0}$;
\item[(iv)] ${\displaystyle \lim_{k\rightarrow\infty} S(A^{(k)})=0}$.
\end{itemize}
\end{Corollary}

\begin{proof}
The implication (iv) $\Rightarrow$ (iii) is obvious. For the converse implication, we use the expressions
$F_k=U_k+(F_k-U_k)$, $k\geq 0$, to transform the recursion (\ref{eq2.38}) into the form
$A^{(k+1)}=U_k^T A^{(k)}U_k+T^{(k)}+E^{(k)}$ with
\[
T^{(k)} = (F_k-U_k)^T A^{(k)}U_k+U_k^TA^{(k)}(F_k-U_k)+(F_k-U_k)^TA^{(k)}( F_k-U_k).
\]
By the boundedness of the sequence $(A^{(k)};k\geq 0)$ combined with the assumption \textbf{A2}, we have
\[
\|T^{(k)}\|\leq \sup \left\{ \|A^{(k)}\|_2; k\geq 0\right\}\cdot \left(2\|F_k-U_k\|+\|F_k-U_k\|^2\right) \rightarrow0 \quad \text{as} \ k\rightarrow\infty.
\]
This confirms relations (\ref{ojantmjedn}), (\ref{ojantmg}), with $g^{(k)}$ that additionally includes the vectors associated with the matrices $T^{(k)}$ and $E^{(k)}$.
The conditions (iii), (\ref{eq2.39}) and the latest relation together imply $\lim_{k\rightarrow\infty}g^{(k)} =0$.
Following the same lines of the proof of Theorem~\ref{opcitm}, one obtains
$\lim_{s\rightarrow\infty}g^{[s]}=0$ and $\lim_{k\rightarrow\infty}a^{(k)}=0$
if and only if $\lim_{s\rightarrow\infty}a^{[s]}=0$.
Therefore, under the conditions of this corollary, the sequence $(a^{[s]}; s\geq 1)$ has all the properties of the sequence (\ref{eq2.23}) from the proof of Theorem~\ref{opcitm}. The rest of the proof closely follows the proof of Theorem~\ref{opcitm}.
$\qquad$ \end{proof}

Let us note that the results in this section hold if the set of pivot orderings $\mathcal{B}_{sg}^{(m)}$ in the assumption \textbf{A1} is replaced with $\bar{\mathcal{B}}_{sg}^{(m)}$. Also, as has already been explained in~\cite{har-15}, it makes sense to rewrite the assumption \textbf{A3} in the equivalent form:
\begin{itemize}
\item[\textbf{A3}] For the diagonal block $F_{ii}^{(k)}$ of $F_k$ one has
\[
 \sigma = \liminf_{t\rightarrow\infty}\sigma^{[t]}>0, \quad \sigma^{[t]}=\min_{(t-1)T\leq k\leq tT-1} \sigma_{\min}(F_{ii}^{(k)}),
\]
where the quantities $\sigma^{[t]}$ are labeled by sweeps (i.e., cycles or quasi-cycles).
\end{itemize}

\subsection{An application to the block $J$-Jacobi methods}

The main purpose of Theorem~\ref{opcitm} is to be used in the global convergence considerations of the block Jacobi methods for the generalized eigenvalue problem, say for the HZ method from~\cite{nov+sin-15}. However, further research is needed to achieve this goal. Hence, we will choose yet another block method, which is well-understood, important in practice, and for which the newly obtained results can be applied  straightforwardly. It is the \emph{full block} $J$-\emph{Jacobi method} from \cite{har+sin-11}, for the pair $(A,J)$, where $A$ is symmetric positive definite and $J=\mbox{diag}(I_{\nu}, -I_{n-\nu})$. The main application of the method is to solving the simple eigenvalue problem $Hx=\lambda x$, where $H$ is indefinite symmetric matrix, with high relative accuracy. The partition $\pi=(n_1,\ldots ,n_m)$ has to comply with the partition $(\nu ,n-\nu)$, i.e., the first has to be a subpartition of the latter. After preliminary transformations, the problem $Hx=\lambda x$ is reduced to the generalized eigenvalue problem $Ax=\lambda Jx$. All details can be found in~\cite{sla-02,har+sin-11}. The method uses $J$-orthogonal elementary block matrices $F_k$, which leave $J$ intact under congruence transformations $F_k^TJF_k=J$, $k\geq 0$. The iteration process has the form (\ref{opciproces}) with a positive definite matrix $A=A^{(0)}$.

In~\cite{har+sin-11} the global convergence of this method was proved under the weak-wavefront strategies and here we prove it for the much larger class of generalized serial strategies.

\begin{Theorem}\label{tm:last}
Let $\pi=(n_1,\ldots,n_m)$ be a partition of $n$, so that $\pi$ refines $(\nu,n-\nu)$.
The full block $J$-Jacobi method defined by the cyclic pivot strategy $I_{\mathcal O}$, $\mathcal{O}\in\mathcal{B}_{sg}^{(m)}$, which uses $\ubce$ $J$-orthogonal transformation matrices is globally convergent.
\end{Theorem}

\begin{proof}
Since the iterates generated by the full block $J$-Jacobi method are bounded \cite[(3.18)]{har+sin-11}, we can apply
Corollary~\ref{tm:corr5.3} with the matrices $E^{(k)}$, $k\geq 0$, set to zero.
The method is called the full block $J$-Jacobi method, because at each step it diagonalizes the pivot submatrix. This implies that condition (iii) of Corollary~\ref{tm:corr5.3} is fulfilled.

All we have to do is follow the same lines of the proof of~\cite[Proposition~3.3]{har+sin-11}, which in turn reduces to checking validity of the assumptions \textbf{A1} -- \textbf{A3} from Theorem~\ref{opcitm}.

The first assumption is presumed. The second one follows from~\cite[Proposition~3.2]{har+sin-11}. Assumption \textbf{A3} holds for two reasons. First, for each hyperbolic elementary block transformation $F_k$ one has $\sigma^{(k)}\geq 1$, and we only have to check \textbf{A3} for the orthogonal elementary block transformations. However, they are exactly the same as those in the block Jacobi method for symmetric matrices from Sections 3 and 4 in this paper. Relations (\ref{uv_1}) and (\ref{uv_1a}) hold for them, even with $\varrho =1$. Since condition (iii) of Corollary~\ref{tm:corr5.3} is fulfilled, we have $S(A^{(k)})\rightarrow 0$ as $k\rightarrow\infty$. The rest of the proof requires an analogue of Theorem~\ref{tm:gl_conv-opp} for the $J$-Jacobi method. However, all that is needed in the proof is an estimate similar to (\ref{eigs3}) for the diagonal elements of $A^{(k)}$. Such an estimate is established in~\cite[Lemma~1.1]{Drmac-Hari-93}.
$\qquad$ \end{proof}

By using the results from Section~4 one can easily show that Theorem~\ref{tm:last} holds for any quasi-cyclic strategy $I_{\mathcal{O}}$, $\mathcal{O}\in\bar{\mathcal{B}}_{sg}^{(m)}$.

\section{Conclusion and future work}

So far, a satisfactory research goal has been the global convergence of the block Jacobi method for symmetric matrices, established for the serial pivot strategies or those that are equivalent or weak equivalent to them, so-called wavefront or weak wavefront strategies. All those strategies were obtained from a single cyclic strategy, say the column-cyclic one. Here we have shown how to further enlarge the class of convergent strategies by using the notion of the reverse strategy and that of permutation equivalent strategies. Hence, with each convergent pivot strategy we have associated the whole large class of convergent strategies, obtained from it by using four equivalence relations, $\sim$, $\stackrel{s}{\sim}$, $\stackrel{w}{\sim}$, $\stackrel{\mathsf{p}}{\sim}$, and by using reverse strategies. The next step was to increase the number of classes of convergent strategies obtained this way. For large $m$, we have obtained at least $2!3!\cdots (m-1)!$ such classes of convergent strategies and we have named their union, the class of generalized serial strategies. Furthermore, the convergence results for that class are stated and proved in the stronger form, which enables us to formulate and prove similar results for the block Jacobi operators. This makes the block Jacobi operators a tool for proving the global convergence of the block Jacobi methods for other eigenvalue problems, in particular for the generalized eigenvalue problem. As an immediate result, we have proved the global convergence of the full block $J$-Jacobi method under any generalized serial pivot strategy.

Future work will include proving the global convergence of the (block-wise and element-wise) HZ method~\cite{har-84,nov+sin-15} for the generalized eigenvalue and singular value problem under the class of generalized serial strategies. We also intend to prove the global convergence of the block Paardekooper method for skew-symmetric matrices. An immediate consequence of the results from this paper is the proof that in the case $n=4$ all $720$ cyclic strategies for the symmetric Jacobi method are convergent~\cite{beg+har15a}. The future research will also be concentrated on the complex block Jacobi methods, first for a single Hermitian matrix and then for a positive definite pair of Hermitian matrices (the complex block $J$-Jacobi and the complex HZ method).

\section*{Acknowledgments}
The authors are indebted to the referees and the editor for their valuable suggestions that improved the paper.

\appendix
\setcounter{equation}{0}

\section{Proofs omitted in the main text}

\subsection{Proof of Theorem~\ref{tm:gl_conv-opp}}

Suppose that $A$ is not a multiple of $I_n$ and assume that for the eigenvalues of $A$ we have
$$\lambda_1=\cdots =\lambda_{\mathbf{s}_1} > \lambda_{\mathbf{s}_1+1}=\cdots =\lambda_{\mathbf{s}_2} > \cdots > \lambda_{\mathbf{s}_{\omega -1}+1}=\cdots =\lambda_{\mathbf{s}_{\omega}}, \quad \mathbf{s}_{\omega}=n,$$
$\mathbf{s}_r=\nu_1+\cdots +\nu_r$, $1\leq r\leq \omega$, where $\nu_1$,\ldots ,$\nu_{\omega}$ are the multiplicities of the eigenvalues. Let
$$3\delta = \min_{1\leq r\leq \omega-1} (\lambda_{\mathbf{s}_r} - \lambda_{\mathbf{s}_{r+1}}).$$
Obviously, there is an integer $k_0$ such that
$$S(A^{(k)})<\delta, \quad k\geq k_0.$$
From \cite[Lemma~2.1]{har-91} we conclude that for $k\geq k_0$ all diagonal elements of $A^{(k)}$ lie in the union of shrinking segments
\[
{\mathcal{D}}_r^{(k)} \equiv \big{\{} t \ \big{|} \ |t-\lambda_{\mathbf{s}_r}|\leq 0.22S(A^{(k)}) \big{\}} \subset
\big{\{} t \ \big{|} \ |t-\lambda_{\mathbf{s}_r}|\leq 0.22\delta\}\equiv {\mathcal{D}}_r, \quad 1\leq r\leq \omega.
\]
Furthermore, by the same lemma, each ${\mathcal{D}}_r^{(k)}$ contains at least $\nu_r$ diagonal elements of $A^{(k)}$ and then it straightforwardly follows that ${\mathcal{D}}_r^{(k)}$ contains exactly $\nu_r$ diagonal elements of $A^{(k)}$. In particular, this implies that for any two diagonal elements of $A^{(k)}$ we have
\begin{equation}\label{eigs3}
\text{either} \ |a_{ll}^{(k)}-a_{mm}^{(k)}| \leq0.44\delta, \qquad \text{or} \ |a_{ll}^{(k)}-a_{mm}^{(k)}|>2.56\delta, \quad k\geq k_0.
\end{equation}
Since the sequence $S(A^{(k)})$, $k\geq 0$, converges to zero, the proof will be completed if we show that for $k\geq k_0$ the diagonal elements cannot change their affiliation to eigenvalues. Afterwards, we will also show how the diagonal elements of $\Lambda$ are ordered along the diagonal.

To establish assertion (i) it is sufficient to prove the first claim only since the second one can be proved in a similar way. Even if the blocks $A_{ii}^{(k)}$, $1\leq i\leq m$, were not diagonalized at the beginning, we can increase $k_0$, if needed, so that the assumption of the first claim reads: each $A_{ii}^{(k)}$, $1\leq i\leq m$, $k\geq k_0$, is diagonal with diagonal elements ordered nonincreasingly. To this end we denote $A_{ii}^{(k)}$ by  $\Lambda_{ii}^{(k)}$, $1\leq i\leq m$. We can also assume that $k_0= t_0 T$, where $T\geq M=m(m-1)/2$ is the period of the strategy. The proof will be completed if we can find $k_0' \geq k_0$ such that
\begin{equation}\label{eigs4}
\{a_{11}^{(k)},a_{22}^{(k)},\ldots,a_{\mathbf{s}_r,\mathbf{s}_r}^{(k)}\} \subset {\bigcup_{p=1}^r} {\mathcal{D}}_p, \quad 1\leq r\leq\omega, \ k\geq k_0'.
\end{equation}
Let us consider step $k$ of the block method with $k\geq k_0$. Let $i=i(k)$, $j=j(k)$, and let
$\widehat{\widehat{A}}_{ij}^{(k)} = \text{diag}(\Lambda_{ii}^{(k+1)},\Lambda_{jj}^{(k+1)})$ be the transformed pivot submatrix $\widehat{A}_{ij}^{(k)}$.
From relation (\ref{blokpivotmatr}) and the perturbation theorem for the symmetric matrices we conclude that
\begin{align}
\label{eigs5} \|\text{diag}(\Lambda_{ii}^{(k+1)},\Lambda_{jj}^{(k+1)})- P_k^T \text{diag}(\Lambda_{ii}^{(k)},\Lambda_{jj}^{(k)})P_k\|_2 & \leq \|A_{ij}^{(k)}\|_2 \\
& \leq \frac{\sqrt{2}}{2}S(A^{(k)}) < \frac{\sqrt{2}}{2}\delta \nonumber
\end{align}
holds for $k\geq k_0$, where $P_k$ are some permutation matrices.
From relations (\ref{eigs5}) and (\ref{eigs3}) we obtain the following geometric interpretation of the movement of the diagonal elements during one step of the method.
The diagonal elements of $A^{(k)}$ are points on the real axis situated within small segments around the eigenvalues. After the completion of step $k$ they have moved (as points) within the same segment, but (as diagonal elements) their subscripts may have changed.

What happens with a diagonal element $a_{qq}^{(k_0)}$, which lies in the segment ${\mathcal{D}}_1$?

First, suppose that $\nu_1 \leq n_p$ for all  $1\leq p\leq m$. Then each time $a_{qq}^{(k_0)}$ is affected, it will be a diagonal element of $\Lambda_{ii}^{(k+1)}$. Let $a_{qq}^{(k_0)}$ lie in $A_{ll}^{(k_0)}$. Then, for $k\geq k_0$, $a_{qq}^{(k)}$ is affected when $i(k)=l$ or $j(k)=l$.

If $i(k)=l$, then $a_{qq}^{(k)}$ will remain in the same diagonal block, which is $\Lambda_{ii}^{(k+1)}=A_{ll}^{(k+1)}$. If $\Lambda_{jj}^{(k)}$ contains some diagonal
elements from ${\mathcal{D}}_1$, then they will move to $\Lambda_{ii}^{(k+1)}$ and thus the number of diagonal elements from ${\mathcal{D}}_1$ in  $A_{ll}^{(k+1)}$
will be larger than in $A_{ll}^{(k)}$.

If $j(k)=l$, then $a_{qq}^{(k)}$ will move to the new diagonal block $\Lambda_{ii}^{(k+1)}$, $i<l$, hence its subscripts will become smaller than or equal to $\mathbf{s}_{l-1}$. Since the pivot strategy is cyclic or quasi-cyclic, the case $j(k)=l$ must occur within the current sweep, unless $l=1$. Hence, during the next sweep $a_{qq}^{(k_0)}$ will move to some diagonal block which lies closer to $A_{11}^{(k_0)}$, unless $l=1$.

Since $a_{qq}^{(k_0)}$ is an arbitrary diagonal element of ${\mathcal{D}}_1$, we conclude that within one sweep all diagonal elements from ${\mathcal{D}}_1$ not belonging to $A_{11}^{(k)}$ will decrease their subscripts to such an extent that they become the diagonal elements of other diagonal blocks. This analysis shows that within the first $m-1$ sweeps all diagonal elements belonging to ${\mathcal{D}}_1$ will be the elements of the first diagonal block.

Now, let $\nu_1$ be such that $1\leq\nu_1<n$ holds. Then the same analysis shows that during  $m-1$ sweeps the diagonal elements affiliated with $\lambda_1$ will be filling in the first $\nu_1$ diagonal positions of the matrix. Hence, there is a number $k_1\geq (m-1)T+k_0$ such that the first $\nu_1$ diagonal elements in $A^{(k_1)}$ are affiliated with $\lambda_1$.

Almost the same analysis shows that within the first $m-1$ sweeps the diagonal elements affiliated with $\lambda_n$ will be filling in the last $\nu_{\omega}$ diagonal positions of the matrix. By increasing $k_1$ if necessary, we can assume that the last $\nu_{\omega}$ diagonal elements of $A^{(k_1)}$ are affiliated with $\lambda_n$.

The rest of the proof considers the matrix $A^{(k_1)}$. In $A^{(k_1)}$ the first $\nu_1$  and the last $\nu_{\omega}$ diagonal positions are occupied by the diagonal elements from ${\mathcal{D}}_1$ and ${\mathcal{D}}_{\omega}$, respectively.  The situation is described by the following block-matrix partition
\[
A^{(k_1)} = \left[
\begin{array}{ccc}
A_{s_{\varrho}}^{(k_1)}& B^{(k_1)} &C^{(k_1)}\\
{B^{(k_1)}}^T &A_{n-s_{\varrho}-\widetilde{s}_{\rho}}^{(k_1)} & G^{(k_1)}\\
{C^{(k_1)}}^T & {G^{(k_1)}}^T & A_{\widetilde{s}_{\rho}}^{(k_1)}
\end{array}
\right],\ \quad \begin{array}{l}
s_{\varrho}=n_1+\cdots +n_{\varrho},\\[1ex]
\widetilde{s}_{\rho} = n_m+\cdots +n_{m-\rho+1},
\end{array}
\]
where $s_{\varrho}\leq  \nu_1<s_{\varrho +1}$ and $\widetilde{s}_{\rho}\leq \nu_{\omega}<\widetilde{s}_{\rho +1}$. In this situation, if the pivot blocks are within $A_{s_{\varrho}}^{(k)}$, $B^{(k)}$, $C^{(k)}$, $G^{(k)}$, $A_{\widetilde{s}_{\rho}}^{(k)}$, $k\geq k_1$, the corresponding steps will make no subscript change in the diagonal elements, or the change will only mean repositions within the same diagonal block. Therefore, our analysis will only consider the central block $A_{n-s_{\varrho}-\widetilde{s}_{\rho}}^{(k_1)}$. The diagonal elements of $A_{n-s_{\varrho}-\widetilde{s}_{\rho}}^{(k_1)}$ from ${\mathcal{D}}_1$ (resp.\@ ${\mathcal{D}}_{\omega}$), if there are any, have already settled, within the first (resp.\@ last) positions of $A_{\varrho+1,\varrho+1}^{(k_1)}$ (resp.\@ $A_{m-\rho,m-\rho}^{(k_1)}$). They will not leave these positions during the next steps.

Then after the following $m-\varrho-\rho-1$ or more sweeps, the diagonal elements from ${\mathcal{D}}_2$ and ${\mathcal{D}}_{\omega-1}$ will settle in. Continuing this consideration we finally obtain the matrix $A^{(k')}$, $k'>k_0$, for which relation (\ref{eigs4}) holds. We note that $k'$ depends on the pivot strategy. For the serial ones, say for the row-cyclic one, the above analysis shows that we can take $k'=k_0+M$.

In order to prove (ii), we consider the diagonalization of the pivot submatrix, which is generally described by relation (\ref{blokpivotmatr}). Since we use a globally convergent element-wise Jacobi method, we know that the off-norm sequence of this submatrix tends to zero. In~\cite{mas-90} it was proved that the diagonal elements always converge. Thus, the limit $\mbox{diag}(\Lambda_{ii}^{(k+1)},\Lambda_{jj}^{(k+1)})$ exists. Finally, it is known that the diagonal elements are updated by $\pm\tan\phi_k a_{lm}^{(k)}$, where $(l,m)$ is the pivot pair. Therefore, the change is smaller then $1\cdot |a_{lm}^{(k)}|\leq \|A_{ij}^{(k)}\|_2 < \sqrt{2}/2\,\delta$. Hence, the diagonal elements cannot change their affiliation to the eigenvalues. This means that the permutation $P_k$ from relation (\ref{eigs5}) can be taken to be identity.

\subsection{Proof of Corollary~\ref{tm:cor_bjann}}

Set $\bm{\mathcal{R}} = \bm{\mathcal{R}}(\widehat{U})$, where $\widehat{U}$ is as in (\ref{piv_mat_E}).
It is sufficient to check that the transpose of each of the three types of submatrices appearing in Theorem~\ref{tm:bl_jac_ann} is of the same type and possesses the same properties. For the first and the third type, the proof is straightforward:
\begin{align*}
\left[\begin{array}{cc}
       U_{ii}^T\otimes I_{n_r} & U_{ji}^T\otimes I_{n_r} \\
       U_{ij}^T\otimes I_{n_r} & U_{jj}^T\otimes I_{n_r}
\end{array}\right]^T & =
\left[\begin{array}{cc}
       U_{ii}\otimes I_{n_r} & U_{ij}\otimes I_{n_r} \\
       U_{ji}\otimes I_{n_r} & U_{jj}\otimes I_{n_r}
\end{array}\right]
=
\left[\begin{array}{cc}
       V_{ii}^T\otimes I_{n_r} & V_{ji}^T\otimes I_{n_r} \\
       V_{ij}^T\otimes I_{n_r} & V_{jj}^T\otimes I_{n_r}
\end{array}\right], \\
\left[\begin{array}{cc}
       I_{n_r}\otimes U_{ii}^T & I_{n_r}\otimes U_{ji}^T \\
       I_{n_r}\otimes U_{ij}^T & I_{n_r}\otimes U_{jj}^T
\end{array}\right]^T & =
\left[\begin{array}{cc}
       I_{n_r}\otimes U_{ii} & I_{n_r}\otimes U_{ij} \\
       I_{n_r}\otimes U_{ji} & I_{n_r}\otimes U_{jj}
\end{array}\right]
   =
\left[\begin{array}{cc}
       I_{n_r}\otimes V_{ii}^T & I_{n_r}\otimes V_{ji}^T \\
       I_{n_r}\otimes V_{ij}^T & I_{n_r}\otimes V_{jj}^T
\end{array}\right],
\end{align*}
with
\[
V=\left[\begin{array}{cc}
        V_{ii} &  V_{ij}\\
        V_{ji} & V_{jj}
\end{array}\right] =U^T.
\]
Note that $V$ has the same essential properties as $U$: dimension and orthogonality (belonging to $\ubce(\varrho)$).
For the second type of submatrices we have
\begin{align*}
\left[\begin{array}{cc}
       I_{n_r}\otimes U_{ii}^T & S(U_{ji}^T\otimes I_{n_r}) \\
       \widetilde{S}(I_{n_r}\otimes U_{ij}^T) & U_{jj}^T\otimes I_{n_r}
\end{array}\right]^T & =
\left[\begin{array}{cc}
       I_{n_r}\otimes U_{ii} & (I_{n_r}\otimes U_{ij})\widetilde{S}^T \\
     (U_{ji}\otimes I_{n_r})S^T   & U_{jj}\otimes I_{n_r}
\end{array}\right]\\
& = \left[\begin{array}{cc}
       I_{n_r}\otimes V_{ii}^T & S(V_{ji}^T \otimes I_{n_r}) \\
     \widetilde{S}(I_{n_r}\otimes V_{ij}^T)  & V_{jj}^T\otimes I_{n_r}
\end{array}\right].
\end{align*}
To prove the second equality we need some extra work. It is obvious that this equality holds for the corresponding diagonal blocks.
To prove that the corresponding $(1,2)$ blocks are equal, let $U_{ij}=(u_{st})$ and note that it is an $n_i\times n_j$ matrix. If $e_k^T$ denotes the $k$th row of $I_{n_r}$, then $U_{ij} (I_{n_j}\otimes e_k^T)$ is an $n_i\times n_j n_r$ matrix and we have
\begin{align*}
U_{ij} (I_{n_j}\otimes e_k^T) & =
\left[\begin{array}{ccc} e_k^T & &\\ & \ddots & \\ & & e_k^T\end{array}\right]
\left[\begin{array}{ccc} u_{11}I_{n_r}&\cdots & u_{1n_j}I_{n_r}\\ &\ddots & \\ u_{n_i 1}I_{n_r} &\cdots & u_{n_i n_j}I_{n_r}\end{array}\right] \\
& = (I_{n_i}\otimes e_k^T)\,(U_{ij}\otimes I_{n_r})
= (I_{n_i}\otimes e_k^T)\,(V_{ji}^T\otimes I_{n_r}),\quad 1\leq k\leq n_r.
\end{align*}
Hence,
$$(I_{n_r}\otimes U_{ij})\widetilde{S}^T  = \left[\begin{array}{c}U_{ij}(I_{n_j}\otimes e_1^T) \\ \vdots \\ U_{ij}(I_{n_j}\otimes e_{n_r}^T) \end{array}\right]
= \left[\begin{array}{c} (I_{n_i}\otimes e_1^T)\,(V_{ji}^T\otimes I_{n_r}) \\ \vdots \\ (I_{n_i}\otimes e_{n_r}^T)\,(V_{ji}^T\otimes I_{n_r}) \end{array}\right]
= S (V_{ji}^T\otimes I_{n_r}).$$
Next, let us prove that the corresponding $(2,1)$ blocks are equal. Note that $U_{ji}=V_{ij}^T$ is an $n_j\times n_i$ matrix. If $e_k$ denotes the $k$th column of $I_{n_r}$, then for each $1\leq k\leq n_r$ we have
\[
(U_{ji}\otimes I_{n_r}) (I_{n_i}\otimes e_k)= U_{ji}I_{n_i}\otimes I_{n_r}e_k= U_{ji}\otimes e_k= (I_{n_j}\otimes e_k) (U_{ji}\otimes I_1) =(I_{n_j}\otimes e_k) U_{ji}
\]
and all $(U_{ji}\otimes I_{n_r}) (I_{n_i}\otimes e_k)$ are $n_j n_r\times n_i$ matrices. Now we have
\begin{align*}
(U_{ji}\otimes I_{n_r})S^T & = (U_{ji}\otimes I_{n_r})[I_{n_i}\otimes e_1 \cdots I_{n_i}\otimes e_{n_r} ] \\
& = [(U_{ji}\otimes I_{n_r})(I_{n_i}\otimes e_1) \cdots (U_{ji}\otimes I_{n_r})(I_{n_i}\otimes e_{n_r}) ] \\
& = [(I_{n_j}\otimes e_1) U_{ji} \cdots (I_{n_j}\otimes e_{n_r}) U_{ji}] \\
& = [I_{n_j}\otimes e_1 \cdots I_{n_j}\otimes e_{n_r}] \left[\begin{array}{ccc}U_{ji} & &\\ &\ddots &\\ & & U_{ji}\end{array}\right] \\
& = \widetilde{S}(I_{n_r}\otimes U_{ji}) = \widetilde{S}(I_{n_r}\otimes V_{ij}^T).
\end{align*}

\subsection{Proof of Theorem~\ref{tm:b1}}

First we list some inequalities that will be used.
Let $\sigma_{min}(X)$ and $\sigma_{max}(X)$ ($=\|X\|_2$) denote the smallest and largest singular value of $X$. Recall that $\|X\|_F$ denotes the Frobenius norm. Let $\|X\|$ denote any matrix norm. We have
\begin{align}
\label{bn1} \big{|}\|X\|-\|Y\|\big{|} & \leq \|X+Y\| \leq \|X\|+\|Y\|, \\
\nonumber \|X_1+\cdots+X_t\| & \leq\|X_1\|+\cdots+\|X_t\|, \quad t\geq1, \\
\nonumber \max\{\sigma_{min}(F)\|G\|_F,\sigma_{min}(G)\|F\|_F\} & \leq \|FG\|_F \leq \min\{\sigma_{max}(F)\|G\|_F,\sigma_{max}(G)\|F\|_F\}, \\
\nonumber \sigma_{min}(X_1\ldots X_t) & \geq\sigma_{min}(X_1)\cdots\sigma_{min}(X_t), \quad t\geq1, \\
\label{bn5} \sigma_{max}(X_1\ldots X_t) & \leq\sigma_{max}(X_1)\cdots\sigma_{max}(X_t), \quad t\geq1.
\end{align}
Besides, if $X=(X_{rs})$ is a block matrix as in relation (\ref{blokmatrica}), then both for the operator matrix norm and also for the Frobenius norm we have
\begin{equation}\label{bloknorma}
\|X_{rs}\|\leq\|X\|,\qquad 1\leq r,s\leq m.
\end{equation}
To prove Theorem~\ref{tm:b1jop}, we start with the partition $\pi=(n_1,\ldots,n_m)$ and denote by $\pi_l=(n_1,\ldots,n_l)$ the partition of $s_l=n_1+\cdots +n_l$. Obviously, for $l=m$ we have $\pi_m=\pi$ and $s_m=n$.
The set associated with $\pi_l$ is $\mathcal{B}_c^{(l)}$ from (\ref{b1def}) where $m$ is replaced by $l$.

We will prove the following statement.
Let $l\in\{2,\ldots,m\}$ and let $A$ be any symmetric matrix  of order $s_l$, carrying the block-matrix partition defined by $\pi_l$. Apply to $A$ the cyclic block Jacobi method defined by the pivot strategy $I_{\mathcal{O}}$, $\mathcal{O}\in\mathcal{B}_c^{(l)}$, thus obtaining the symmetric matrices $A^{(0)}=A$, $A^{(1)},\ldots$ defined by the recursion (\ref{jacobiblokagm}). If the transformation matrices are from the class $\ubce_{\pi_l}(\varrho)$, then
\begin{equation}\label{b1t}
S^2(A^{(L)})\leq\eta_{\pi_l,\varrho}S^2(A), \quad 0\leq \eta_{\pi_l,\varrho}<\widetilde{\eta}_{s_l,\varrho}<1, \quad L=\frac{l(l-1)}{2},
\end{equation}
where the constants $\eta_{\pi_l,\varrho}$ and $\widetilde{\eta}_{s_l,\varrho}$ in~\eqref{b1t} depend only on $\pi_l$, $\varrho$ and $s_l$, $\varrho$, respectively.

Obviously, for $l=m$ we obtain the assertion of Theorem~\ref{tm:b1jop}.
The proof of (\ref{b1t}) uses mathematical induction on $l$, $2\leq l\leq m$.

For $l=2$,  $A$ is of order $s_2=n_1+n_2$. Its only pivot block is $A_{12}$, so that $\widehat{A}=A$.
One step of the block Jacobi method is needed to diagonalize $A$. We have $S^2(A^{(1)})=0$, so $\eta_{\pi_2,\varrho}=\widetilde{\eta}_{s_2,\varrho}=0$.

Assume that assertion~\eqref{b1t} holds for $l-1$, $l\in\{3,\ldots,m\}$, and for the partition $\pi_{l-1}$, with constants $0\leq\eta_{\pi_{l-1},\varrho}\leq\widetilde{\eta}_{s_{l-1},\varrho}<1$. In the induction step, we will prove that~\eqref{b1t} holds for $l$.

Set $A=(A_{rs})$ be a symmetric matrix of order $s_l$, partitioned according to $\pi_l$.
For an arbitrary ordering $\mathcal{O}$ from $\mathcal{B}_c^{(l)}$ apply the cyclic block Jacobi method defined by $I_{\mathcal{O}}$. Let the transformation matrices be from the class $\ubce_{\pi_l}(\varrho)$ and let the obtained sequence of matrices be denoted by $A^{(0)}=A,A^{(1)},\ldots$

Let $\widetilde{L}=(l-1)(l-2)/2$. Let $\widetilde{A}$ ($A_{l-1}$) be the leading submatrix of $A^{(\widetilde{L})}$ ($A$) of order $s_{l-1}$. In other words, $\widetilde{A}$ is obtained from $A_{l-1}$ after completing one full sweep of $\widetilde{L}$ Jacobi steps. During these steps, the last, $l$th block-column of $A$, has been affected only by the left transformations. Therefore, we have
\[
\sum_{i=1}^{l-1}\|\widetilde{A}_{il}\|^2_F=\sum_{i=1}^{l-1}\|A_{il}\|^2_F,
\]
where $A^{(\widetilde{L})}=(\widetilde{A}_{rs})$. Let $0\leq\epsilon\leq1$ be such that
\begin{equation}\label{epsilon}
(1-\epsilon^2)S^2(A)=\sum_{i=1}^{l-1}\|A_{il}\|^2_F.
\end{equation}
Thus, $S(A_{l-1})=\epsilon S(A)$.
The submatrix $A_{l-1}$ is of order $s_{l-1}$ and carries the block-matrix partition defined by $\pi_{l-1}$. Therefore, the induction hypothesis can be applied. It follows that
\begin{align} \label{b1mi}
S^2(A^{(\widetilde{L})}) & = S^2(\widetilde{A}) + \sum_{i=1}^{l-1}\|\widetilde{A}_{il}\|^2_F \leq \eta_{\pi_{l-1},\varrho} S^2(A_{l-1}) + \sum_{i=1}^{l-1}\|\widetilde{A}_{il}\|^2_F \\
\nonumber & = \eta_{\pi_{l-1},\varrho}\epsilon^2 S^2(A)+(1-\epsilon^2)S^2(A) \\
\nonumber & = \left(1-\epsilon^2(1-\eta_{\pi_{l-1},\varrho})\right)S^2(A).
\end{align}
Even thought we have $S^2(A^{(L)})\leq S^2(A^{(\widetilde{L})})$, we cannot set $\eta_{\pi_l,\varrho}=1-\epsilon^2(1-\eta_{\pi_{l-1}},\varrho)$ because $\epsilon$ can be arbitrarily small or zero. We still need to estimate the contribution to the off-norm reduction coming from the last $l-1$ steps.

According to relation (\ref{b1def}), the blocks in the $l$th block-column are annihilated in the order: $(\tau_l(1),l),\ldots,(\tau_l(l-1),l)$. Let us consider how the block $\tau_l(i)$ changes until it is annihilated in the $i$th step. To simplify notation in this analysis, we write $\tau $ instead of $\tau_l$ (the permutation of the set $\{1,\ldots ,l-1\}$) until relation (\ref{feps}). We have
\begin{align*}
\widetilde{A}_{\tau(i)l}^{(1)} & = \widetilde{A}_{\tau(i)l}U_{ll}^{(\widetilde{L})}+\widetilde{A}_{\tau(i)\tau(1)}U_{\tau(1)l}^{(\widetilde{L})},\\
\widetilde{A}_{\tau(i)l}^{(2)} & = \widetilde{A}_{\tau(i)l}^{(1)}U_{ll}^{(\widetilde{L}+1)}+
\widetilde{A}_{\tau(i)\tau(2)}U_{\tau(2)l}^{(\widetilde{L}+1)},\\
& \vdots \\
\widetilde{A}_{\tau(i)l}^{(i-1)} & = \widetilde{A}_{\tau(i)l}^{(i-2)}U_{ll}^{(\widetilde{L}+i-2)}+
\widetilde{A}_{\tau(i)\tau(i-1)}U_{\tau(i-1)l}^{(\widetilde{L}+i-2)},\\
\widetilde{A}_{\tau(i)l}^{(i)} & = 0.
\end{align*}
The contribution to the off-norm reduction comes from $\|\widetilde{A}_{\tau(i)l}^{(i-1)}\|_F$, so we have to estimate it from below. To express $\widetilde{A}_{\tau(i)l}^{(i-1)}$ in terms of the blocks from $\widetilde{A}$, we multiply the equation for $\widetilde{A}_{\tau(i)l}^{(i-2)}$ from the right by $U_{ll}^{(\widetilde{L}+i-2)}$, then multiply the equation for $\widetilde{A}_{\tau(i)l}^{(i-3)}$ by $U_{ll}^{(\widetilde{L}+i-3)}U_{ll}^{(\widetilde{L}+i-2)}$, etc. Finally, we multiply the first equation by $U_{ll}^{(\widetilde{L}+1)}\cdots U_{ll}^{(\widetilde{L}+i-2)}$ from the right. Then we take the sum of the obtained equations. It follows that
\[
\widetilde{A}_{\tau(i)l}^{(i-1)} = \widetilde{A}_{\tau(i)l}U_{ll}^{(\widetilde{L})}U_{ll}^{(\widetilde{L}+1)}\cdots U_{ll}^{(\widetilde{L}+i-2)} + \sum_{k=1}^{i-1}\widetilde{A}_{\tau(i)\tau(k)}U_{\tau(k)l}^{(\widetilde{L}+k-1)}U_{ll}^{(\widetilde{L}+k)}\cdots U_{ll}^{(\widetilde{L}+i-2)}.
\]
Denote
\begin{equation}\label{varsigma}
\zeta_l=\min_{0\leq k\leq l-3}\big{\{}\sigma_{\min}(U_{ll}^{(\widetilde{L}+k)})\big{\}}.
\end{equation}
Using inequalities~\eqref{bn1} -- \eqref{bn5}, for $1\leq i\leq l-1$, we obtain
\begin{eqnarray*}
\lefteqn{
\|\widetilde{A}_{\tau(i)l}^{(i-1)}\|_F  \geq \ \left| \|\widetilde{A}_{\tau(i)l}U_{ll}^{(\widetilde{L})}\cdots U_{ll}^{(\widetilde{L}+i-2)}\|_F - \|\sum_{k=1}^{i-1}\widetilde{A}_{\tau(i)\tau(k)}U_{\tau(k)l}^{(\widetilde{L}+k-1)}\cdots U_{ll}^{(\widetilde{L}+i-2)}\|_F \ \right| } \\
&& \geq \sigma_{\min}\big{(}U_{ll}^{(\widetilde{L})}\cdots U_{ll}^{(\widetilde{L}+i-2)}\big{)}\|\widetilde{A}_{\tau(i)l}\|_F - \sum_{k=1}^{i-1}\sigma_{\max}\big{(}U_{\tau(k)l}^{(\widetilde{L}+k-1)}\cdots U_{ll}^{(\widetilde{L}+i-2)}\big{)}\|\widetilde{A}_{\tau(i)\tau(k)}\|_F \\
&& \geq \sigma_{\min}\big{(}U_{ll}^{(\widetilde{L})}\big{)}\sigma_{\min}\big{(}U_{ll}^{(\widetilde{L}+1)}\big{)}\cdots \sigma_{\min}\big{(}U_{ll}^{(\widetilde{L}+i-2)}\big{)}\|\widetilde{A}_{\tau(i)l}\|_F - \\
&& \quad - \sum_{k=1}^{i-1}\sigma_{\max}\big{(}U_{\tau(k),l}^{(\widetilde{L}+k-1)}\big{)} \sigma_{\max}\big{(}U_{ll}^{(\widetilde{L}+k)}\big{)}\cdots \sigma_{\max}\big{(}U_{ll}^{(\widetilde{L}+i-2)}\big{)}\|\widetilde{A}_{\tau(i)\tau(k)}\|_F \\
&& \geq \zeta_l^{i-1} \|\widetilde{A}_{\tau(i),l}\|_F - \sum_{k=1}^{i-1}\|\widetilde{A}_{\tau(i)\tau(k)}\|_F.
\end{eqnarray*}
Here, we have used (\ref{bloknorma}) for the transformation matrices, which are orthogonal.
Squaring the obtained inequality and then using
$(a-b)^2\geq \frac{1}{2}a^2-b^2$, $a,b\in\mathbb{R}$,
and the Cauchy--Schwarz inequality, we obtain
\begin{equation}\label{here_1}
\|\widetilde{A}_{\tau(i)l}^{(i-1)}\|_F^2 \geq \frac{1}{2}\zeta_l^{2(i-1)} \|\widetilde{A}_{\tau(i)l}\|_F^2 - (i-1)\sum_{k=1}^{i-1}\|\widetilde{A}_{\tau(i)\tau(k)}\|_F^2.
\end{equation}
Now we have the lower bound for the reduction of $S^2(A^{(\widetilde{L})})$ coming from just one annihilated block in the last block-column. The lower bound coming from all blocks in the last block-column is obtained by summing up these. Using relations (\ref{epsilon}), (\ref{b1mi}) and (\ref{here_1}), we have
\begin{eqnarray}
\sum_{i=1}^{l-1} \|\widetilde{A}_{\tau(i)l}^{(i-1)}\|_F^2 & \geq&
\frac{1}{2}\sum_{i=1}^{l-1} \ \zeta_l^{\ 2(i-1)} \|\widetilde{A}_{\tau(i)l}\|_F^2 - (l-2)\sum_{i=1}^{l-1}\sum_{k=1}^{i-1}\|\widetilde{A}_{\tau(i)\tau(k)}\|_F^2 \label{feps}\\
&\geq & \frac{1}{2} \ \zeta_l^{\ 2(l-2)} \sum_{i=1}^{l-1}\|\widetilde{A}_{\tau(i)l}\|_F^2 - (l-2)S^2(\widetilde{A}) \nonumber \\
&\geq& \ \frac{1}{2} \ \zeta_l^{\ 2(l-2)}(1-\epsilon^2)S^2(A)- (l-2)\eta_{\pi_{l-1},\varrho}S^2(A_{l-1})\nonumber \\
& = &\frac{1}{2} \ \zeta_l^{\ 2(l-2)}(1-\epsilon^2)S^2(A)- (l-2)\eta_{\pi_{l-1},\varrho}\epsilon^2S^2(A)\nonumber \\
& = &f(\epsilon)S^2(A),\nonumber
\end{eqnarray}
where the function $f:[0,1]\rightarrow\mathbb{R}$ is defined by
\[
f(\epsilon)=\frac{1}{2}\zeta_l^{2(l-2)}-\left(\frac{1}{2}\zeta_l^{2(l-2)}+l\eta_{\pi_{l-1},\varrho} - 2\eta_{\pi_{l-1},\varrho}\right)\epsilon^2.
\]
The first derivative of $f$ is not positive since $l\geq2$. Recall that the transformation matrices are $\ubce_{\pi_l}(\varrho)$ and satisfy relations (\ref{uv_1}) and (\ref{uv_1a}). Therefore, we have
\begin{equation}\label{zetal_est}
1\geq \zeta_l=\min_{0\leq k\leq l-3}\left\{\sigma_{\min}(U_{ll}^{(\widetilde{L}+k)})\right\}\geq \varrho
\min_{1\leq i<j\leq l}\gamma_{ij} > \frac{3\sqrt{2}\varrho}{\sqrt{4^{s_l}+26}}>0.
\end{equation}
This implies that the first derivative is negative, i.e., $f$ is decreasing. Its maximum is at $\epsilon=0$, and its zero is
\[
\epsilon_l=\sqrt{\frac{\zeta_l^{2(l-1)}}{\zeta_l^{2(l-1)}+2l\eta_{\pi_{l-1},\varrho}-4\eta_{\pi_{l-1},\varrho}}}.
\]
In relation (\ref{feps}) the left-hand side is nonnegative. Therefore, it is better to replace $f$ with a nonnegative function $f_{+}$, such that
\[
f_{+}(\epsilon)=\left\{
                    \begin{array}{ll}
                      f(\epsilon), & \epsilon\in[0,\epsilon_l\rangle, \\
                      0, & \epsilon\in[\epsilon_l,1].
                    \end{array}
                  \right.
\]
Then $f_+'(\epsilon)\leq0$ for $\epsilon\neq\epsilon_l$ and $f_+(\epsilon)\geq f(\epsilon)$ for $0\leq \epsilon \leq 1$. From relations (\ref{b1mi}) and (\ref{feps}), we have
\begin{align*}
S^2(A^{(L)}) & = S^2(A^{(\widetilde{L})})-\sum_{i=1}^{l-1}\|\widetilde{A}_{\tau(i),l}^{(i)}\|_F^2 \leq S^2(A^{(\widetilde{L})})-f_+(\epsilon)S^2(A)\\
& \leq \big{(}1-\epsilon^2(1-\eta_{\pi_{l-1},\varrho})-f_+(\epsilon)\big{)}S^2(A)
= g(\epsilon)S(A),
\end{align*}
where
\[
g(\epsilon)= \left\{
                    \begin{array}{ll}
                      1-\frac{1}{2}\zeta_l^{2(l-1)}+\epsilon^2(l\eta_{\pi_{l-1},\varrho} - \eta_{\pi_{l-1},\varrho} + \frac{1}{2}\zeta_l^{2(l-1)}-1), & \epsilon\in[0,\epsilon_l\rangle, \vspace{1ex} \\
                      1-\epsilon^2(1-\eta_{\pi_{l-1},\varrho}), & \epsilon\in[\epsilon_l,1].
                    \end{array}
                  \right.
\]
The function $g$ is differentiable on $\langle0,1\rangle \ \backslash \ \{\epsilon_l\}$ and one has
$$g'(\epsilon)= \left\{
                    \begin{array}{ll}
                      \epsilon(2l\eta_{\pi_{l-1},\varrho} - 2\eta_{\pi_{l-1},\varrho}+\zeta_l^{2(l-1)}-2), & \epsilon\in\langle0,\epsilon_l\rangle, \\
                      2\epsilon(\eta_{\pi_{l-1},\varrho}-1), & \epsilon\in\langle\epsilon_l,1\rangle.
                    \end{array}
                  \right.
$$
For $\epsilon\in\langle\epsilon_l,1\rangle$ we have $g'(\epsilon)<0$ since $\eta_{\pi_{l-1},\varrho}-1<0$. For $\epsilon\in\langle0,\epsilon_l\rangle$,  $g'(\epsilon)$ is either positive or negative on the whole interval $\langle0,1\rangle$, depending on $l$. We conclude that $g$ is either decreasing on whole segment $[0,1]$, or increasing on $[0,\epsilon_l\rangle$ and decreasing on $\langle\epsilon_l,1]$. Thus, $g$ attains its maximum either at $\epsilon=0$ or at $\epsilon=\epsilon_l$. Therefore, we have
\[
\eta_{\pi_l,\varrho}=\max\big{\{}g(0),g(\epsilon_l)\big{\}} =
\max \left\{1-\frac{1}{2} \ \zeta_l^{\ 2(l-1)} \ , \ 1-\frac{(1-\eta_{\pi_{l-1},\varrho}) \zeta_l^{\ 2(l-1)}}{\zeta_l^{\ 2(l-1)}+2(l-2)\eta_{\pi_{l-1},\varrho}} \right\}.
\]
From relation (\ref{zetal_est}) we see that $\zeta_l$ is bounded from below by a positive constant that depends on $\pi_l$ and $\varrho$. Therefore, the constant $\eta_{\pi_l,\varrho}$ depends on $\pi_l$ and $\varrho$ and $0\leq \eta_{\pi_l,\varrho}<1$.

It remains to show that there is a bound $\widetilde{\eta}_{s_l,\varrho}$ for $S^2(A^{(L)})/S^2(A)$ depending only on $s_l$ and  $\varrho$ such that $\eta_{\pi_l,\varrho}\leq \widetilde{\eta}_{s_l,\varrho}<1$. It will be derived from $\eta_{\pi_l,\varrho}$.

If $\eta_{\pi_l,\varrho}=g(0)$, then relation (\ref{zetal_est}) implies
\[
\eta_{\pi_l,\varrho} < 1-\frac{1}{2}\left(\frac{3\sqrt{2}\varrho}{\sqrt{4^{s_l}+26}}\right)^{2(l-1)}
< 1-\frac{1}{2}\left(\frac{3\sqrt{2}\varrho}{\sqrt{4^{s_l}+26}}\right)^{2(s_l-1)} \equiv\eta_{s_l,\varrho}'.
\]
If $\eta_{\pi_l,\varrho}=g(\epsilon_l)$, then
\[
\eta_{\pi_l,\varrho}= \frac{2(l-2)+\zeta_l^{2(l-1)}}{\zeta_l^{2(l-1)}+2(l-2)\eta_{\pi_{l-1},\varrho}}\eta_{\pi_{l-1},\varrho},
\]
which implies $\eta_{\pi_{l-1},\varrho}\leq \eta_{\pi_{l},\varrho}$. Therefore
\[
\eta_{\pi_l,\varrho}= 1-\frac{(1-\eta_{\pi_{l-1},\varrho})\,\zeta_l^{2(l-1)}}{\zeta_l^{2(l-1)}+2(l-2)\eta_{\pi_{l-1},\varrho}}
\leq 1-\frac{(1-\eta_{\pi_{l},\varrho})\zeta_l^{2(l-1)}}{\zeta_l^{2(l-1)}+2(l-2)} = \frac{2(l-2)+\eta_{\pi_{l},\varrho}\zeta_l^{2(l-1)}}{2(l-2)+\zeta_l^{2(l-1)}},
\]
which is equivalent to $\eta_{\pi_l,\varrho} \leq (l-2)/(\zeta_l^{2(l-1)}+l-2)$. By (\ref{zetal_est}) we have
\[
\eta_{\pi_l,\varrho} \leq \frac{l-2}{\zeta_l^{2(l-1)}+l-2} < \frac{l-2}{\left(\frac{3\sqrt{2}\varrho}{\sqrt{4^{s_l}+26}}\right)^{2(l-1)}+l-2} <
\frac{s_l-2}{\left(\frac{3\sqrt{2}\varrho}{\sqrt{4^{s_l}+26}}\right)^{2(l-1)}+s_l-2} \equiv\eta_{s_l,\varrho}''.
\]
Finally, set
\[
\widetilde{\eta}_{n,\varrho}=\max\left\{\eta_{s_l,\varrho}',\eta_{s_l,\varrho}''\right\}.
\]
This completes the induction step and the proof of assertion (\ref{b1t}).

\subsection{Proof of Corollary~\ref{tm:c1}}

We follow the proof of Theorem~\ref{tm:b1} and use the notation from therein. Now the blocks become the elements and we can use sharper estimates for the rotation angles. In particular, we can replace the lower bound $\varrho\gamma_{ij}$ from relation (\ref{uv_1a}) with $\sqrt{2}/2$. We have $m=n$, $2\leq l\leq n$ and
\[
\zeta_l=\min_{0\leq k\leq l-1}\{\sigma_{\min}(U_{ll}^{(\widetilde{L}+k)})\}=\min_{0\leq k\leq l-1}{\cos\phi_{\widetilde{L}+k}}\geq\frac{\sqrt{2}}{2}.
\]
Hence, applying that lower bound for $\zeta_l$, using notation $\eta_l$ for $\eta_{\pi_{l}}$, we obtain
\[
\eta_{l}=\max\left\{g(0),g(\epsilon_l)\right\},
\]
where
\[
g(0)=1-2^{-l}, \qquad g(\epsilon_l)=1-\frac{(1-\eta_{l-1})2^{-l}}{2^{-l}+(l-2)\eta_{l-1}}.
\]
This yields the constant $\eta_n$ by replacing $l$ with $n$.

If the whole analysis is performed on the elements, a somewhat larger constant $\eta_n$ can be obtained. In~\cite{B} it has been shown that
\begin{equation*}\label{const1}
\eta_n=\max\Big{\{}1-2^{1-n},1-\frac{2^{2-n}(1-\eta_{n-1})}{2^{2-n}+(n-2)\eta_{n-1}}\Big{\}}.
\end{equation*}

\subsection{Proof of Theorem~\ref{kvazib1}}

The proof is similar to the proof of Theorem~\ref{tm:b1}, so let us follow the same lines. The notation $\pi$, $\pi_l$ and $s_l$ has the same meaning as before. The proof uses the mathematical induction on $l$, $2\leq l\leq m$. Relation (\ref{b1t}) remains the same, except for $L$, which has to take the actual number of steps into account. The same is true for $\widetilde{L}$. Since $\mathcal{O}\in\bar{\mathcal{B}}_c^{(m)}$, we have
\[
L=\frac{l(l-1)}{2}+|\Os_3|+\cdots+|\Os_l|,\quad \widetilde{L}=\frac{(l-1)(l-2)}{2}+|\Os_3|+\cdots+|\Os_{l-1}|.
\]
As earlier, the matrix $\widetilde{A}$ ($A_{l-1}$) is the leading submatrix of $A^{(\widetilde{L})}$ ($A$) of order $s_{l-1}$. Until the end, the proof uses the same lines as the proof of Theorem~\ref{tm:b1}. We only note that the final estimate is first obtained for the matrix $A^{(L-|\mathcal{O}_l|)}$, but since $S(A^{(L)})\leq S(A^{(L-|\mathcal{O}_l|)})$, it automatically holds for $S(A^{(L)})$.

\subsection{Proof of Theorem~\ref{opcitm}}

The proof is similar to the proof of \cite[Theorem~5.1]{har-15}. The only difference is in the fact that here the considered matrices are real and therefore the block Jacobi annihilators and operators are of order $K$ (where $K$ is from relation (\ref{vec1})) and not $2K$ as they are in \cite{har-15}. Also, the iterative process (\ref{opciproces}) uses the congruence transformation, while in \cite{har-15} it uses the equivalence transformation. Finally, here we show how the parameter $\varrho$ is used to avoid the assumption that the matrices $U_k$ from \textbf{A2} have to be \ubce.
Hence, for the completeness of the paper, we will present a somewhat shorter version of the proof, often referring to the proof of \cite[Theorem~5.1]{har-15}. The complete proof can be found in the thesis~\cite{B}.

Using the relation $F_k = U_k+(F_k-U_k)$, $k\geq 0$, and the assumption \textbf{A2} it is easy to transform the process (\ref{opciproces}) into the form
\begin{equation}\label{new_process}
A^{(k+1)}=U_k^TA^{(k)}U_k+E^{(k)}, \quad k\geq 0,
\end{equation}
where the ``perturbation'' matrices $E^{(k)}$ satisfy
\begin{equation}\label{lim_ek}
\lim_{k\rightarrow\infty}\frac{E^{(k)}}{\|A^{(k)}\|_F}=0.
\end{equation}
The matrices $A$, $A^{(k)}$, $U_k$, $E^{(k)}$, $k\geq 0$, carry matrix block-partition defined by $\pi$.
Applying the function $\ve_{\pi}$ to both sides of equation (\ref{new_process}) and using (\ref{lim_ek}) together with condition (i), one obtains
(cf. \cite[Lemma~5.2]{har-15})
\begin{equation}\label{ojantmjedn}
a^{(k+1)}=\mathcal{R}^{(k)}a^{(k)} + g^{(k)}, \quad k \geq 0,
\end{equation}
and
\begin{equation}\label{ojantmg}
\lim_{k\rightarrow\infty}\frac{g^{(k)}}{\|A^{(k)}\|_F}=0.
\end{equation}
Here, $\mathcal{R}^{(k)}$ is the block Jacobi annihilator determined by the pivot submatrix $\widehat{U}_k$ of $U_k$ and the pivot pair $(i(k),j(k))$,
while $g^{(k)}=\ve (H^{(k)})+\ve(E^{(k)})$. The matrix $H^{(k)}$ of order $n$ carries the same partition as $A^{(k)}$, which differs from the zero-matrix only in the pivot submatrix of order $n_i+n_j$ where it equals $\widehat{U}_k^T\widehat{A}_{ij}^{(k)}\widehat{U}_k$.

Recall that one cycle consists of $M$ steps. After the first cycle has been completed, relation (\ref{ojantmjedn}) implies that we can write
$a^{[1]} = \mathcal{J}^{[1]}a^{(0)}+g^{[1]}$, where $a^{[1]}=a^{(M)}$, $\mathcal{J}^{[1]}=\mathcal{R}^{(M-1)}\cdots\mathcal{R}^{(0)}$, and
\[
g^{[1]} = g^{(M-1)}+\sum_{k=0}^{M-2}\mathcal{R}^{(M-1)}\cdots\mathcal{R}^{(k+1)}g^{(k)}.
\]
By Theorem~\ref{tm:bl_jac_ann} we have $\|\mathcal{R}^{(k)}\|_2\leq 1$ for all $k$. Hence,
$\|g^{[1]}\| \leq \|g^{(0)}\|+\cdots+\|g^{(M-1)}\|$.
Similarly, after $s$ cycles we have
\begin{equation}\label{eq2.19}
a^{[s]} = \mathcal{J}^{[s]}a^{[s-1]}+g^{[s]}, \quad s\geq 1,
\end{equation}
with $a^{[s]}=a^{(sM)}$,
$\mathcal{J}^{[s]}=\mathcal{R}^{(sM-1)}\cdots\mathcal{R}^{((s-1)M)}$ and
\begin{equation}\label{gs}
\|g^{[s]}\| \leq \|g^{((s-1)M)}\|+\cdots+\|g^{(sM-1)}\|.
\end{equation}
We will also write $A^{[s]}=A^{(sM)}$, so that $a^{[s]} = \ve (A^{[s]})$, $s\geq 0$.

Using assumption \textbf{A2} and condition (i), it is easy to prove (see \cite[Lemma~5.3]{har-15})
\begin{eqnarray}
\lim_{k\rightarrow\infty}\frac{\|A^{(k+i)}\|_F}{\|A^{(k)}\|_F} &=& 1 \quad \text{for each} \ i\geq 0, \label{ass_a} \\
\lim_{s\rightarrow\infty}\frac{g^{[s]}}{\|A^{[s]}\|_F} &=& 0, \label{ass_b} \\
\lim_{k\rightarrow\infty}\frac{a^{(k)}}{\|A^{(k)}\|_F} &=& 0 \quad \text{iff} \ \lim_{s\rightarrow\infty}\frac{a^{[s]}}{\|A^{[s]}\|_F}=0. \label{ass_c}
\end{eqnarray}
Relation (\ref{ass_a}) is implied by relations (\ref{new_process}) and (\ref{lim_ek}), while relation (\ref{ass_b}) follows directly from (\ref{gs}), (\ref{ass_a}) and (\ref{ojantmg}). Relation (\ref{ass_c}) is implied by (\ref{ojantmg}) and (\ref{ass_a}).
From (\ref{ass_c}) it follows that, to prove $\|a^{(k)}\| / \|A^{(k)}\|_F \rightarrow 0$ as $k\rightarrow\infty$, it is sufficient to show that $\lim_{s\rightarrow\infty}b^{[s]}= 0$ for
\begin{equation}\label{eq2.22}
b^{[s]} = \frac{a^{[s]}}{\|A^{[s]}\|_F}, \quad s\geq 0.
\end{equation}
We transform the iterative process (\ref{eq2.19}) into
\begin{equation}\label{eq2.23}
b^{[s]} = \mathcal{J}^{[s]}b^{[s-1]}+c^{[s]}, \quad s\geq 1,
\end{equation}
where
\begin{equation}\label{eq2.24}
c^{[s]} = \left(\frac{\|A^{[s-1]}\|_F}{\|A^{[s]}\|_F} - 1\right) \mathcal{J}^{[s]} b^{[s-1]} + \frac{g^{[s]}}{\|A^{[s]}\|_F}, \quad s\geq 0.
\end{equation}
By taking the norm of both sides of the equation (\ref{eq2.24}), we obtain
\[
\|c^{[s]}\| \leq \left|\frac{\|A^{[s-1]}\|_F}{\|A^{[s]}\|_F} - 1\right| \|\mathcal{J}^{[s]}\|_2 \|b^{[s-1]}\| + \frac{\|g^{[s]}\|}{\|A^{[s]}\|_F}, \quad s\geq 0.
\]
Relation (\ref{eq2.22}) and Theorem~\ref{tm:bl_jac_ann} imply $\|b^{[s-1]}\| \leq 1$ and $\|\mathcal{J}^{\,[s]}\|_2\leq 1$ for all $s\geq 0$. Hence, from (\ref{ass_a}) and (\ref{ass_b}) we have
\begin{equation}\label{a.30}
\lim_{s\rightarrow\infty}c^{[s]} = 0.
\end{equation}
The proofs of the preceding relations also hold for any quasi-cyclic Jacobi-type process satisfying assumption \textbf{A2} and condition (i) of the theorem. To prove
$$\lim_{s\rightarrow\infty}b^{[s]} = 0,$$
we will additionally use assumptions \textbf{A3} and \textbf{A1}.

Since \textbf{A3} holds, it implies $\sigma>0$. From the definition of $\sigma$, we know that there exists $s_0\geq 1$ such that
\[
\|F_k-U_k\|_2\leq \frac{1}{4}\sigma \quad \text{and} \quad \sigma_{\min}\big{(}F_{ii}^{(k)}\big{)}\geq \frac{3}{4}\sigma, \quad s\geq s_0M.
\]
By the perturbation theorem for the singular values, we have
\begin{align*}
\sigma_{\min}\big{(}U_{ii}^{(k)}\big{)} & = \sigma_{\min}\big{(}F_{ii}^{(k)}- (F_{ii}^{(k)}-U_{ii}^{(k)})\big{)}
\geq \frac{3}{4}\sigma -\|F_{ii}^{(k)}-U_{ii}^{(k)}\|_2 \\
& \geq \frac{3}{4}\sigma -\|F_k-U_k\|_2 \geq \frac{3}{4}\sigma - \frac{1}{4}\sigma = \frac{1}{2}\sigma, \quad k\geq s_0M.
\end{align*}
Set $\varrho =\frac{1}{2}\sigma$. We have proved that $\widehat{U}_k$ belongs to the class $\ubce(\varrho)$ provided that $k\geq s_0 M$.
Then, by Definition~\ref{def:jop}, the block Jacobi operators $\mathcal{J}^{[s]}$, $s\geq s_0$, from (\ref{eq2.23}) are in $\Jl_{\!\!\! \mathcal{O}}^{\ubce(\varrho)}$.

Next, we have to use assumption \textbf{A1}.
Since $\mathcal{O}\in\mathcal{B}_{sg}^{(m)}$, we can presume that the  chain connecting $\mathcal{O}$ to $\mathcal{O''}\in\mathcal{B}_{sp}^{(m)}$ (see Definition~\ref{tm:gss})
is in the canonical form and contains $d$ shift equivalences. Without the loss of generality, we may assume (as in Theorem~\ref{tm:tm3-9}) that $\mathcal{O}\stackrel{\mathsf{p}}{\sim}\mathcal{O'}\stackrel{\mathsf{w}}{\sim}\mathcal{O''} \in\mathcal{B}_{sp}^{(m)}$ with $\mathcal{O'}=\mathcal{O}(\mathsf{q})$, or $\mathcal{O}\stackrel{\mathsf{w}}{\sim}\mathcal{O'}\stackrel{\mathsf{p}}{\sim}\mathcal{O''} \in\mathcal{B}_{sp}^{(m)}$ with $\mathcal{O''}=\mathcal{O'}(\mathsf{q})$, for some permutation $\mathsf{q}$ of the set $\mathcal{S}_m$.

Applying Theorem~\ref{tm:tm3-9}, one concludes that
there are constants $\mu_{\pi_{\mathsf{q}},\varrho}$ and $\widetilde{\mu}_{n,\varrho}$ depending only on $\pi_{\mathsf{q}}$, $\varrho$ and $n$, $\varrho$, respectively, such that
\begin{equation}\label{a.32}
\|\mathcal{J}^{[s+d]}\cdots\mathcal{J}^{[s+1]}\mathcal{J}^{[s]}\|_2 \leq \mu_{\pi_{\mathsf{q}},\varrho},
\quad 0\leq \mu_{\pi_{\mathsf{q}}\!,\varrho}<\widetilde{\mu}_{n,\varrho}<1, \ s\geq s_0.
\end{equation}
By unfolding the recursion (\ref{eq2.23}) $d$ times, similarly as in the proof of \cite[Theorem~5.1]{har-15}, one obtains
$$b^{[s+d]} = \mathcal{J}^{[s+d]}\mathcal{J}^{[s+d-1]}b^{[s-1]}\cdots \mathcal{J}^{[s]} b^{[s-1]} + h^{[s]}, \quad s\geq s_0,$$
where $ h^{[s]}\rightarrow 0$ as $s\rightarrow\infty$. Here, we have used (\ref{a.30}) and Theorem~\ref{tm:bl_jac_ann}. Taking the Euclidean norm, it follows that
$$\|b^{[s+d]}\| \leq \|\mathcal{J}^{[s+d]}\mathcal{J}^{[s+d-1]} b^{[s-1]}\cdots \mathcal{J}^{[s]}\|_2\|b^{[s-1]}\|+\|h^{[s]}\|, \quad s\geq s_0.$$
This inequality together with (\ref{a.32}) implies
\begin{equation}\label{A.34}
\beta_{s+d} \leq \mu_{\pi_{\mathsf{q}},\varrho}\beta_{s-1} + \varepsilon_s, \quad s\geq s_0, \quad \text{with} \ \lim_{s\rightarrow\infty}\varepsilon_s = 0,
\end{equation}
where $0\leq \mu_{\pi_{\mathsf{q}},\varrho} <1$ and $\beta_{s} = \|b^{[s]}\|$, $\varepsilon_s=\|h^{[s]}\|$,
$s\geq s_0$. Set $\alpha_t=\beta_{s_0-1+t(d+1)}$ and $\eta_t=\varepsilon_{s_0+t(d+1)}$, $t\geq 0$. Then from relation (\ref{A.34}) it follows that
\[
\alpha_{t+1}\leq \mu_{\pi_{\mathsf{q}},\varrho} \alpha_{t}+\eta_t, \quad t\geq 0, \quad \text{with} \ \lim_{t\rightarrow\infty}\eta_t =0.
\]
This enables us to apply \cite[Lemma~1]{har-82} to obtain $\lim_{t\rightarrow\infty}\alpha_t=0$, i.e.,
$$\lim_{t\rightarrow\infty} \beta_{s_0-1+t(d+1)} =0.$$
Relations (\ref{eq2.23}) and (\ref{a.30}) imply
$$\beta_{s_0-1+t(d+1)+r} \leq \beta_{s_0-1+t(d+1)} + \vartheta_{r,t}, \quad 0\leq r\leq d, \quad t\geq 0,$$
with $\lim_{t\rightarrow\infty}\vartheta_{r,t} = 0$ for any $0\leq r\leq d$. This proves
$\lim_{s\rightarrow\infty}b^{[s]}= 0$ and, because of (\ref{ass_c}), it also proves
$\lim_{k\rightarrow\infty}\,a^{(k)}/\|A^{(k)}\|_F=0$.

Therefore, it suffices to show that $\displaystyle \lim_{k\rightarrow\infty}\frac{\sum_{l=1}^m S^2(A_{ll}^{(k)})}{\|A^{(k)}\|_F^2}=0$. Let $\varepsilon>0$. Then by condition (i) of the theorem and by relation (\ref{ass_a}), there is an integer $k_{\varepsilon}$ such that (cf. \cite[Theorem~5.1]{har-15})
\begin{equation}\label{kraj1}
\frac{S(\widehat{A}_{ij}^{(k+1)})}{\|A^{(k)}\|_F}\leq \varepsilon, \quad \frac{\|A^{(p)}\|_F}{\|A^{(k)}\|_F}\leq 1+\varepsilon,
\quad k-M\leq p< k, \ k\geq k_{\varepsilon}.
\end{equation}
Here $M$ is the number of steps within one cycle. For given $k\geq k_{\varepsilon}+M$ and $l\in\{1,\ldots,m\}$,  let  $q<k$ denote the last step when $A_{ll}^{(q)}$ was a part of some pivot submatrix. Obviously, $k_{\varepsilon}\leq q < k$. Relation (\ref{kraj1}) implies
\[
\frac{S(A_{ll}^{(k)})}{\|A^{(k)}\|_F}=\frac{S(A_{ll}^{(q+1)})}{\|A^{(k)}\|_F} = \frac{S(A_{ll}^{(q+1)})}{\|A^{(q)}\|_F}\frac{\|A^{(q)}\|_F}{\|A^{(k)}\|_F}
\leq \varepsilon (1+\varepsilon),
\]
for any $1\leq l\leq m$ and any $k\geq k_{\varepsilon}+M$. This proves the theorem.


\begin{thebibliography}{10}

\bibitem{B}
E. Begovi\'{c}:
\emph{Convergence of Block Jacobi Methods}.
Ph.D. thesis, University of Zagreb, 2014.
%
\bibitem{beg+har15a}
E. Begovi\'{c}~Kova\v{c}, V. Hari:
\emph{Jacobi method for symmetric matrices of order $4$ converges for every cyclic pivot strategy}.
arXiv:1701.02387 [math.NA]
%
\bibitem{buj+drm-12}
Z. Bujanovi\'{c}, Z. Drma\v{c}:
\emph{A contribution to the theory and practise of the block Kogbetliantz method for computing the SVD}.
BIT 52 (4) (2012) 827--849.
%
\bibitem{dem-ves-92}
J. Demmel, K. Veseli\'{c}:
\emph{Jacobi's method is more accurate than QR}.
SIAM J. Matrix Anal. Appl. 13 (4) (1992) 1204--1245.
%
\bibitem{drm-07}
Z. Drma\v{c}:
\emph{A global convergence proof of cyclic Jacobi methods with block rotations}.
SIAM J. Matrix Anal. Appl. 31 (3) (2009) 1329--1350.
%
\bibitem{Drmac-Hari-93}
Z. Drma\v{c}, V. Hari:
\emph{On quadratic convergence bounds for the $J$-symmetric {J}acobi method}.
Numer. Math. 64 (1) (1993) 147--180.
%
\bibitem{drm+ves-04a}
Z. Drma\v{c}, K. Veseli\'{c}:
\emph{New fast and accurate Jacobi SVD algorithm I}.
SIAM J. Matrix Anal. Appl. 29 (4) (2008) 1322--1342.
%
\bibitem{drm+ves-04b}
Z. Drma\v{c}, K. Veseli\'{c}:
\emph{New fast and accurate Jacobi SVD algorithm II}.
SIAM J. Matrix Anal. Appl. 29 (4) (2008) 1343--1362.
%
\bibitem{fer-89}
K.~V. Fernando:
\emph{Linear convergence of the row cyclic Jacobi and Kogbetliantz methods}.
Numer. Math. 56 (1) (1989) 73--91.
%
\bibitem{for+hen-60}
G.~E. Forsythe, P. Henrici:
\emph{The cyclic Jacobi method for computing the principal values of a complex matrix}.
Trans. Amer. Math. Soc. 94 (1960) 1--23.
%
\bibitem{han-63}
E.~R. Hansen:
\emph{On cyclic Jacobi methods}.
SIAM J. Appl. Math. 11 (2) (1963) 448--459.
%
\bibitem{har-82}
V. Hari:
\emph{On the global convergence of Eberlein method for real matrices}.
Numer. Math. 39 (3) (1982) 361--369.
%
\bibitem{har-84}
V. Hari:
\emph{On Cyclic Jacobi Methods for the Positive Definite Generalized Eigenvalue Problem}.
Ph.D. thesis, University of Hagen, 1984.
%
\bibitem{har-86}
V. Hari:
\emph{On the convergence of cyclic Jacobi-like processes}.
Linear Algebra Appl. 81 (1986) 105--127.
%
\bibitem{har-91}
V. Hari:
\emph{On sharp quadratic convergence bounds for the serial Jacobi methods}.
Numer. Math. 60 (1) (1991) 375--406.
%
\bibitem{har-07}
V. Hari:
\emph{Convergence of a block-oriented quasi-cyclic {J}acobi method}.
SIAM J. Matrix Anal. Appl. 29 (2) (2007) 349--369.
%
\bibitem{har-09}
V. Hari:
\emph{On block Jacobi annihilators}.
Proceedings of Algoritmy, Publishing House of STU, Slovak University of Technology in Bratislava (2009) 429--439.
%
\bibitem{har-15}
V. Hari:
\emph{Convergence to diagonal form of block Jacobi-type methods}.
Numer. Math. 129 (3) (2015) 449--481.
%
\bibitem{har+sin-10}
V. Hari, S. Singer, S. Singer:
\emph{Block-oriented $J$-Jacobi methods for Hermitian matrices}.
Linear Algebra Appl. 433 (8-10) (2010) 1491--1512.
%
\bibitem{har+sin-11}
V. Hari, S. Singer, S. Singer:
\emph{Full block $J$-Jacobi method for Hermitian matrices}.
Linear Algebra Appl. 444 (2014) 1--27.
%
\bibitem{har+ves-87}
V. Hari, K. Veseli\'{c}:
\emph{On Jacobi methods for singular value decompositions}.
SIAM J. Sci. and Stat. Comput. 8 (5) (1987) 741--754.
%
\bibitem{hen+zim-68}
P. Henrici, K. Zimmermann:
\emph{An estimate for the norms of certain cyclic Jacobi operators}.
Linear Algebra Appl. 1 (4) (1968) 489--501.
%
\bibitem{jac-1846}
C.~G.~J. Jacobi:
\emph{\"{U}ber ein leichtes Verfahren die in der Theorie der S\"{a}cularst\"{o}rungen vorkommenden Gleichungen numerisch aufzul\"{o}sen}.
Crelle's Journal 30 (1846) 51–-95.
%
\bibitem{luk+par-89}
F.~T. Luk, H. Park:
\emph{On parallel Jacobi orderings}.
SIAM J. Sci. and Stat. Comput. 10 (1) (1989) 18--26.
%
\bibitem{mas-90}
W.~F. Mascarenhas:
\emph{On the convergence of the Jacobi method for arbitrary orderings}.
SIAM J. Matrix Anal. Appl. 16 (4) (1995) 1197–-1209.
%
\bibitem{mat-08}
J. Mateja\v{s}:
\emph{Accuracy of the Jacobi method on scaled diagonally dominant symmetric matrices}.
SIAM J. Matrix Anal. Appl. 31 (1) (2009) 133--153.
%
\bibitem{math-95}
R. Mathias:
\emph{Accurate eigensystem computations by Jacobi methods}.
SIAM J. Matrix Anal. Appl. 16 (3) (1995) 977–-1003.
%
\bibitem{naz-75}
L. Nazareth:
\emph{On the convergence of the cyclic Jacobi methods}.
Linear Algebra Appl. 12 (2) (1975) 151--164.
%
\bibitem{nov+sin-15}
V. Novakovi\'{c}, S. Singer, S. Singer:
\emph{Blocking and parallelization of the Hari–-Zimmermann variant of the Falk–-Langemeyer algorithm for the generalized SVD}.
Parallel Comput. 49 (2015) 136--152.
%
\bibitem{rhe+har-93}
H.~N. Rhee, V. Hari:
\emph{On the global and cubic convergence of a quasi-cyclic Jacobi method}.
Numer. Math. 66 (1) (1993) 97--122.
%
\bibitem{sam-71}
A.~H. Sameh:
\emph{On Jacobi and Jacobi-like algorithms for parallel computer}.
Math. Comp. 25 (1971) 579--590.
%
\bibitem{shr+sch-89}
G. Shroff, R. Schreiber:
\emph{On the convergence of the cyclic Jacobi method for parallel block orderings}.
SIAM J. Matrix Anal. Appl. 10 (3) (1989) 326--346.
%
\bibitem{sla-02}
I. Slapni\v{c}ar:
\emph{Highly accurate symmetric eigenvalue decomposition and hyperbolic SVD}.
Linear Algebra Appl. 358 (1-3) (2003) 387--424.
%
\bibitem{ves+har-89}
K. Veseli\'{c}, V. Hari:
\emph{A note on a one-sided Jacobi algorithm}.
Numer. Math. 56 (6) (1989) 627--633.
%
\bibitem{wil-62}
J.~H. Wilkinson:
\emph{Note on the quadratic convergence of the cyclic Jacobi process}.
Numer. Math. 4 (1) (1962) 296--300.

\end{thebibliography}
\end{document}